# FIBRATIONS AND YONEDA'S LEMMA IN AN $\infty$-COSMOS

EMILY RIEHL AND DOMINIC VERITY

Abstract. We use the terms $\infty$-*categories* and $\infty$-*functors* to mean the objects and morphisms in an $\infty$-*cosmos*: a simplicially enriched category satisfying a few axioms, reminiscent of an enriched category of fibrant objects. Quasi-categories, Segal categories, complete Segal spaces, marked simplicial sets, iterated complete Segal spaces, $\theta_n$-spaces, and fibered versions of each of these are all $\infty$-categories in this sense. Previous work in this series shows that the basic category theory of $\infty$-categories and $\infty$-functors can be developed only in reference to the axioms of an $\infty$-cosmos; indeed, most of the work is internal to the *homotopy 2-category*, a strict 2-category of $\infty$-categories, $\infty$-functors, and natural transformations. In the $\infty$-cosmos of quasi-categories, we recapture precisely the same category theory developed by Joyal and Lurie, although our definitions are 2-categorical in natural, making no use of the combinatorial details that differentiate each model.

In this paper, we introduce cartesian fibrations, a certain class of $\infty$-functors, and their groupoidal variants. Cartesian fibrations form a cornerstone in the abstract treatment of "category-like" structures a la Street and play an important role in Lurie's work on quasi-categories. After setting up their basic theory, we state and prove the Yoneda lemma, which has the form of an equivalence between the quasi-category of maps out of a representable fibration and the quasi-category underlying the fiber over its representing element. A companion paper will apply these results to establish a calculus of *modules* between $\infty$-categories, which will be used to define and study pointwise Kan extensions along $\infty$-functors.

## Contents

## 1. Introduction

$(\infty, 1)$-*categories* are infinite-dimensional categories with non-invertible morphisms only in dimension one. Equivalently, $(\infty, 1)$-categories are categories weakly enriched over $\infty$-groupoids, i.e., topological spaces. These schematic definitions are realized by a number of concrete *models* of $(\infty, 1)$-categories. Important independent work of Töen and of Barwick and Schommer-Pries proves that all models of $(\infty, 1)$-categories "have the same homotopy theory," in the sense of being connected by a zig-zag of Quillen equivalences of model categories [26] or having equivalent quasi-categories [1]. Inspired by this result, the dream is to be able to work with $(\infty, 1)$-categories "model independently," which begs the question: can the category theory, and not just the homotopy theory, of $(\infty, 1)$-categories be developed model independently?

This paper describes one possible direction to take in pursuit of that goal. We introduce the notion of an $\infty$-*cosmos*, a simplicially enriched category whose objects we call $\infty$-*categories* and whose morphisms we call $\infty$-*functors* or simply *functors*. A quotient defines a strict 2-category which we call the *homotopy 2-category* of the $\infty$-cosmos, whose objects are again $\infty$-categories, whose morphisms are functors between them, and whose 2-cells are natural transformations of a suitable variety. The homotopy 2-category should be thought of as a categorification of the usual notion of homotopy category spanned by the fibrant-cofibrant objects in a model category that is analogous to the 2-category of ordinary categories, functors, and natural transformations — which indeed is the homotopy 2-category of a suitable $\infty$-cosmos.

Previous work [18, 21, 19] shows that a large portion of the category theory of quasi-categories—one model of $(\infty, 1)$-categories that has been studied extensively by Joyal, Lurie, and others—can be developed in the homotopy 2-category of the $\infty$-cosmos of quasi-categories. Indeed, nearly all of the results in these papers, which develop the basic theory of adjunctions, limits and colimits, and monadicity, apply in the homotopy 2-category of any $\infty$-cosmos. In particular, complete Segal spaces, Segal categories, and marked simplicial sets all have their own $\infty$-cosmoi; not coincidentally, these are the models of $(\infty, 1)$-categories whose model categories are the best behaved. Thus each of these varieties of $(\infty, 1)$-categories are examples of $\infty$-categories, in our sense. The axioms imply that the 2-categorical notion of equivalence, interpreted in the homotopy 2-category, precisely coincides with the model categorical notion of weak equivalence. Thus the category theory developed here is appropriately "homotopical," i.e., weak equivalence invariant.

Unlike the work of Töen and Barwick-Schommer-Pries, an $\infty$-cosmos is not meant to axiomatize a "simplicially enriched category of $(\infty, 1)$-categories." For instance, slices of an $\infty$-cosmos again define an $\infty$-cosmos. Indeed, $\theta_n$-spaces and iterated complete Segal spaces, two of the most prominent models of $(\infty, n)$-categories, also define $\infty$-cosmoi. Thus, our work begins to develop the basic category theory of $(\infty, n)$-categories as well.

There is a good notion of functor between $\infty$-cosmoi that preserves all of the structure specified by the axiomatization. Examples include "underlying $(\infty, 1)$-category" functors from the cosmoi for $\theta_n$-spaces or iterated complete Segal spaces to the $\infty$-cosmos for quasi-categories. There is also a functor from the $\infty$-cosmos for strict 1-categories (whose homotopy 2-category is the usual 2-category of categories) to the $\infty$-cosmos of quasi-categories or of complete Segal spaces, and also a functor from the $\infty$-cosmos for Kan complexes to the $\infty$-cosmos for quasi-categories. A certain special class of functors of $\infty$-cosmoi, coming from enriched right Quillen equivalences of model categories, both preserve the structures in the $\infty$-cosmoi and reflect equivalences. These functors give a strong meaning to the sense in which the basic category theory of $(\infty, 1)$-categories developed in this framework is "model independent": basic categorical notions are both preserved and reflected by the functors between the $\infty$-cosmoi of quasi-categories, complete Segal spaces, Segal categories, and marked simplicial sets. Furthermore, a theorem of Low implies that the induced 2-functors between their homotopy 2-categories define bicategorical equivalences [11].

In §2 we define $\infty$-cosmoi and functors between them and demonstrate that all of the examples listed above can be realized as the underlying 1-category of a suitable $\infty$-cosmos. In §3, we define the homotopy 2-category of an $\infty$-cosmos and explore its relevant 2-categorical structure. In fact, the majority of the results in this paper can be stated and proven in an *abstract homotopy 2-category*, which we define to be a (strict) 2-category equipped with *comma objects* and *iso-comma objects* of a suitably weak variety.

The second half of this paper continues the project of developing the basic category theory of $\infty$-categories — including, for the reasons just explained, models of $(\infty, n)$-categories and their sliced variants — by introducing a suitably model independent notion of cartesian fibration, an important cornerstone in an abstract theory of "category-like" structures [22]. Functors valued in $\infty$-categories are most efficiently encoded as *cartesian* or *cocartesian fibrations*, the "co" signaling that the functor so-encoded is covariant. We also study a "groupoidal" variant of these notions whose fibers are $\infty$-categories that are "representably groupoidal." For the models of $(\infty, 1)$-categories, these groupoidal fibers are the $\infty$-groupoids in the sense appropriate to each model. In general, the underlying quasi-category of a representably groupoidal $\infty$-category is a Kan complex.

Cartesian fibrations, as introduced by Lurie [12], play an important role in the theory of quasi-categories. In a model-independent context, we are not able to make use of his definition which refers to the fact that quasi-categories are simplicial sets. Instead, we present a new definition of cartesian fibration, defined in any homotopy 2-category, that when interpreted in the homotopy 2-category of quasi-categories coincides precisely with Lurie's notion. In the special case of quasi-categories, the groupoidal cartesian fibrations are precisely the *right fibrations* introduced by Joyal, while the groupoidal cocartesian fibrations are the dual *left fibrations*.

In §4 we introduce cartesian and cocartesian fibrations. The main theorem presents three equivalent characterizations of cartesian fibration in an abstract homotopy 2-category making use of the comma constructions mentioned above.

4.1.10. **Theorem.** *If $p\colon E \twoheadrightarrow B$ is an isofibration, then the following are equivalent:*

*(i) $p$ is a cartesian fibration.*
*(ii) The induced functor $i\colon E \to B \downarrow p$ admits a right adjoint which is fibred over $B$:*

$$\begin{array}{ccc} B\downarrow p & \underset{r}{\overset{i}{\rightleftarrows}}\ \bot & E \\ & {\scriptstyle p_0}\searrow \quad \swarrow{\scriptstyle p} & \\ & B & \end{array}$$

*(iii) The induced functor $k\colon E^2 \to B \downarrow p$ admits a right adjoint right inverse:*

$$B\downarrow p \underset{\bar{r}}{\overset{k}{\rightleftarrows}}\ \bot\ E^2$$

Groupoidal cartesian fibrations, defined to be cartesian fibrations with groupoidal fibers, admit a similar characterization.

4.2.7. **Proposition.** *An isofibration $p\colon E \twoheadrightarrow B$ is a groupoidal cartesian fibration if and only if the functor $k\colon E^2 \to B \downarrow p$ is an equivalence.*

An important corollary of these results is that cartesian and groupoidal cartesian fibrations are representably-defined notions. Indeed, they are preserved by any functors of $\infty$-cosmoi, up to possibly replacing the map in the image by an equivalent isofibration.

So as to not interrupt the overall narrative flow, we start with a sketched proof of Theorem 4.1.10, deferring full details to the appendix. In §5, we prove an analogous result characterizing cartesian functors between cartesian fibrations. We also prove that cartesian fibrations are stable under pullback, a somewhat delicate result because the sort of pullbacks that are present in the homotopy 2-category are somewhat weak.

Finally, in §6, we formulate and prove the Yoneda lemma in a form inspired by Street [23].

6.0.1. **Theorem** (Yoneda lemma)**.** *Given any cartesian fibration $p\colon E \twoheadrightarrow B$ and any point $b\colon 1 \to B$, restriction along the terminal object $t\colon 1 \to B \downarrow b$ induces an equivalence of quasi-categories*

$$\mathrm{map}_B^{\mathrm{cart}}(p_0\colon B\downarrow b \twoheadrightarrow B, p\colon E \twoheadrightarrow B) \simeq \mathrm{map}_B(b\colon 1 \to B, p\colon E \twoheadrightarrow B).$$

Here $p_0 \colon B \downarrow b \twoheadrightarrow B$ is the groupoidal cartesian fibration represented by the point $b$. The quasi-category $\mathrm{map}_B(b, p)$ can be thought of as the underlying quasi-category of the fiber of $p$ over $b$.

A companion paper will use the theory developed here to establish a calculus of two-sided groupoidal fibrations, which are meant to encode $\infty$-groupoid-valued bifunctors—prototypically the "$\infty$-category of arrows"—that are covariant in one variable and contravariant in another. Following the Australian school, we call these *modules*; synonyms include *profunctor*, *correspondence*, or *distributor*. Using modules, we define and develop the basic theory of pointwise Kan extensions along functors between $\infty$-categories.

References to [18, 21, 19] will have the form I.x.x.x, II.x.x.x, and III.x.x.x, respectively. We refer the reader to §I.2 for an account of the quasi-categorical notational conventions to be used, which are standard.

1.1. **Acknowledgments.** This material is based upon work supported by the National Science Foundation under Award No. DMS-1103790 and by the Australian Research Council under Discovery grant number DP130101969. A substantial portion of this work was completed while the second-named author was visiting Harvard University, during which time he was partially supported by an NSF grant DMS-0906194 and a DARPA grant HR0011-10-1-0054-DOD35CAP held by Mike Hopkins. We are particularly grateful for his support. Some of the writing took place while the first-named author was in residence at the Hausdorff Research Institute for Mathematics, supported by an AMS-Simons Travel Grant. Peter May suggested the name "$\infty$-cosmos," a substantial improvement upon previous informal terminology.

## 2. $\infty$-COSMOI

An $\infty$-*cosmos*, introduced below, is an axiomatization of the basic properties of the simplicial category of quasi-categories and its simplicial slices, which together present some form of categorified derivator. A particular quotient of an $\infty$-cosmos defines a strict 2-category, which we refer to as its *homotopy 2-category*, with certain properties. An $\infty$-cosmos and its homotopy 2-category each have the same underlying 1-category, whose objects we call $\infty$-categories and whose morphisms we call $\infty$-functors. For example, quasi-categories, Segal categories, complete Segal spaces, general *Rezk objects* such as $n$-fold complete Segal spaces, and fibered versions of each of these objects all define the $\infty$-categories in suitable $\infty$-cosmoi.

The remainder of this paper, and indeed the work contained in the previous papers in this series, relies only upon these axioms. In this way, we can develop the basic theory of $\infty$-categories and $\infty$-functors uniformly across various models, i.e., the particular examples mentioned above. In fact, much of the basic theory can be developed in an *abstract homotopy 2-category*, a notion which will be defined in §3.

In §2.1, we introduce $\infty$-cosmoi and develop a bit of the general theory, while in §2.2 we present a number of examples.

2.1. **$\infty$-cosmoi and functors.** The prototypical example of an $\infty$-cosmos is $\underline{\mathrm{qCat}}$, the simplicially enriched category of quasi-categories. The underlying 1-category is comprised of the quasi-categories and functors between them, which are simply maps of simplicial sets.

2.1.1. **Definition** ($\infty$-cosmos)**.** Suppose that $\mathcal{K}$ is a simplicially enriched category, equipped with two specified classes of its 0-arrows $\mathcal{W}$ called its class of *weak equivalences* and $\mathcal{F}$ called its class of *isofibrations*. As usual, we shall assume that $\mathcal{W}$ contains all isomorphisms and satisfies the *2-of-6 property* [6] and that $\mathcal{F}$ also contains all isomorphisms and is closed under composition. A functor of $\mathcal{K}$ which is both an isofibration and an equivalence is called a *trivial fibration*. As usual we shall use the decorated arrow symbols $\twoheadrightarrow$, $\xrightarrow{\sim}$, and $\overset{\sim}{\twoheadrightarrow}$ to distinguish arrows which are isofibrations, equivalences, and trivial fibrations respectively. We shall use the notation $\mathrm{map}(A, B)$ to denote the hom-space between the objects $A$ and $B$ of $\mathcal{K}$.

We say that $\mathcal{K}$ together with its classes of weak equivalences and isofibrations is an $\infty$-*cosmos* if and only if it satisfies the following axioms:

(a) (completeness) the category $\mathcal{K}$ possesses a terminal object 1, cotensors $U \pitchfork A$ of all objects $A$ by all finitely presented simplicial sets $U$ (those that only have a finite set of non-degenerate simplices), and pullbacks of isofibrations along any functor;
(b) (fibrant objects) all of the functors $!\colon A \to 1$ with codomain 1 are isofibrations;
(c) (pullback stability) the classes of isofibrations and trivial fibrations are both stable under pullback along all functors;
(d) (Leibniz stability) if $p\colon E \twoheadrightarrow B$ is an isofibration in $\mathcal{K}$ and $i\colon U \hookrightarrow V$ is an inclusion of finitely presented simplicial sets then the Leibniz cotensor $i \,\widehat{\pitchfork}\, p\colon V \pitchfork E \to U \pitchfork E \times_{U \pitchfork B} V \pitchfork B$ is an isofibration and it is a trivial fibration when $p$ is a trivial fibration in $\mathcal{K}$ or $i$ is trivial cofibration in the Joyal model structure on $\underline{\mathrm{sSet}}$; and
(e) (cofibrant replacement) if $A$ is an arbitrary object of $\mathcal{K}$ then there exists a trivial fibration $r_A\colon \tilde{A} \overset{\sim}{\twoheadrightarrow} A$ for which the object $\tilde{A}$ is *cofibrant* in the sense that it enjoys the left lifting property with respect to all trivial fibrations in $\mathcal{K}$.

2.1.2. *Observation.* It is a routine consequence of axioms 2.1.1(a) and 2.1.1(b) that the binary product $A \times B$ of any pair of objects exists in an $\infty$-cosmos $\mathcal{K}$. What is more, axiom 2.1.1(c) implies that each of the projections $\pi_0\colon A \times B \twoheadrightarrow B$ and $\pi_1\colon A \times B \twoheadrightarrow A$ is an isofibration. It follows that if we are given an isofibration $p\colon E \twoheadrightarrow A \times B$ then both of its components $p_0 := \pi_0 p\colon E \to B$ and $p_1 := \pi_1 p\colon E \to A$ are also isofibrations. However, it should be noted that the converse implication cannot be assumed to hold.

2.1.3. *Observation.* We also ought to justify the existence of the pullback used to define the codomain of the Leibniz product $i \,\widehat{\pitchfork}\, p$ when stating condition 2.1.1(d). To this end, observe that in the case of an inclusion $\emptyset \hookrightarrow V$, the pullback in question is a trivial one along the identity on 1, and $i \pitchfork p$ is isomorphic to the cotensor $V \pitchfork p\colon V \pitchfork E \to V \pitchfork B$. So in this particular situation condition 2.1.1(d) is well defined and it simply postulates that $V \pitchfork p$ is an isofibration for any isofibration $p$ and any simplicial set $V$. Now the existence of all of the other pullbacks mentioned in condition 2.1.1(d) follows from the completeness condition 2.1.1(a), simply because they are all pullbacks along isofibrations of the form $V \pitchfork p$ for some isofibration $p$.

2.1.4. **Example** (the $\infty$-cosmos of quasi-categories)**.** Following the conventions established in [18], we shall let $\underline{\mathrm{qCat}}$ denote the simplicially enriched category whose objects are the quasi-categories and whose hom-space $\mathrm{map}(A, B)$ is the simplicial function space $B^A$. It is a standard result then that each hom-space of $\underline{\mathrm{qCat}}$ is itself a quasi-category. Now we may define a (weak) equivalence $w\colon A \to B$ of quasi-categories to be a functor for which there is some functor $w'\colon B \to A$ and isomorphisms $ww' \cong \mathrm{id}_B$ in the quasi-category $B^B$ and $w'w \cong \mathrm{id}_A$ in the quasi-category $A^A$. We also take the isofibrations of $\underline{\mathrm{qCat}}$ to be those functors $p\colon E \to B$ of quasi-categories which have the right lifting property with respect to all inner horns $\Lambda^{n,k} \hookrightarrow \Delta^n$ $(0 < k < n)$ and either of the inclusions $\Delta^0 \hookrightarrow J$. Here $J$ denotes the (nerve of) the generic isomorphism category $\{0 \cong 1\}$. These choices make $\underline{\mathrm{qCat}}$ into an $\infty$-cosmos in which every object is cofibrant. Specifically the axioms laid out above follow from standard results used in the construction of Joyal's model structure on simplicial sets (see [8], [27], or I.2.2.5-I.2.2.9 for example).

Other examples of $\infty$-cosmoi will follow shortly after we first develop a bit of the general theory.

2.1.5. *Recall.* Suppose that $\mathbb{I}$ is the simplicial subset of the isomorphism category $J = \{0 \cong 1\}$ generated by its non-degenerate 3-simplex determined by the vertex sequence $\{1, 0, 1, 0\}$. The following facts, which the reader may glean from [8] or [27], are entirely standard foundational results in the theory of quasi-categories:

(a) An arrow $f$ in a quasi-category $A$ is an isomorphism, in the sense that it is carried to an isomorphism in the homotopy category $h(A)$ under the quotient functor $A \to h(A)$, if and only if there exists a simplicial map $e\colon \mathbb{I} \to A$ which carries the 1-simplex of $\mathbb{I}$ with vertex sequence $\{0, 1\}$ to $f$;
(b) The inclusions $i_0, i_1\colon \Delta^0 \hookrightarrow \mathbb{I}$ and the unique map $!\colon \mathbb{I} \to \Delta^0$ are all weak equivalences in the Joyal model structure;
(c) The inclusion $\mathbb{I} \hookrightarrow J$ may be presented as a countable composite of pushouts of inner horn inclusions (one at each dimension $n > 3$); and
(d) A functor $p\colon E \to B$ between quasi-categories is an isofibration if and only if it possesses the right lifting property with respect to all inner horn inclusions $\Lambda^{n,k} \hookrightarrow \Delta^n$ and either of the inclusions $\Delta^0 \hookrightarrow \mathbb{I}$.

2.1.6. **Lemma** (Brown's factorisation lemma)**.** *Fix an $\infty$-cosmos $\mathcal{K}$.*

*(i) Any object $B$ has a path object defined by cotensoring with the simplicial set $\mathbb{I}$.*

$$\begin{array}{ccc} B & \xrightarrow{\;i\;} & \mathbb{I} \pitchfork B \\ & {\scriptstyle (\mathrm{id}_B, \mathrm{id}_B)} \searrow & \downarrow\!\!\downarrow {\scriptstyle (p_1, p_0)} \\ & & B \times B \end{array}$$

*(ii) Any functor* $f \colon A \to B$ *may be factored as*

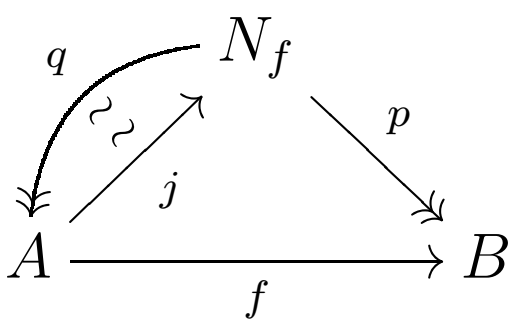

*where* $p \colon N_f \twoheadrightarrow B$ *is an isofibration and* $j \colon A \xrightarrow{\sim} N_f$ *is right inverse to a trivial fibration* $q \colon N_f \overset{\sim}{\twoheadrightarrow} A$.

It is traditional to refer to the object $N_f$ as the *mapping path space* of $f$.

*Proof.* Apply the contravariant functor $- \pitchfork B$ to the inclusion $\Delta^0 + \Delta^0 \hookrightarrow \mathbb{I}$ to construct a functor $(p_1, p_0) \colon \mathbb{I} \pitchfork B \to B \times B$ which, by application of axiom 2.1.1(d), is an isofibration. Furthermore, we may also express the individual functors $p_0, p_1 \colon \mathbb{I} \pitchfork B \to B$ as those obtained by applying $- \pitchfork B$ to the inclusions $i_0, i_1 \colon \Delta^0 \hookrightarrow \mathbb{I}$. Recollection 2.1.5(b) tells us that these latter inclusions are trivial cofibrations in the Joyal model structure, so applying the Leibniz stability axiom again we see that $p_0$ and $p_1$ are trivial fibrations. The diagonal map $B \to B \times B$ factors as a composite of $i :=! \pitchfork B \colon B \to \mathbb{I} \pitchfork B$ and $(p_1, p_0) \colon \mathbb{I} \pitchfork B \to B \times B$, so it follows that $i$ is also a weak equivalence, since it is right inverse to each of the trivial fibrations $p_0$ and $p_1$.

The details of the construction of the mapping path space factorization using a path object are standard; see, e.g., [3, p. 421]. □

When combined with the other $\infty$-cosmos axioms, this demonstrates that (the underlying category of) an $\infty$-cosmos is a *category of fibrant objects* in the sense introduced by Brown [3].

2.1.7. *Observation.* In the case where $A$ and $B$ are cofibrant we may pick the mapping path space $N_f$ in this factorisation so that it too is cofibrant. To do this form the factorisation of Lemma 2.1.6, take the cofibrant replacement $\tilde{N}_f$ of the mapping path space using axiom 2.1.1(e), and lift the map $j$ along the trivial fibration $r \colon \tilde{N}_f \overset{\sim}{\twoheadrightarrow} N_f$ using the assumption that $A$ is cofibrant.

2.1.8. **Lemma.** *Suppose that* $X$ *is a cofibrant object and that* $p \colon E \twoheadrightarrow B$ *is an isofibration (respectively trivial fibration) in an* $\infty$*-cosmos* $\mathcal{K}$*. Then the hom-space* $\mathrm{map}(X, B)$ *is a quasi-category and* $\mathrm{map}(X, p) \colon \mathrm{map}(X, E) \to \mathrm{map}(X, B)$ *is an isofibration (respectively trivial fibration) of quasi-categories.*

*Proof.* A standard duality argument tells us that $\mathrm{map}(X, p) \colon \mathrm{map}(X, E) \to \mathrm{map}(X, B)$ has the right lifting property with respect to some inclusion $i \colon U \hookrightarrow V$ of finitely presented simplicial sets if and only if the Leibniz cotensor $i \mathbin{\widehat{\pitchfork}} p \colon V \pitchfork E \to U \pitchfork E \times_{U \pitchfork B} V \pitchfork B$ has the right lifting property with respect to $X$. Now axiom 2.1.1(d) also tells us that the latter Leibniz cotensor is a trivial fibration whenever $p$ is a trivial fibration in $\mathcal{K}$ or $i$ is a trivial cofibration in the Joyal model structure. So under either of those conditions, we find that $i \mathbin{\widehat{\pitchfork}} p$ has the right lifting property with respect to the cofibrant object $X$, and thus that $\mathrm{map}(X, p)$ has the right lifting property with respect to $i$.

Now we may prove that $\mathrm{map}(X,B)$ is a quasi-category by applying the result of the last paragraph to the isofibration $!\colon B \twoheadrightarrow 1$ of axiom 2.1.1(b) and each of the inner horn inclusions $\Lambda^{n,k} \hookrightarrow \Delta^n$ in turn. Then we may prove that $\mathrm{map}(X,p)\colon \mathrm{map}(X,E) \to \mathrm{map}(X,B)$ is an isofibration of quasi-categories by applying the result of the last paragraph to the isofibration $p\colon E \twoheadrightarrow B$ and each of the inner horn inclusions $\Lambda^{n,k} \hookrightarrow \Delta^n$ and the inclusion $1 \hookrightarrow \mathbb{I}$ of Recollection 2.1.5 in turn.

Finally if $p\colon E \xrightarrow{\sim}\!\!\!\!\twoheadrightarrow B$ is a trivial fibration then we may apply the same argument with respect to each of the boundary inclusions $\partial\Delta^n \hookrightarrow \Delta^n$ to show that $\mathrm{map}(X,p)$ is a trivial fibration of quasi-categories. □

2.1.9. **Definition.** A *functor of $\infty$-cosmoi* $F\colon \mathcal{K} \to \mathcal{L}$ is a simplicial functor that carries isofibrations (resp. trivial fibrations) of $\mathcal{K}$ to isofibrations (resp. trivial fibrations) in $\mathcal{L}$ and preserves the limits listed in 2.1.1(a).

2.1.10. **Proposition.** *For any cofibrant object $X$ in an $\infty$-cosmos $\mathcal{K}$, the simplicial representable* $\mathrm{map}(X,-)$ *is a functor of $\infty$-cosmoi from $\mathcal{K}$ to* $\underline{\mathrm{qCat}}$.

*Proof.* Under the assumption that $X$ is cofibrant, Lemma 2.1.8 tells us that $\mathrm{map}(X,-)$ carries each object of $\mathcal{K}$ to a quasi-category, so it provides us with a simplicial functor $\mathrm{map}(X,-)\colon \mathcal{K} \to \underline{\mathrm{qCat}}$. Furthermore, that lemma also tells us that $\mathrm{map}(X,-)$ carries isofibrations (resp. trivial fibrations) in $\mathcal{K}$ to isofibrations (resp. trivial fibrations) in $\underline{\mathrm{qCat}}$. Its preservation of the various limits possessed by $\mathcal{K}$ is simply the familiar result that (enriched) covariant representables preserve (weighted) limits. □

2.1.11. **Example** (sliced $\infty$-cosmoi)**.** Suppose that $B$ is a fixed object in an $\infty$-cosmos $\mathcal{K}$. Then we shall let $\mathcal{K}/B$ denote the *sliced $\infty$-cosmos over $B$*: the full simplicial subcategory of the usual simplicial slice category whose objects are those functors $p\colon E \twoheadrightarrow B$ that are isofibrations. Explicitly, if $p\colon E \twoheadrightarrow B$ and $q\colon F \twoheadrightarrow B$ are two objects of this slice then the simplicial hom-space $\mathrm{map}_B(p,q)$ between them is formed by taking the pullback

$$\begin{array}{ccc} \mathrm{map}_B(p,q) & \longrightarrow & \mathrm{map}(E,F) \\ \downarrow & \lrcorner & \downarrow{\scriptstyle \mathrm{map}(E,q)} \\ \Delta^0 & \xrightarrow[p]{} & \mathrm{map}(E,B) \end{array} \tag{2.1.12}$$

in simplicial sets. In order to equip $\mathcal{K}/B$ with the rest of the structure of an $\infty$-cosmos, a functor in there is taken to be an isofibration or a weak equivalence if and only if its underlying functor in $\mathcal{K}$ is such. With these definitions in place it is now easily verified that $p\colon E \twoheadrightarrow B$ is a cofibrant object in the slice $\mathcal{K}/B$ if and only if $E$ is a cofibrant object in $\mathcal{K}$.

As we might expect this slice $\mathcal{K}/B$ has the identity map $\mathrm{id}_B\colon B \to B$ as its terminal object and its pullbacks are constructed by taking pullbacks of underlying diagrams in $\mathcal{K}$. Combining these observations with the fact that isofibrations and trivial fibrations in the slice $\mathcal{K}/B$ are also defined in terms of the corresponding property of their underlying functors in $\mathcal{K}$ we immediately verify axioms 2.1.1(b), (c), and (e). All that remains of axiom 2.1.1(a) is to construct the cotensor of an object $p\colon E \twoheadrightarrow B$ of $\mathcal{K}/B$ by a finitely presented simplicial set $U$, but this is just the left-hand vertical arrow in the following pullback in $\mathcal{K}$:

$$\begin{array}{ccc} U \pitchfork_p E & \longrightarrow & U \pitchfork E \\ \downarrow & \lrcorner & \downarrow {\scriptstyle U \pitchfork p} \\ B & \xrightarrow[\Delta]{} & U \pitchfork B \end{array}$$

Here the arrow $\Delta$ appearing along the bottom is the adjoint transpose of the constant map $U \to \Delta^0 \xrightarrow{\mathrm{id}_B} \mathrm{map}(B,B)$ at the identity for $B$. Finally, with these observations in place it is only a matter of a routine pullback calculation to show that the Leibniz stability axiom 2.1.1(d) also holds for the slice.

2.1.13. **Proposition** (pulling back between sliced $\infty$-cosmoi)**.** *Pulling back along a functor $f\colon A \to B$ in an $\infty$-cosmos $\mathcal{K}$ induces a functor of $\infty$-cosmoi $f^*\colon \mathcal{K}/B \to \mathcal{K}/A$.*

*Proof.* The isofibrations and trivial fibrations of $\mathcal{K}/A$ and $\mathcal{K}/B$ are simply those functors whose underlying functor in $\mathcal{K}$ is such. So the pullback stability property of the isofibrations and trivial fibrations of $\mathcal{K}$, axiom 2.1.1(c) immediately implies that $f^*\colon \mathcal{K}/B \to \mathcal{K}/A$ preserves these classes of functors. To complete this proof, we deploy a Yoneda embedding argument to show that $f^*$ also preserves the limits specified in 2.1.1(a).

We start by choosing our universe of sets in order to make $\mathcal{K}$ a small simplicial category and we consider the simplicial Yoneda embedding $\mathscr{Y}\colon \mathcal{K} \to \underline{\mathrm{sSet}}^{\mathcal{K}^{\mathrm{op}}}$. We know that this preserves and reflects all of the simplicial limits that exist in $\mathcal{K}$ and thus we know that the sliced functor $\mathscr{Y}_{/B}\colon \mathcal{K}/B \to \underline{\mathrm{sSet}}^{\mathcal{K}^{\mathrm{op}}}/\mathscr{Y}(B)$ preserves and reflects all of the limit types of Definition 2.1.1(a). It should be noted here that the codomain of this sliced functor is the full simplicially enriched slice whereas its domain is the restricted slice discussed in Example 2.1.11. It follows, therefore, that we get an essentially commutative diagram of simplicial functors

$$\begin{array}{ccc} \mathcal{K}/B & \xrightarrow{\quad f^* \quad} & \mathcal{K}/A \\ {\scriptstyle \mathscr{Y}_{/B}} \downarrow & \cong & \downarrow {\scriptstyle \mathscr{Y}_{/A}} \\ \underline{\mathrm{sSet}}^{\mathcal{K}^{\mathrm{op}}}/\mathscr{Y}(B) & \xrightarrow[\mathscr{Y}(f)^*]{} & \underline{\mathrm{sSet}}^{\mathcal{K}^{\mathrm{op}}}/\mathscr{Y}(A) \end{array}$$

in which the the bottom horizontal simplicial functor is pullback along $\mathscr{Y}(f)\colon \mathscr{Y}(A) \to \mathscr{Y}(B)$ between the full simplicial slice categories depicted; its existence is assured by the fact that $\underline{\mathrm{sSet}}^{\mathcal{K}^{\mathrm{op}}}$ admits all pullbacks. This pullback functor has a simplicial left adjoint $\Sigma_{\mathscr{Y}(f)}\colon \underline{\mathrm{sSet}}^{\mathcal{K}^{\mathrm{op}}}/\mathscr{Y}(A) \to \underline{\mathrm{sSet}}^{\mathcal{K}^{\mathrm{op}}}/\mathscr{Y}(B)$ given by composition with $\mathscr{Y}(f)\colon \mathscr{Y}(A) \to \mathscr{Y}(B)$, so in particular $\mathscr{Y}(f)^*$ preserves all limits. As the vertical functors preserve and reflect the limits that are assured to exist in an $\infty$-cosmos, it follows that $f^*\colon \mathcal{K}/B \to \mathcal{K}/A$ also preserves those as required. □

2.2. **Examples of $\infty$-cosmoi.** Some general model category theoretic results will be used to produce examples of $\infty$-cosmoi.

2.2.1. **Lemma.** *If $\mathcal{M}$ is a model category that is enriched as such over the Joyal model structure on simplicial sets, then the simplicial subcategory $\mathcal{M}_{\mathrm{fib}}$ spanned by its fibrant objects is an $\infty$-cosmos.* □

2.2.2. *Observation.* Of course, if $\mathcal{M}$ is a model category satisfying the conditions of the last lemma then so is its dual $\mathcal{M}^{\mathrm{op}}$. Consequently, it follows in that case that the dual $(\mathcal{M}_{\mathrm{cof}})^{\mathrm{op}}$ of the category of cofibrant objects in $\mathcal{M}$ is also an $\infty$-cosmos.

Indeed, the axioms of an $\infty$-cosmos and particularly the construction of its homotopy 2-category in §3.1 are considerably simplified with the additional hypothesis that "all objects are cofibrant." Observation 2.2.2, which will be exploited in a future paper, motivates our decision not to require this condition.

2.2.3. **Proposition.** *If $\mathcal{M}$ is a cartesian closed model category equipped with a Quillen adjunction to the Joyal model structure on $\underline{\mathrm{sSet}}$*

$$\mathcal{M} \underset{R}{\overset{L}{\underset{\longrightarrow}{\overset{\longleftarrow}{\bot}}}} \underline{\mathrm{sSet}}$$

*whose left adjoint preserves binary products, then $\mathcal{M}$ is enriched as a model category over the Joyal model structure, with hom-spaces defined by applying $R$ to the internal homs of $\mathcal{M}$ and simplicial cotensors defined by applying $L$ to the simplicial set and then forming the internal hom in $\mathcal{M}$. Consequently the full simplicial subcategory spanned by the fibrant objects in $\mathcal{M}$ defines an $\infty$-cosmos and the right adjoint restricts to define a functor $R\colon \mathcal{M}_{\mathrm{fib}} \to \underline{\mathrm{qCat}}$ of $\infty$-cosmoi.*

*Proof.* By a standard result in enriched category theory, a product-preserving left adjoint between cartesian closed categories provides its codomain with the structure of a category that is enriched, tensored, and cotensored over the domain [16, 3.7.11]; proofs of these facts follow easily from the Yoneda lemma. The enriched model category axioms are elementary consequences of the cartesian closed model category axioms on $\mathcal{M}$ and the fact that the adjunction defining the enrichment, tensors, and cotensors is Quillen. □

2.2.4. **Example** (the $\infty$-cosmos of categories)**.** Applying Proposition 2.2.3 to the strong monoidal Quillen adjunction

$$\underline{\mathrm{Cat}} \underset{N}{\overset{h}{\underset{\longrightarrow}{\overset{\longleftarrow}{\bot}}}} \underline{\mathrm{sSet}}$$

between the folk model structure on $\underline{\mathrm{Cat}}$ and the Joyal model structure on $\underline{\mathrm{sSet}}$, the category of small categories becomes an $\infty$-cosmos. The isofibrations are the usual isofibrations between categories. The hom-space between small categories $\mathbf{C}$ and $\mathbf{D}$ is the nerve of the functor category $\mathbf{D}^{\mathbf{C}}$, which is isomorphic to the quasi-category of maps between the nerves of $\mathbf{C}$ and $\mathbf{D}$. Weak equivalences are the usual notion of equivalence of categories. All categories are both fibrant and cofibrant in this model structure, so it follows that both of $\underline{\mathrm{Cat}}$ and $\underline{\mathrm{Cat}}^{\mathrm{op}}$ are $\infty$-cosmoi under this enrichment.

2.2.5. **Example** (the $\infty$-cosmos of complete Segal spaces)**.** Precomposing with the adjoint pair of functors

$$\overset{p_1}{\underset{i_1}{\underset{\bot}{\leftrightarrows}}} \quad \times$$

defined by $i_1([n]) = [n] \times [0]$ and $p_1([n] \times [m]) = [n]$ induces an adjunction

$$\underline{\mathrm{ssSet}} \overset{p_1^*}{\underset{i_1^*}{\underset{\bot}{\leftrightarrows}}} \underline{\mathrm{sSet}}$$

The right adjoint $i_1^*$ takes a bisimplicial set to its $0^{\text{th}}$ row; for a complete Segal space, the $0^{\text{th}}$ row defines its underlying quasi-category. The pair $p_1^* \dashv i_1^*$ defines a Quillen equivalence between Joyal's model structure for quasi-categories and Rezk's model structure for complete Segal spaces [10, 4.11]. Both model structures are cartesian closed and the left adjoint $p_1^*$ preserves products. By Proposition 2.2.3, it follows that this strong monoidal Quillen adjunction can be used to convert the cartesian closed model structure on $\underline{\mathrm{ssSet}}$ into an enrichment of the model category of complete Segal spaces over the Joyal model structure. In this way, the full subcategory of complete Segal spaces defines an $\infty$-cosmos whose hom-spaces are simply the underlying quasi-categories of the internal homs. With respect to this definition, the underlying quasi-category functor $i_1^* \colon \underline{\mathrm{CSS}} \to \underline{\mathrm{qCat}}$ defines a functor of $\infty$-cosmoi.

2.2.6. **Example.** Joyal and Tierney also describe a second Quillen equivalence

$$\underline{\mathrm{sSet}} \overset{t_!}{\underset{t^!}{\underset{\bot}{\leftrightarrows}}} \underline{\mathrm{ssSet}}$$

between the model structure for quasi-categories and the model structure for complete Segal spaces whose right adjoint sends a simplicial set $A$ to the bisimplicial set

$$(t^! A)_{m,n} := \underline{\mathrm{sSet}}(\Delta^m \times \widetilde{\Delta^n}, A)$$

whose value at a pair of objects $[m], [n] \in$ is the set of simplicial maps to $A$ from the product of the ordinal category $[m]$ with the groupoid obtained by freely inverting the morphisms in $[n]$ [10, 4.12]. As a right Quillen adjoint between model categories whose fibrant objects are cofibrant, $t^! \colon \underline{\mathrm{qCat}} \to \underline{\mathrm{CSS}}$ preserves isofibrations and the conical limits listed in Definition 2.1.1. To show that it is a functor of $\infty$-cosmoi, it remains only to show that $t^!$ is simplicially enriched and preserves simplicial cotensors, or equivalently, to show that the adjunction $t_! \dashv t^!$ is simplicially enriched.

We adopt notation from [10], and in particular, make use of the "external product" bifunctor $-\square- \colon \underline{\mathrm{sSet}} \times \underline{\mathrm{sSet}} \to \underline{\mathrm{ssSet}}$. Note that $\Delta^n \square \Delta^m$ is the representable bisimplicial set at the object $([n], [m]) \in \quad \times \quad$. Using exponential notation for the internal hom in both $\underline{\mathrm{sSet}}$ and $\underline{\mathrm{ssSet}}$, we wish to show that for any simplicial set $A$ and bisimplicial set $X$ that the hom-simplicial sets $A^{t_! X} \cong i_1^*((t^! A)^X)$ are isomorphic; recall $i_1^*$ returns the 0th row of the bisimplicial set $(t^!A)^X$. Because $X$ may be recovered as a conical colimit of representable bisimplicial sets, it suffices to show this in the case $X = \Delta^n\Box\Delta^m$. We have

$$\begin{aligned} i_1^*((t^!A)^{\Delta^n\Box\Delta^m})_k &= ((t^!A)^{\Delta^n\Box\Delta^m})_{k,0} \\ &= \underline{\mathrm{ssSet}}((\Delta^n\Box\Delta^m)\times(\Delta^k\Box\Delta^0), t^!A) \\ &\cong \underline{\mathrm{sSet}}(t_!((\Delta^n\Box\Delta^m)\times(\Delta^k\Box\Delta^0)), A) \\ &\cong \underline{\mathrm{sSet}}(t_!((\Delta^n\times\Delta^k)\Box(\Delta^m\times\Delta^0)), A) \\ &\cong \underline{\mathrm{sSet}}(t_!((\Delta^n\times\Delta^k)\Box\Delta^m), A). \end{aligned}$$

By [10, 2.11] this is isomorphic to:

$$\begin{aligned} &\cong \underline{\mathrm{sSet}}((\Delta^n\times\Delta^k)\times\widetilde{\Delta^m}, A) \\ &\cong \underline{\mathrm{sSet}}((\Delta^n\times\widetilde{\Delta^m})\times\Delta^k, A) \end{aligned}$$

which again by [10, 2.11] is isomorphic to:

$$\cong \underline{\mathrm{sSet}}(t_!(\Delta^n\Box\Delta^m)\times\Delta^k, A) = (A^{t_!(\Delta^n\Box\Delta^m)})_k.$$

2.2.7. **Example** (the $\infty$-cosmos of Segal categories)**.** The category $\underline{\mathrm{ssSet}}$ of bisimplicial sets has a full subcategory $\underline{\mathrm{PCat}}$ of *precategories*, those bisimplicial sets whose $0^{\text{th}}$ column is discrete. The category of precategories bears a cartesian closed model structure whose fibrant objects are the *Segal categories* (precategories satisfying the Segal condition) by results of Hirschowitz-Simpson, Pellisier, Bergner, and Joyal; see [10, §5]. Joyal and Tierney demonstrate that there is a Quillen equivalence between this model category and the model structure for quasi-categories [10, 5.6]. The right adjoint $\underline{\mathrm{PCat}} \to \underline{\mathrm{sSet}}$ carries a Segal category to its $0^{\text{th}}$ row, the underlying quasi-category. The left adjoint carries a simplicial set to the external product with the terminal simplicial set. This functor preserves products, so again by Proposition 2.2.3 this strong monoidal Quillen adjunction induces an enrichment of the model structure for Segal categories over the model structure for quasi-categories. So the full subcategory spanned by the Segal categories defines an $\infty$-cosmos such that the underlying quasi-category functor is a functor of $\infty$-cosmoi.

2.2.8. **Example** (the $\infty$-cosmos of marked simplicial sets)**.** The category $\underline{\mathrm{msSet}}$ of marked simplicial sets bears a model structure in which the cofibrations are the monomorphisms and the weak equivalences are those maps $X \to Y$ for which the induced map

$$\mathrm{Map}^\flat(Y, Z^\natural) \to \mathrm{Map}^\flat(X, Z^\natural)$$

is an equivalence of quasi-categories for all naturally marked quasi-categories $Z^\natural$; here $\mathrm{Map}^\flat$ denotes the underlying simplicial set of the internal hom in the category of marked simplicial sets; see [12, 3.1.3]. With these mapping spaces, this model category is enriched over the Joyal model structure [12, 3.1.4.5]. Thus the full simplicial subcategory spanned by the fibrant objects, the naturally marked quasi-categories, defines an $\infty$-cosmos. This enrichment of the marked model structure can be understood as an application of Proposition 2.2.3 to the strong monoidal Quillen adjunction

$$\underline{\mathrm{msSet}} \underset{U}{\overset{(-)^\flat}{\underset{\perp}{\rightleftarrows}}} \underline{\mathrm{sSet}}$$

in which the left adjoint marks only degenerate 1-simplices and the right adjoint forgets the markings. Thus the underlying quasi-category functor, forgetting the markings, is a functor of $\infty$-cosmoi.

In all of the examples to this point all objects in the $\infty$-cosmoi discussed have been cofibrant. In the next proposition we generalise Example 2.2.5 to define an $\infty$-cosmos of Rezk objects in any left-proper combinatorial model category $\mathcal{M}$. At that level of generality we cannot show that all objects will be cofibrant in the resulting $\infty$-cosmos.

Special cases of the next result include Rezk's model structure for complete Segal spaces and Barwick's model structure for $n$-fold complete Segal spaces. It is actually the case in both of those examples that all objects of the resulting $\infty$-cosmos are cofibrant.

2.2.9. **Proposition** (the $\infty$-cosmos of Rezk objects). *If $\mathcal{M}$ is a left-proper combinatorial model category, then the Reedy model structure on $\mathcal{M}^{\Delta^{\mathrm{op}}}$ admits a left Bousfield localization whose fibrant objects are* Rezk objects*: Reedy fibrant simplicial objects satisfying the Segal and completeness conditions. This model category is enriched over the Joyal model structure for quasi-categories and so its full subcategory spanned by fibrant objects is an $\infty$-cosmos.*

*Proof.* Consider the set $J = J_{\mathrm{Segal}} \cup J_{\mathrm{completeness}}$, where

$$J_{\mathrm{Segal}} = \{\Lambda^{n,k} \hookrightarrow \Delta^n \mid n \geq 2, 0 \leq k \leq n\} \quad \text{and} \quad J_{\mathrm{completeness}} = \{\Delta^0 \hookrightarrow \mathbb{I}\}.$$

We say that a Reedy fibrant simplicial object $X$ is a *Rezk object* if the induced maps $\{j, X\}$ on weighted limits are trivial fibrations in $\mathcal{M}$ for all $j \in J$; note that Reedy fibrancy implies already that these maps are fibrations. The "spine inclusions" $\Delta^1 \cup_{\Delta^0} \Delta^1 \cup_{\Delta^0} \cdots \cup_{\Delta^0} \Delta^1 \hookrightarrow \Delta^n$ are contained in the weak saturation of the inner horn inclusions. Thus, the maps in $J_{\mathrm{Segal}}$ impose the Segal condition on the Reedy fibrant objects. The map in $J_{\mathrm{completeness}}$ imposes the completeness condition, by the 2-of-3 property; see [14, §6].

The category $\mathcal{M}^{\Delta^{\mathrm{op}}}$ is enriched, tensored, and cotensored over simplicial sets in such a way that Leibniz tensors monomorphisms of simplicial sets with (trivial) Reedy cofibrations are (trivial) Reedy cofibrations; see e.g., [4, 4.4]. We write $\otimes$ for the simplicial tensor and Map for the hom-spaces. We will apply Jeff Smith's theorem [2, 1.7] to prove that $\mathcal{M}^{\Delta^{\mathrm{op}}}$ admits a model structure whose

- cofibrations are the Reedy cofibrations,
- fibrant objects are the *Rezk objects*, those Reedy fibrant simplicial objects $X$ for which $\{j, X\}$ is a trivial fibration in $\mathcal{M}$ for all $j \in J$,
- fibrations are the *Rezk fibrations*, those Reedy fibrations $p\colon X \to Y$ for which the Leibniz weighted limits $\widehat{\{j, p\}}$ are trivial fibration in $\mathcal{M}$ for all $j \in J$,

- weak equivalences are the *Rezk weak equivalences*, those maps $w\colon U \to V$ so that $w^*\colon \mathrm{Map}(V,X) \to \mathrm{Map}(U,X)$ is an equivalence of quasi-categories for all Rezk objects $X$.

We call the localized model structure on $\mathcal{M}^{\Delta^{\mathrm{op}}}$ the *model structure for Rezk objects*.

By adjunction, $p\colon X \to Y$ is a Rezk fibration between Rezk objects if and only if the map $\widehat{\mathcal{M}}(i,p)\colon \mathcal{M}(B,X) \to \mathcal{M}(A,X) \times_{\mathcal{M}(A,Y)} \mathcal{M}(B,Y)$ is an isofibration of quasi-categories for all generating cofibrations $i\colon A \to B$ in $\mathcal{M}$. This is the case if and only if these maps lift against the set $J \mathbin{\widehat{\times}} B$, where $B = \{\partial\Delta^n \hookrightarrow \Delta^n \mid n \geq 0\}$; one direction of this implication is obvious and proofs of the other can be found in [5, §§A-B]. Transposing again, we see that $p$ is a Rezk fibration between Rezk objects if and only if $p$ has the right lifting property against the set of maps $(j \mathbin{\widehat{\times}} b) \mathbin{\widehat{*}} i \cong j \mathbin{\widehat{\otimes}} (b \mathbin{\widehat{*}} i)$ for all $j \in J$, $b \in B$, and generating cofibrations $i$ of $\mathcal{M}$, where $*$ denotes the pointwise tensor functor $- * -\colon \underline{\mathrm{sSet}} \times \mathcal{M} \to \mathcal{M}^{\;\;\mathrm{op}}$. The set $B \mathbin{\widehat{*}} i$, where $i$ ranges over the generating cofibrations of $\mathcal{M}$, defines a set of generating cofibrations for the Reedy model structure on $\mathcal{M}^{\;\;\mathrm{op}}$ [17, 7.7]. So, by adjunction, $p$ is a Rezk fibration between Rezk objects if and only if $\widehat{\mathrm{Map}}(f,p)$ is an isofibration of quasi-categories, for every Reedy cofibration $f$.

We easily verify the conditions of Jeff Smith's theorem. By definition, the Rezk weak equivalences form an accessible subcategory of the arrow category $\mathcal{M}^{\;\;\mathrm{op}}$ and this class is closed under weak equivalences and satisfies the 2-of-3 property. For any Reedy trivial cofibration $w\colon U \to V$ and Reedy fibrant object $X$, $w^*\colon \mathrm{Map}(V,X) \to \mathrm{Map}(U,X)$ is a trivial fibration of simplicial sets. By Ken Brown's lemma, $w^*\colon \mathrm{Map}(V,X) \to \mathrm{Map}(U,X)$ is still an equivalence of quasi-categories of $w$ is only a Reedy weak equivalence. So the Reedy weak equivalences are Rezk weak equivalences, proving that the Reedy trivial fibrations are weak equivalences. The intersection of the Reedy cofibrations and Rezk weak equivalences are those maps $w$ so that $w^*\colon \mathrm{Map}(V,X) \to \mathrm{Map}(U,X)$ is a trivial fibration between quasi-categories for all fibrant objects $X$. This condition is closed under pullbacks and limits of towers, so the class of Reedy cofibrations and weak equivalences is closed under pushout and transfinite composition. Jeff Smith's theorem now implies that the model structure for Rezk objects exists.

To prove that this model category is enriched over the Joyal model structure, recall that the Reedy model structure on $\mathcal{M}^{\;\;\mathrm{op}}$ satisfies the $2/3^{\mathrm{rd}}$'s of the SM7 axiom having to do with Leibniz products with monomorphisms between simplicial sets. It follows easily that the localized model structure again has this property: apply [4, 3.2.a], the fact that the cofibrations are unchanged, and a general fact about Bousfield localizations: that the Rezk fibrations between Rezk objects are the Reedy fibrations between Rezk objects [7, 3.3.16]. To show that the model structure for Rezk objects is enriched as a model category over the Joyal model structure, we apply [4, 3.2.b], making use of the fact that a Bousfield localization of a Reedy model structure on a left property model category is again left proper. Now it suffices to prove that for every Rezk trivial cofibration $w\colon U \to V$ and every Rezk fibrant object $X$, the map $w^*\colon \mathrm{Map}(V,X) \to \mathrm{Map}(U,X)$ is an isofibration of quasi-categories. But this is immediate from how these classes are defined.

Now for any model categories $\mathcal{M}$ for which the Rezk objects are Reedy cofibrant, it follows immediately that the full sub simplicial category spanned by the Rezk objects defines an $\infty$-cosmos. □

2.2.10. **Example** (the $\infty$-cosmos of $\theta_n$-spaces)**.** Another model for $(\infty, n)$-categories is given by the fibrant objects in a cartesian closed model structure due to Charles Rezk in the category of $\underline{\mathrm{sSet}}$-valued presheaves on Joyal's category $\theta_n$ [15]. The cofibrations in this model structure are the monomorphisms, so to prove that the fibrant objects form an $\infty$-cosmos it suffices, by Proposition 2.2.3, to find a strong monoidal Quillen left adjoint from Joyal's model structure on $\underline{\mathrm{sSet}}$ to Rezk's model structure on $\underline{\mathrm{sSet}}^{\theta_n^{\mathrm{op}}}$. The desired adjunction is a composite of the adjunction of Example 2.2.5 and a strong monoidal Quillen adjunction between complete Segal spaces and $\theta_n$-spaces.

Recall that the category $\theta_n$ is defined inductively to be the categorical wreath product $\Delta \wr \theta_{n-1}$, where $\theta_1 = \Delta$. Each category $\theta_n$ has a terminal object. For any small category $\mathbf{C}$ with a terminal object $t$, the adjunction

$$1 \underset{t}{\overset{!}{\underset{\longrightarrow}{\overset{\longleftarrow}{\bot}}}} \mathbf{C}$$

induces an adjunction

$$\Delta = \Delta \wr 1 \underset{\Delta \wr t}{\overset{\Delta \wr !}{\underset{\longrightarrow}{\overset{\longleftarrow}{\bot}}}} \Delta \wr \mathbf{C}$$

upon applying the 2-functor $\Delta \wr - \colon \underline{\mathrm{Cat}} \to \underline{\mathrm{Cat}}$. In the case $\mathbf{C} = \theta_{n-1}$, $\Delta \wr t \colon \Delta \hookrightarrow \theta_n$ defines an inclusion of $\Delta$ into $\theta_n$ as "pasting diagrams comprised only of composable 1-cells." On taking simplicial presheaves, this defines an adjunction

$$\underline{\mathrm{sSet}}^{\theta_n^{\mathrm{op}}} \underset{(\Delta \wr t)^*}{\overset{(\Delta \wr !)^*}{\underset{\longrightarrow}{\overset{\longleftarrow}{\bot}}}} \underline{\mathrm{sSet}}^{\Delta^{\mathrm{op}}} \tag{2.2.11}$$

The left adjoint here has a further left adjoint, defined by left Kan extension along $\Delta \wr ! \colon \theta_n^{\mathrm{op}} \to \Delta^{\mathrm{op}}$, and so it preserves products.

It remains to argue that this adjunction is Quillen. Both the model structure for $\theta_n$-spaces and for complete Segal spaces are defined as left Bousfield localizations of the injective (or, equivalently, the Reedy) model structures on the functor categories, defined relative to Quillen's model structure on $\underline{\mathrm{sSet}}$. Because the left adjoint (2.2.11) is simply a restriction functor, it is manifestly left Quillen between these injective model structures. Thus, (2.2.11) descends to a Quillen adjunction between the localized model structures if and only if the right adjoint preserves local (fibrant) objects. A functor $X \colon \theta_n^{\mathrm{op}} \to \underline{\mathrm{sSet}}$ is fibrant if it satisfies "Segal" and "completeness" conditions, defined using weighted limits [15, §7]. It follows that the underlying bisimplicial sets of these "complete Segal objects" are complete Segal spaces. Thus, Proposition 2.2.3 applies and we conclude that $\theta_n$-spaces define an $\infty$-cosmos, and that the underlying complete Segal space functor defines a functor of $\infty$-cosmoi.

## 3. 2-CATEGORY THEORY IN AN $\infty$-COSMOS

Our aim in this section is to demonstrate that all of the basic 2-category theory of quasi-categories developed in sections 3, 4, and 5 of [18] now generalises to provide corresponding results in the *homotopy 2-category* of any $\infty$-cosmos $\mathcal{K}$. We introduce this 2-category in §3.1: it is the strict 2-category spanned by the cofibrant objects whose hom-categories are defined to be homotopy categories of the hom-quasi-categories.

Indeed, almost all of the work in those sections can be repeated in an *abstract homotopy 2-category*, an arbitrary 2-category which possesses certain weak 2-dimensional limits. In §3.2, we introduce suitably weak notions of *comma objects* and *iso-comma objects* in a general 2-category. An abstract homotopy 2-category is simply a strict 2-category possessing weak commas and weak iso-commas. In §3.3, we explain how comma objects and iso-comma objects are constructed in the homotopy 2-category of an $\infty$-cosmos. This proves that the homotopy 2-category is an abstract homotopy 2-category, as our terminology would suggest. In §3.4, we develop the basic general theory of comma and iso-comma constructions in an abstract homotopy 2-category.

In §3.5, we consider slices of an abstract homotopy 2-category and describe a pullback operation. Given a functor $f\colon A \to B$, any isofibration over $B$ pulls back to define an isofibration over $A$. This process is not 2-functorial, as one might expect, because pullbacks, which are closely related to iso-commas, are also weak. Nonetheless important categorical structures, notably fibred equivalences and adjunctions, can be pulled back along any functor, as we prove in §3.6. Thus pullback provides a sufficiently well-behaved reindexing operation, which will frequently be exploited.

### 3.1. **The homotopy 2-category of an $\infty$-cosmos.** 
In this section we introduce the *homotopy 2-category* associated to an $\infty$-cosmos $\mathcal{K}$.

3.1.1. **Definition** (the 2-categorical quotient of an $\infty$-cosmos)**.** An $\infty$-cosmos $\mathcal{K}$ is a simplicially enriched category so we may apply Observation I.3.1.2 to construct a 2-category $h_*\mathcal{K}$. This has the same objects as $\mathcal{K}$ and has hom-categories given by taking the homotopy categories of its hom-spaces, that is $\hom(A,B) := h(\mathrm{map}(A,B))$. The horizontal composition operation of this 2-category is constructed by applying $h\colon \underline{\mathrm{sSet}} \to \underline{\mathrm{Cat}}$ to the composition functions of $\mathcal{K}$, using the fact that the homotopy category construction preserves finite products

3.1.2. **Definition** (the homotopy 2-category)**.** The *homotopy 2-category* of an $\infty$-cosmos $\mathcal{K}$, is the full sub-2-category $\mathcal{K}_2$ of $h_*\mathcal{K}$ spanned by the cofibrant objects in $\mathcal{K}$. That is:

- the objects in the homotopy 2-category are the cofibrant objects in the $\infty$-cosmos, and
- the hom-categories in the homotopy 2-category are the homotopy categories of the hom quasi-categories between cofibrant objects in the $\infty$-cosmos
$$\hom(A,B) := h(\mathrm{map}(A,B)).$$

Of course, if all objects in the $\infty$-cosmos $\mathcal{K}$ are cofibrant, then the homotopy 2-category and the 2-category $h_*\mathcal{K}$ coincide. Before moving on to explore the properties of the homotopy 2-category $\mathcal{K}_2$ we recall a few 2-categorical notions:

3.1.3. *Recall* (equivalences, isofibrations, and surjective equivalences). A 1-cell $p\colon E \to B$ in a 2-category $\mathcal{C}$ is a *(representable) isofibration* if the functor $\hom(X,p)\colon \hom(X,E) \to \hom(X,B)$ is an isofibration of categories for all objects $X$. This latter condition is equivalent to postulating that any invertible 2-cell as depicted on the left of the following diagram

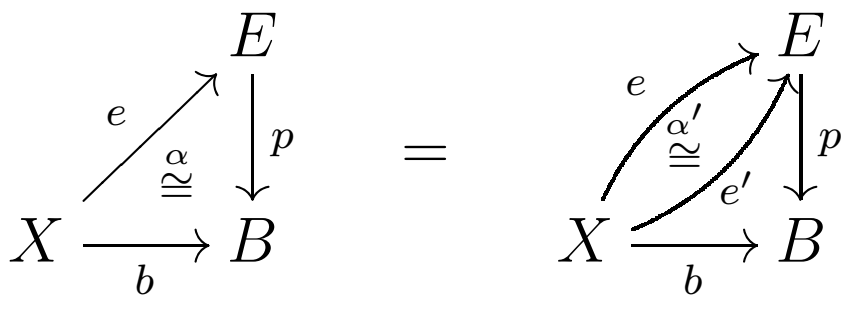

has a lift along $p$ as depicted on the right.

A 1-cell $w\colon A \to B$ is an *equivalence* in $\mathcal{C}$ if and only if there is a 1-cell $w'\colon B \to A$ and invertible 2-cells $ww' \cong \mathrm{id}_B$ and $w'w \cong \mathrm{id}_A$. The equivalences in $\mathcal{C}$ are precisely those 1-cells $w\colon A \to B$ for which $\hom(X,w)\colon \hom(X,A) \to \hom(X,B)$ is an equivalence of categories.

More explicitly, a 1-cell $w\colon A \to B$ is an equivalence if and only if it satisfies two representably defined properties. First it must be (representably) fully faithful, in the sense that if we have a 2-cell $\lambda\colon wa \Rightarrow wa'$ then there exists a unique 2-cell $\bar{\lambda}\colon a \Rightarrow a'$ with $w\bar{\lambda} = \lambda$. Secondly it must be (representably) essentially surjective, in the sense that if we have a 1-cell $b\colon X \to B$ then there exists a 1-cell $a\colon X \to A$ and an invertible 2-cell $wa \cong b$. A well-known result demonstrates that any equivalence in a 2-category can be promoted to an adjoint equivalence by changing one of the 2-cells (cf. [13, IV.4.1]).

Combining these notions, a 1-cell $p\colon E \to B$ of the 2-category $\mathcal{C}$ is a *(representable) surjective equivalence* if it is both an equivalence and an isofibration. This condition holds precisely when for all objects $X$ the functor $\hom(X,p)\colon \hom(X,E) \to \hom(X,B)$ is a surjective equivalence of categories, that is to say an equivalence of categories which acts surjectively on objects. We may use isomorphism lifting to show that any representable trivial fibration $p\colon E \to B$ has a right inverse $r\colon B \to E$ for which there is an isomorphism $\gamma\colon rp \cong \mathrm{id}_E$ with $\gamma r = \mathrm{id}_r$ and $p\gamma = \mathrm{id}_p$.

3.1.4. **Lemma.** *If $X$ is a cofibrant object and $p\colon E \twoheadrightarrow B$ is an isofibration (resp. trivial fibration) in an $\infty$-cosmos $\mathcal{K}$, then the functor* $\hom(X,p)\colon \hom(X,E) \to \hom(X,B)$ *of hom-categories of $h_*\mathcal{K}$ is an isofibration (resp. surjective equivalence) of categories.*

*Proof.* Using the explicit description of the homotopy category of a quasi-category, as presented in I.2.2.2 for example, it is easy to check that if a functor $p\colon E \to B$ of quasi-categories is an isofibration (resp. trivial fibration) then the functor $h(p)\colon h(E) \to h(B)$ of homotopy categories is an isofibration (resp. surjective equivalence). By Lemma 2.1.8 we know that the functor $\mathrm{map}(X,p)\colon \mathrm{map}(X,E) \to \mathrm{map}(X,B)$ is an isofibration (resp. trivial fibration) of quasi-categories so by taking homotopy categories it follows that the functor $\hom(X,p)\colon \hom(X,E) \to \hom(X,B)$ is an isofibration (resp. surjective equivalence) of categories as required. □

3.1.5. **Corollary.** *Isofibrations (resp. trivial fibrations) between cofibrant objects in an $\infty$-cosmos $\mathcal{K}$ are representable isofibrations (resp. representable surjective equivalences) in its homotopy 2-category $\mathcal{K}_2$.* □

3.1.6. *Observation* (functors representing invertible 2-cells). Consider an invertible 2-cell $\alpha\colon f \cong g\colon A \to B$ in $h_*\mathcal{K}$ whose domain is cofibrant. Then Lemma 2.1.8 tells us that the hom-space $\mathrm{map}(A, B)$ is a quasi-category, so we may apply Recollection 2.1.5(a) to show that there exists some simplicial map $\bar{\alpha}\colon \mathbb{I} \to \mathrm{map}(A, B)$ that maps the 1-simplex with vertex sequence $\{0, 1\}$ in $\mathbb{I}$ to a 1-simplex in $\mathrm{map}(A, B)$ which is a representative of the 2-cell $\alpha$ in $\hom(A, B) = h(\mathrm{map}(A, B))$. Transposing, we obtain a corresponding functor $\hat{\alpha}\colon A \to \mathbb{I} \pitchfork B$ for which $p_0\hat{\alpha} = f$ and $p_1\hat{\alpha} = g$. Indeed, invertible 2-cells $\alpha\colon f \Rightarrow g$ are in bijective correspondence with isomorphism classes of 1-cells

$$\begin{array}{ccc} A & \xrightarrow{\hat{\alpha}} & \mathbb{I} \pitchfork B \\ & {\scriptstyle (g,f)} \searrow \quad \swarrow {\scriptstyle (p_1,p_0)} & \\ & B \times B & \end{array}$$

in the slice 2-category $(h_*\mathcal{K})/B \times B$, a fact which we leave as an exercise for the reader.

3.1.7. **Lemma.** *If $\alpha\colon f \cong g\colon A \to B$ is an invertible 2-cell in $h_*\mathcal{K}$ and $A$ is cofibrant, then $f$ is a weak equivalence if and only if $g$ is a weak equivalence.*

*Proof.* Suppose $f$ is a weak equivalence. By the 2-of-3 property, the equality $f = p_0\hat{\alpha}$, and the fact that $f$ and $p_0$ are weak equivalences, we see that the functor $\hat{\alpha}\colon A \to K \pitchfork B$ is also a weak equivalence. Applying, 2-of-3 a second time, we infer that $g$ is a weak equivalence. $\square$

3.1.8. **Proposition.** *A functor $w\colon A \to B$ between cofibrant objects is a weak equivalence in the $\infty$-cosmos $\mathcal{K}$ if and only if it is an equivalence in the homotopy 2-category $\mathcal{K}_2$.*

*Proof.* If $w\colon A \to B$ is an equivalence in the 2-category $\mathcal{K}_2$ with equivalence inverse $w'\colon B \to A$ then applying Lemma 3.1.7 to the isomorphisms $ww' \cong \mathrm{id}_B$ and $w'w \cong \mathrm{id}_A$ we see that $ww'$ and $w'w$ are both weak equivalences since all identities are such. By the 2-of-6 property, it follows that $w$ is a weak equivalence in our $\infty$-cosmos $\mathcal{K}$.

For the converse, Brown's factorisation lemma, Lemma 2.1.6, tells us in particular that any weak equivalence $w\colon A \to B$ in the $\infty$-cosmos $\mathcal{K}$ can be factored as $w = pj$ where $p\colon N_f \xrightarrow{\sim} B$ is a trivial fibration and $j\colon A \xrightarrow{\sim} N_f$ is a weak equivalence right inverse to a trivial fibration. Applying Lemma 3.1.4, we know that if $X$ is a cofibrant object then the functor $\hom(X, p)\colon \hom(X, N_f) \xrightarrow{\sim} \hom(X, B)$ is a surjective equivalence of hom-categories and that $\hom(X, i)\colon \hom(X, A) \xrightarrow{\sim} \hom(X, N_f)$ is right inverse to a surjective equivalence and is thus an equivalence of categories. On composing these see that $\hom(X, w)\colon \hom(X, A) \xrightarrow{\sim} \hom(X, B)$ is also an equivalence of categories, so it follows that if $w\colon A \to B$ is a weak equivalence of cofibrant objects then it is also an equivalence in the 2-category $\mathcal{K}_2$. $\square$

We next describe the 2-universal properties of the cofibrant replacement process.

3.1.9. *Observation* (cofibrant replacement as a *bi-coreflection*). For any object $A$ of an $\infty$-cosmos $\mathcal{K}$ the cofibrant replacement axiom 2.1.1(e) supplies us with a cofibrant object $\tilde{A}$ and a trivial fibration $r_A\colon \tilde{A} \xrightarrow{\sim\!\!\!\!\to\!\!\!\!\to} A$. Now if $X$ is any cofibrant object then Lemma 2.1.8 tells us that $\hom(X, r_A)\colon \hom(X, \tilde{A}) \to \hom(X, A)$ is a surjective equivalence of categories.

What this tells us, in particular, is that the inclusion 2-functor $\mathcal{K}_2 \hookrightarrow h_*\mathcal{K}$ has a *right bi-adjoint* determined by cofibrant replacement. In other words, we may extend the cofibrant replacement operation to a *pseudo-functor*, also called a *homomorphism*, $\widetilde{(-)}\colon h_*\mathcal{K} \to \mathcal{K}_2$ in such a way that the $r_A\colon \tilde{A} \xrightarrow{\sim\!\!\!\!\to\!\!\!\!\to} A$ become the components of a 2-natural counit transformation which induces surjective equivalences between the adjoint hom-categories $\hom(X, \tilde{A})$ of $\mathcal{K}_2$ and $\hom(X, A)$ of $h_*\mathcal{K}$.

It should be noted that we will not require much of the general theory of bi-adjunctions here, and certainly the reader should be able to follow our arguments in this regard without any formal preparation in the yoga of such things. However, the inquisitive reader may wish to consult §1 of Street's classic paper [24] for a complete account.

3.1.10. **Definition.** If $F\colon \mathcal{K} \to \mathcal{L}$ is a functor of $\infty$-cosmoi, then the *induced homomorphism of homotopy 2-categories* is defined to be the composite

$$F_2 := \mathcal{K}_2 \hookrightarrow h_*\mathcal{K} \xrightarrow{h_*F} h_*\mathcal{L} \xrightarrow{\widetilde{(-)}} \mathcal{L}_2.$$

In practice, there often exists a simplicially functorial cofibrant replacement satisfying Definition 2.1.1(e), in which case $F_2\colon \mathcal{K}_2 \to \mathcal{L}_2$ is a 2-functor; for instance, this is the case when all objects in $\mathcal{L}_2$ are cofibrant. But in practice, it makes no great difference if $F_2$ is only a pseudo-functor: in any case, it preserves adjunctions and equivalences in $\mathcal{K}_2$ in addition to other important structures.

3.2. **Abstract homotopy 2-categories.** In much of what follows, all of our arguments will essentially be 2-categorical in nature. To stress this point we will often work in an abstract 2-category satisfying very simple axioms. The following axiomatisation is not intended to be a complete or exhaustive account of the 2-categorical structures possessed by the homotopy 2-category associated with an $\infty$-cosmos. It simply encapsulates the essential 2-categorical properties and constructions we shall need here in order to develop the basic theory of cartesian fibrations.

An *abstract homotopy 2-category* is a strict 2-category admitting *comma* and *iso-comma* constructions characterised by suitably weakened 2-universal properties. In 2-category theory such constructions are usually required to possess a *strict* 2-universal property which demands that certain canonical comparison functors between hom-categories are isomorphisms. However, herein we ask only for a *weak* 2-universal property under which these canonical comparisons possess the weaker property of being *smothering functors*, that is to say they are surjective on objects, full, and conservative (see Definition I.3.3.1).

We study this weak 2-universal property because it is this, and not its strong counterpart, which characterises comma and iso-comma constructions in the the homotopy 2-category associated with an $\infty$-cosmos. Consequently, from hereon the term *comma object* in a 2-category $\mathcal{C}$ will refer exclusively to the *weak comma objects* introduced in §I.3.3, whose definition we now recall.

3.2.1. **Definition** (comma object)**.** Given morphisms $f \colon B \to A$ and $g \colon C \to A$, a *comma object* is given by the data

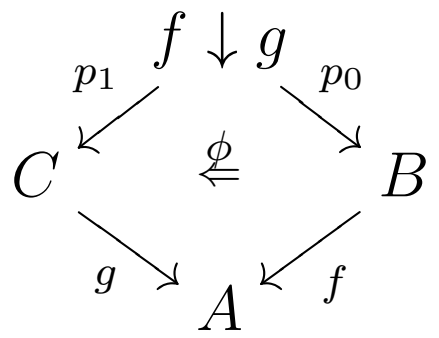

(3.2.2)

so that for each object $X \in \mathcal{C}$ the induced comparison functor of hom-categories

$$\hom(X, f \downarrow g) \longrightarrow \hom(X, f) \downarrow \hom(X, g)$$

is a *smothering functor*: a functor which is surjective on objects, full, and conservative.

Explicitly, this weak universal property supplies us with three operations:

(i) (1-cell induction) Given a *comma cone* $\alpha \colon fb \Rightarrow gc$

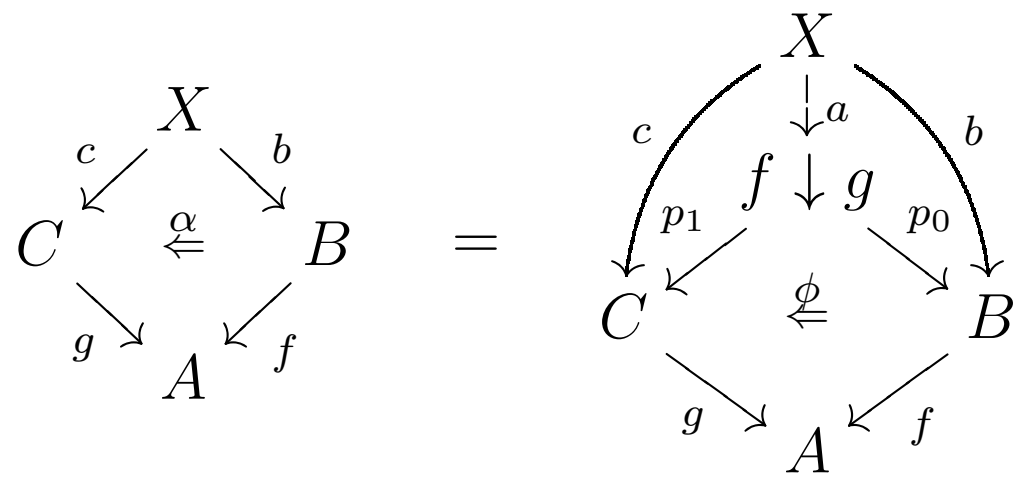

over the pair of functors $f$ and $g$, there exists a 1-cell $a \colon X \to f \downarrow g$ so that $p_0 a = b$, $p_1 a = c$, and $\alpha = \phi a$.

(ii) (2-cell induction) Given a pair of functors $a, a' \colon X \to f \downarrow g$ and a pair of 2-cells $\tau_0 \colon p_0 a \Rightarrow p_0 a'$ and $\tau_1 \colon p_1 a \Rightarrow p_1 a'$ which are compatible in the sense that $\phi a' \cdot f\tau_0 = g\tau_1 \cdot \phi a$, then there exists a 2-cell $\tau \colon a \Rightarrow a'$ with $p_0\tau = \tau_0$ and $p_1 \tau = \tau_1$.

(iii) (conservativity) Any 2-cell $\tau \colon a \Rightarrow a' \colon X \to f \downarrow g$ with the property that the whiskered 2-cells $p_0\tau$ and $p_1\tau$ are both isomorphisms is also an isomorphism.

We refer to (3.2.2) as a *comma square* and $C \xleftarrow{p_1} f \downarrow g \xrightarrow{p_0} B$ as a *comma span*.

When $f$ or $g$ is an identity, we write $A \downarrow g$ or $f \downarrow A$, respectively, for the comma object. In the case where both $f$ and $g$ are identities, we write $A^2$ for $A \downarrow A$ because this object is a *weak* 2*-cotensor*, in the sense introduced in §I.3.3.

A *iso-comma object* in $\mathcal{C}$ will refer exclusively to what might be called a *weak iso-comma object*, defined in analogy with 3.2.1 except with an additional requirement that the 2-cells in iso-comma cones are always invertible. Iso-comma objects are closely related to pullbacks, as we shall see in Lemma 3.5.6, hence our choice of notation.

3.2.3. **Definition** (iso-comma object)**.** Given morphisms $f\colon B \to A$ and $g\colon C \to A$, an *iso-comma object* is given by the data

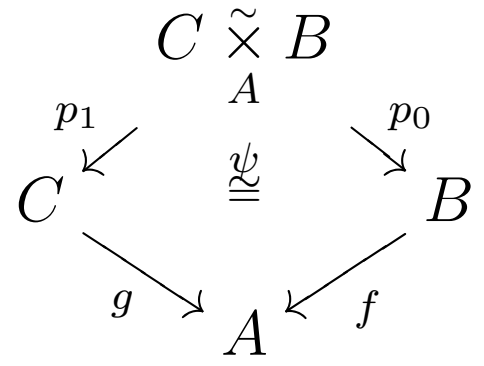

(3.2.4)

i.e., a span together with an invertible 2-cell $\psi\colon fp_0 \cong gp_1$, so that for each object $X \in \mathcal{C}$ the induced comparison functor of hom-categories

$$\hom(X, C \underset{A}{\tilde{\times}} B) \longrightarrow \hom(X, C) \underset{\hom(X,A)}{\tilde{\times}} \hom(X, B)$$

is a smothering functor. Here the category on the right is the full subcategory of the comma category $\hom(X, g) \downarrow \hom(X, f)$ spanning those objects whose underlying map in $\hom(X, A)$ is an isomorphism. Explicitly, this weak universal property supplies us with three operations:

(i) (1-cell induction) Given an *iso-comma cone*, an invertible 2-cell $\alpha\colon fb \cong gc$

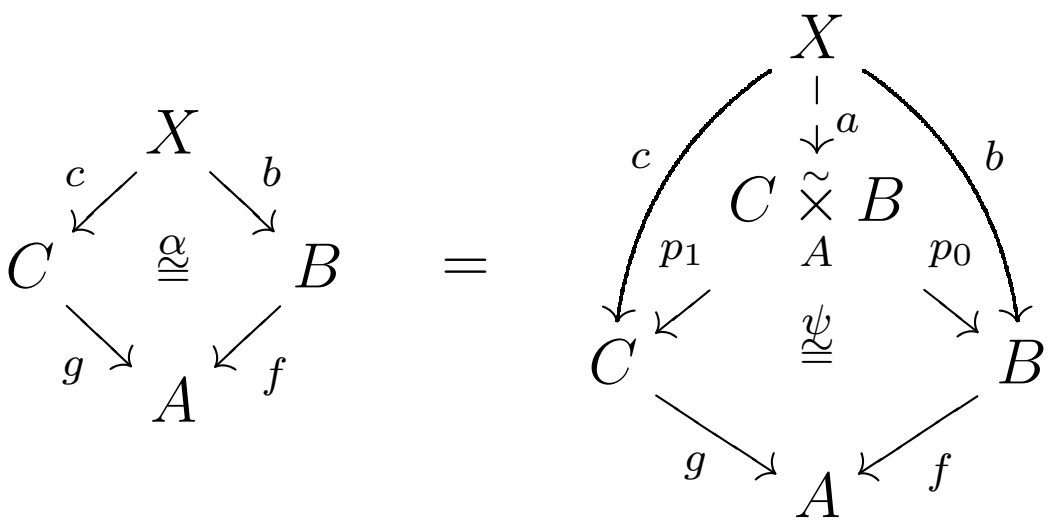

over the pair of functors $f$ and $g$, there exists a 1-cell $a\colon X \to C \tilde{\times}_A B$ so that $p_0 a = b$, $p_1 a = c$, and $\alpha = \psi a$.

(ii) (2-cell induction) Given a pair of functors $a, a'\colon X \to C \tilde{\times}_A B$ and a pair of 2-cells $\tau_0\colon p_0 a \Rightarrow p_0 a'$ and $\tau_1\colon p_1 a \Rightarrow p_1 a'$ which are compatible in the sense that $\psi a' \cdot f\tau_0 = g\tau_1 \cdot \psi a$, then there exists a 2-cell $\tau\colon a \Rightarrow a'$ with $p_0\tau = \tau_0$ and $p_1\tau = \tau_1$.

(iii) (conservativity) Any 2-cell $\tau\colon a \Rightarrow a'\colon X \to C \tilde{\times}_A B$ with the property that the whiskered 2-cells $p_0\tau$ and $p_1\tau$ are both isomorphisms is also an isomorphism.

We refer to (3.2.4) as an *iso-comma square* and $C \xleftarrow{p_1} C \tilde{\times}_A B \xrightarrow{p_0} A$ as an *iso-comma span*.

For definiteness, we have chosen a direction for the 2-cell in an iso-comma cone to be compatible with the direction of the 2-cell in a comma cone, but this choice plays no essential role in the theory of iso-comma squares, which is why it is not indicated in the display (3.2.4).

3.2.5. **Definition** (abstract homotopy 2-category)**.** An *abstract homotopy 2-category* is a strict 2-category $\mathcal{C}$ admitting comma objects and iso-comma objects, in the sense of Definitions 3.2.1 and 3.2.3.

3.3. **Comma and iso-comma objects in a homotopy 2-category.** Our aim in this section is to show that the homotopy 2-category $\mathcal{K}_2$ associated with an $\infty$-cosmos $\mathcal{K}$ is an *abstract homotopy 2-category* in the sense of Definition 3.2.5. This is accomplished by the following pair of lemmas, which demonstrate that $\mathcal{K}_2$ possesses comma objects and iso-comma objects.

3.3.1. **Lemma.** *Given a pair of functors $f\colon B \to A$ and $g\colon C \to A$ in the homotopy 2-category of an $\infty$-cosmos $\mathcal{K}$, their comma object may be constructed by forming the cofibrant replacement of the pullback formed in $\mathcal{K}$:*

$$\begin{array}{ccc} f \downarrow g & \longrightarrow & \Delta^1 \pitchfork A \\ {\scriptstyle (p_1,p_0)}\big\downarrow & \lrcorner & \big\downarrow{\scriptstyle (p_1,p_0)} \\ C \times B & \xrightarrow[g \times f]{} & A \times A \end{array} \tag{3.3.2}$$

*Proof.* The data of the pullback (3.3.2) defines a canonical square

$$\begin{array}{ccccc} & & f \downarrow g & & \\ & {\scriptstyle p_1}\swarrow & & \searrow{\scriptstyle p_0} & \\ C & & \overset{\psi}{\Leftarrow} & & B \\ & {\scriptstyle g}\searrow & & \swarrow{\scriptstyle f} & \\ & & A & & \end{array}$$

in the 2-category $h_*\mathcal{K}$. For any cofibrant object $X$, the proof of Proposition I.3.3.18, while stated in $\underline{\text{qCat}}_2$, applies equally in $h_*\mathcal{K}$ to show that the induced comparison functor of hom-categories

$$\hom(X, f \downarrow g) \longrightarrow \hom(X, f) \downarrow \hom(X, g)$$

is smothering.

It should be noted that our insistence here on $X$ being cofibrant is absolutely necessary. The pertinent point is that the theory developed in §I.3.3 relies upon the assumption that each simplicial hom-space $\text{map}(X, D)$ of maps out of $X$ is a quasi-category. It follows then that each homotopy category $h(\text{map}(X, D)) = \hom(X, D)$ admits a simple explicit description which is exploited repeatedly to make computations in loc. cit. That this supposition holds whenever $X$ is a cofibrant object in an $\infty$-cosmos is simply a consequence of Lemma 2.1.8.

When $f\colon B \to A$ and $g\colon C \to A$ are functors in $\mathcal{K}_2$ we have no reason believe that the object $f \downarrow g$ will also be in there. However, we can take its cofibrant replacement $r_{f\downarrow g}\colon (f \downarrow g)^\sim \xrightarrow{\sim\!\!\!\!\twoheadrightarrow} f \downarrow g$ and use the fact that for any cofibrant $X$ the post-composition functor $\hom(X, r_{f\downarrow g})\colon \hom(X, (f \downarrow g)^\sim) \xrightarrow{\sim\!\!\!\!\twoheadrightarrow} \hom(X, f \downarrow g)$ is a surjective equivalence, as discussed in Lemma 2.1.8, to show that $(f \downarrow g)^\sim$ also satisfies the weak 2-universal property discussed above. □

3.3.3. **Lemma.** *Given a pair of functors $f\colon B \to A$ and $g\colon C \to A$ in the homotopy 2-category of an $\infty$-cosmos $\mathcal{K}$, their iso-comma object may be constructed by forming the* *cofibrant replacement of the pullback formed in $\mathcal{K}$:*

$$\begin{array}{ccc} B \tilde{\times}_A C & \longrightarrow & \mathbb{I} \pitchfork A \\ {\scriptstyle (p_1,p_0)} \downarrow & \lrcorner & \downarrow {\scriptstyle (p_1,p_0)} \\ C \times B & \xrightarrow[g \times f]{} & A \times A \end{array} \tag{3.3.4}$$

*Proof.* By Observation 3.1.6, the data of the pullback (3.3.4) defines a canonical square

$$\begin{array}{ccccc} & & B \tilde{\times}_A C & & \\ & {\scriptstyle p_1} \swarrow & & \searrow {\scriptstyle p_0} & \\ C & & \overset{\psi}{\cong} & & B \\ & {\scriptstyle g} \searrow & & \swarrow {\scriptstyle f} & \\ & & A & & \end{array}$$

in the 2-category $h_*\mathcal{K}$. For any cofibrant object $X$, the functor $\mathrm{map}(X,-)\colon \mathcal{K} \to \underline{\mathrm{qCat}}$ of $\infty$-cosmoi provides a pullback of quasi-categories

$$\begin{array}{ccc} \mathrm{map}(X, B \tilde{\times}_A C) & \longrightarrow & \mathbb{I} \pitchfork \mathrm{map}(X,A) \\ {\scriptstyle (p_1,p_0)} \downarrow & \lrcorner & \downarrow {\scriptstyle (p_1,p_0)} \\ \mathrm{map}(X,C) \times \mathrm{map}(X,B) & \xrightarrow[g \times f]{} & \mathrm{map}(X,A) \times \mathrm{map}(X,A) \end{array}$$

Applying $h\colon \underline{\mathrm{qCat}} \to \underline{\mathrm{Cat}}$, Proposition I.3.3.14 demonstrates that the canonical comparison functor

$$\hom(X, B \underset{A}{\tilde{\times}} C) \longrightarrow (\hom(X,C) \times \hom(X,B)) \underset{\hom(X,A)\times\hom(X,A)}{\times} h(\mathbb{I} \pitchfork \mathrm{map}(X,A))$$

is smothering. Proposition I.3.3.13 tells us that the canonical comparison functor

$$h(\mathbb{I} \pitchfork \mathrm{map}(X,A)) \longrightarrow \hom(X,A)^{\mathbb{I}}$$

is also smothering. A pullback of this defines a smothering functor, which composes with the first smothering functor to demonstrate that the desired functor

$$\hom(X, C \underset{A}{\tilde{\times}} B) \longrightarrow \hom(X,C) \underset{\hom(X,A)}{\tilde{\times}} \hom(X,B)$$

is smothering.

Again, when $f\colon B \to A$ and $g\colon C \to A$ are functors in $\mathcal{K}_2$ we have no reason to believe that the object $C \tilde{\times}_A B$ will also be in there. However, we can take its cofibrant replacement $r_{C\tilde{\times}_A B}\colon (C \tilde{\times}_A B)^\sim \xrightarrow{\sim} C \tilde{\times}_A B$ and use the fact that for any cofibrant $X$ the post-composition functor $\hom(X, r_{C\tilde{\times}_A B})\colon \hom(X, (C \tilde{\times}_A B)^\sim) \xrightarrow{\sim} \hom(X, C \tilde{\times}_A B)$ is a surjective equivalence, as discussed in Lemma 2.1.8, to show that $(C \tilde{\times}_A B)^\sim$ also satisfies the weak 2-universal property discussed above: the composite functor

$$\hom(X, r_{C\tilde{\times}_A B})\colon \hom(X, (C \underset{A}{\tilde{\times}} B)^\sim) \xrightarrow{\sim} \hom(X, C \underset{A}{\tilde{\times}} B) \to \hom(X,C) \underset{\hom(X,A)}{\tilde{\times}} \hom(X,B)$$

remains smothering. □

Immediately from Lemmas 3.3.1 and 3.3.3 we have:

3.3.5. **Corollary.** *The homotopy 2-category of an $\infty$-cosmos is an abstract homotopy 2-category.* □

3.3.6. *Observation.* When working in the homotopy 2-category of an $\infty$-cosmos, we drop extraneous tildes and write simply $f\downarrow g$ and $C \tilde{\times}_A B$ for cofibrant replacements of the simplicial pullbacks defined in (3.3.2) and (3.3.4). Note that the legs of the comma and iso-commas spans produced by these constructions are isofibrations in the $\infty$-cosmos. Lemma 3.4.2 will show that in any abstract homotopy 2-category the legs of a comma span or iso-comma span are always representable isofibrations. By Corollary 3.1.5, isofibrations in the homotopy 2-category $\mathcal{K}_2$ of an $\infty$-cosmos are also representable isofibrations, but the property of being an isofibration in $\mathcal{K}$ is a stronger condition.

3.4. **Stability and uniqueness of comma and iso-comma constructions.** In this section, we develop some of the basic theory of comma and iso-comma constructions in an abstract homotopy 2-category. When we are working in a generic 2-category without a specifically designated class of isofibrations or trivial fibrations then the unqualified terms isofibration and surjective equivalence will be taken to refer to the representably defined concept.

Our first aim is to show that the legs of any comma span or iso-comma span

$$C \xleftarrow{p_1} f \downarrow g \xrightarrow{p_0} B \qquad\qquad C \xleftarrow{p_1} C \underset{A}{\tilde{\times}} B \xrightarrow{p_0} B$$

are isofibrations. In fact more is true: the 2-cell lifts defined with respect to one leg can be chosen to live in the fiber over an identity along the other. We first introduce terminology for this sort of situation.

3.4.1. **Definition.** A *two-sided isofibration* is a span

$$A \xtwoheadleftarrow{q} E \xtwoheadrightarrow{p} B$$

so that

(i) $p\colon E \twoheadrightarrow B$ and $q\colon E \twoheadrightarrow A$ are isofibrations
(ii) The span $(q,p)$ is an *isofibration on the right*: lifts of 2-isomorphisms along $p$ can be chosen to project to identities along $q$.
(iii) The span $(q,p)$ is an *isofibration on the left*: lifts of 2-isomorphisms along $q$ can be chosen to project to identities along $p$.

3.4.2. **Lemma.** *Comma spans $C \xleftarrow{p_1} f \downarrow g \xrightarrow{p_0} B$ and iso-comma spans $C \xleftarrow{p_1} C \tilde{\times}_A B \xrightarrow{p_0} B$ are two-sided isofibrations.*

*Proof.* We prove this for comma spans. The argument for iso-comma spans is analogous, and slightly easier. To show that $(p_1, p_0)$ defines an isofibration on the right consider a

2-cell $\beta\colon b \cong p_0 e$ and form the displayed comma cone and its induced 1-cell:

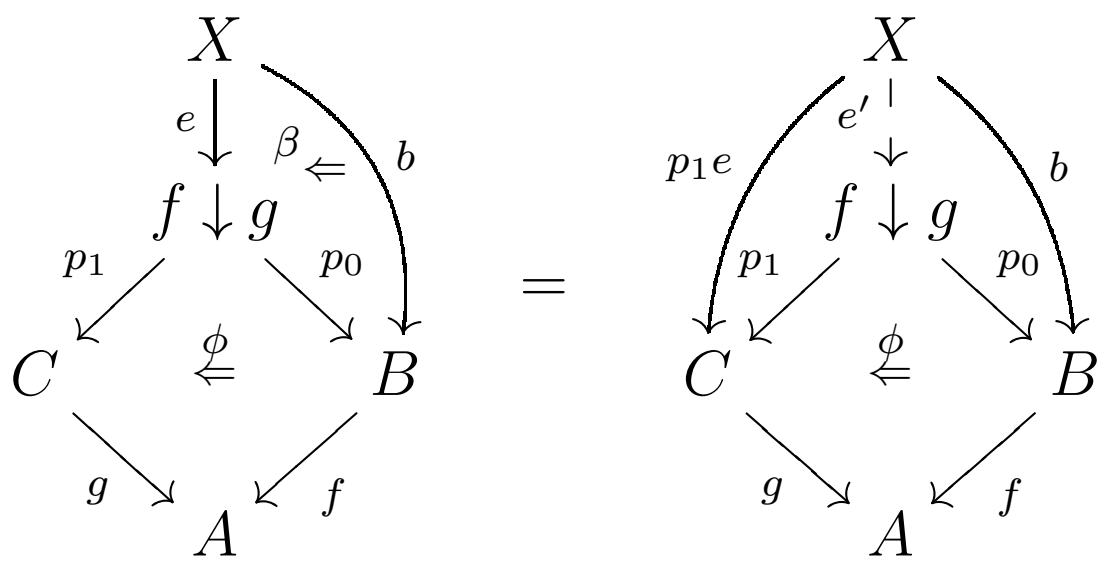

Now, by 2-cell induction and conservativity, there exists a 2-cell isomorphism $\alpha\colon e' \cong e$ defined via the pair $p_0\alpha := \beta\colon p_0 e' \cong p_0 e$ and $p_1\alpha := \mathrm{id}_{p_1 e}\colon p_1 e' \cong p_1 e$. Thus, $\alpha$ is an isomorphism lifting $\beta$ and projecting along $p_1$ to an identity. □

3.4.3. **Definition.** Two-sided isofibrations between a fixed pair of objects form the objects of a strict 2-category $\underline{\mathrm{Span}}_{\mathcal{C}}(A, B)$, the *2-category of spans from $A$ to $B$*, whose:

- objects are two-sided isofibrations $A \overset{p_1}{\twoheadleftarrow} E \overset{p_0}{\twoheadrightarrow} B$,
- 1-cells are maps of spans

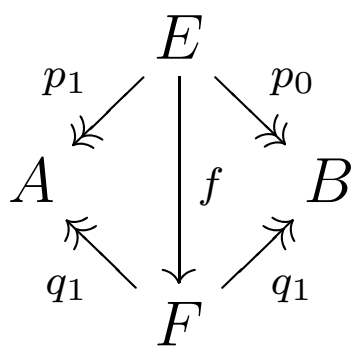

(3.4.4)

  that is to say a 1-cell $f\colon E \to F$ satisfying the equations $q_0 f = p_0$ and $q_1 f = p_1$, and
- 2-cells are 2-cells between maps of spans

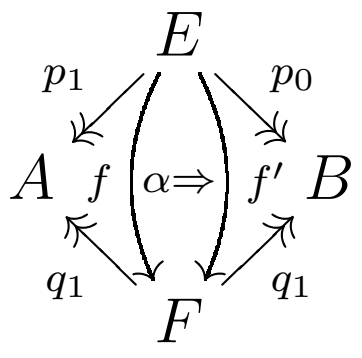

  satisfying the identities $q_0\alpha = \mathrm{id}_{p_0}$ and $q_1\alpha = \mathrm{id}_{p_1}$.

Importantly:

3.4.5. **Lemma.** *A map of spans (3.4.4) is an equivalence in the 2-category $\underline{\mathrm{Span}}_{\mathcal{C}}(A, B)$ if and only if its underlying 1-cell $f\colon E \to F$ defines an equivalence in the 2-category $\mathcal{C}$.*

*Proof.* It is clear that the forgetful 2-functor $\underline{\mathrm{Span}}_{\mathcal{C}}(A, B) \to \mathcal{C}$ preserves equivalences. To demonstrate the converse, choose an adjoint equivalence inverse $g\colon F \to E$ in $\mathcal{C}$ with unit and counit isomorphisms $\eta\colon \mathrm{id}_E \cong gf$ and $\epsilon\colon fg \cong \mathrm{id}_F$. There exist 2-cell isomorphisms:

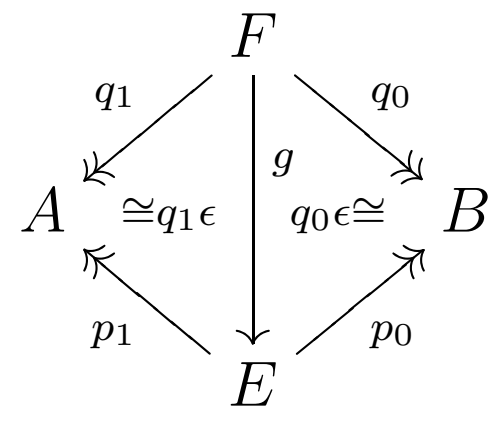

Lifting $q_1\epsilon$ along $p_1$, we obtain a 1-cell $g'\colon F \to E$ and an isomorphism $\alpha_1\colon g \cong g'$ for which $p_1\alpha_1 = q_1\epsilon$ (so $p_1g' = q_1$) and $p_0\alpha_1$ is an identity 2-cell (so $p_0g = p_0q'$). Furthermore, lifting the isomorphism $q_0\epsilon\colon p_0g' = p_0g \cong q_0$ along $p_0$ we obtain a second isomorphism $\alpha_0\colon g' \cong \bar{g}$ for which $p_0\alpha_0 = q_0\epsilon$ (so $p_0\bar{g} = q_0$) and $p_1\alpha_0$ is an identity 2-cell (so $p_1g' = p_1\bar{g}$). It follows, from the equations listed that $p_1\bar{g} = p_1g' = q_1$ and $p_0\bar{g} = q_0$ and consequently that $\bar{g}$ is a map of spans. It is now straightforward to verify that the composite isomorphisms $\mathrm{id}_E \cong gf \cong g'f \cong \bar{g}f$ and $f\bar{g} \cong fg \cong \mathrm{id}_F$ are actually 2-cells in $\underline{\mathrm{Span}}_{\mathcal{C}}(A,B)$. □

Our next aim is to show that the universal properties defining commas and iso-commas characterize unique equivalence classes of objects in the 2-category of spans.

3.4.6. *Observation* (essential uniqueness of induced 1-cells). The 1-cells induced by a cone under these weakly 2-universal properties may not be unique but they are, at least, unique up to isomorphism. For example consider a weakly 2-universal comma object $f \downarrow g$ in a 2-category $\mathcal{C}$ and a pair of 1-cells $\ell, m\colon X \to f \downarrow g$ that both enjoy the same defining properties as 1-cells induced by the weak 2-universal property of $f \downarrow g$, i.e., they satisfy $p_0\ell = p_0m$, $p_1\ell = p_1m$, and $\phi\ell = \phi m$. Then the pair of identities on $p_0\ell = p_0m$ and $p_1\ell = p_1m$ induce a 2-cell $\gamma\colon \ell \Rightarrow m$ and the conservativity property of the comma object $f \downarrow g$ implies that $\gamma$ is an isomorphism. So we have shown that $\ell$ and $m$ are isomorphic via an isomorphism $\gamma$ which is a 2-cell of span maps.

3.4.7. *Recall.* Lemma I.3.3.5 tells us that weak 2-limits in a 2-category $\mathcal{C}$ are unique up to equivalence. Specifically it says that the summits of any two weak 2-limit cones over the same diagram are related by an equivalence which commutes with the legs of those cones.

In our abstract homotopy 2-categories this means that:

- if $C \overset{p_1}{\twoheadleftarrow} f \downarrow g \overset{p_0}{\twoheadrightarrow} B$ and $C \overset{p_1'}{\twoheadleftarrow} (f \downarrow g)' \overset{p_0'}{\twoheadrightarrow} B$ are both comma spans associated with the same pair of 1-cells $f\colon B \to A$ and $g\colon C \to A$ then there exists an equivalence $e\colon f \downarrow g \xrightarrow{\sim} (f \downarrow g)'$ which makes the following triangles

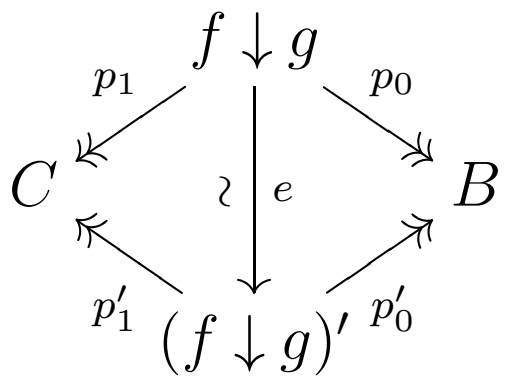

  commute and which defines a factorization of the comma square for $f \downarrow g$ through the comma square for $(f \downarrow g)'$, and

- if $C \overset{p_1}{\twoheadleftarrow} C \tilde{\times}_A B \overset{p_0}{\twoheadrightarrow} B$ and $C \overset{p_1'}{\twoheadleftarrow} (C \tilde{\times}_A B)' \overset{p_0'}{\twoheadrightarrow} B$ are both iso-comma spans associated with the same pair of 1-cells $f\colon B \to A$ and $g\colon C \to A$ then there exists an equivalence $e\colon C \tilde{\times}_A B \xrightarrow{\sim} (C \tilde{\times}_A B)'$ which makes the following triangles

$$\begin{array}{ccccc} & & C \tilde{\times}_A B & & \\ & \swarrow^{p_1} & & \searrow^{p_0} & \\ C & & \wr \downarrow e & & B \\ & \nwarrow_{p_1'} & & \nearrow_{p_0'} & \\ & & (C \tilde{\times}_A B)' & & \end{array}$$

commute and which defines a factorization of the iso-comma square for $C \tilde{\times}_A B$ through the iso-comma square for $(C \tilde{\times}_A B)'$.

Recollection 3.4.7 and Lemma 3.4.5 show that any pair of (iso-)comma spans over the same pair of functors are equivalent in the 2-category of spans. The following lemma proves the converse: that any two-sided isofibration equipped with an equivalence to a comma or iso-comma is again a comma or iso-comma.

3.4.8. **Lemma** (stability of (iso-)comma objects under equivalence)**.** *Suppose that we are given an equivalence*

$$\begin{array}{ccccc} & & E & & \\ & \swarrow^{r_1} & & \searrow^{r_0} & \\ C & & \wr \downarrow e & & B \\ & \nwarrow_{p_1} & & \nearrow_{p_0} & \\ & & f \downarrow g & & \end{array}$$

*in* $\underline{\mathrm{Span}}_{\mathcal{C}}(C, B)$*, where* $f \downarrow g$ *is a comma object displayed by the data in* (3.2.2)*. Then the square*

$$\begin{array}{ccccc} & & E & & \\ & \swarrow^{r_1} & & \searrow^{r_0} & \\ C & & \Downarrow \psi e & & B \\ & \searrow_{g} & & \swarrow_{f} & \\ & & A & & \end{array}$$

*is a comma square. The same result is true, mutatis mutandis, for iso-comma spans.*

Lemma 3.4.5 implies that the direction of the given equivalence in Lemma 3.4.8 is immaterial.

*Proof.* For any $X$, the canonical functor

$$\hom(X, E) \to \hom(X, f \downarrow g) \to \hom(X, f) \downarrow \hom(X, g)$$

induced by the square $\psi e$ is the composite of an equivalence with a smothering functor, and as such is immediately full and conservative. It remains only to show that the composite, which is clearly essentially surjective on objects, is in fact surjective on objects.

To this end, observe that any object in $\hom(X,f)\downarrow\hom(X,g)$ has a preimage in $\hom(X,f\downarrow g)$. By Lemma 3.4.5, this object is isomorphic, via some isomorphism projecting to an identity in $\hom(X,C)\times\hom(X,B)$ to an object in the image of $\hom(X,E)\to\hom(X,f\downarrow g)$. This is where we make use of the hypothesis that $E$ is a two-sided isofibration. But any pair of objects in $\hom(X,f\downarrow g)$, which are isomorphic over an identity in $\hom(X,C)\times\hom(X,B)$, have the same image in $\hom(X,f)\downarrow\hom(X,g)$ by Observation 3.4.6. Thus $\hom(X,E)\to\hom(X,f)\downarrow\hom(X,g)$ is surjective on objects, as desired. □

In summary:

3.4.9. **Corollary.** *Given a fixed cospan $C\xrightarrow{g}A\xleftarrow{f}B$ of 1-cells in $\mathcal{C}$*

*(i) Any two comma spans over $f$ and $g$ are equivalent as objects in the 2-category of spans from $C$ to $B$.*

*(ii) Any two iso-comma spans $f$ and $g$ are equivalent as objects in the 2-category of spans from $C$ to $B$.*

*And conversely:*

*(iii) Any two-sided isofibration that is equivalent to a comma span over $f$ and $g$ is again a comma span over that pair.*

*(iv) Any two-sided isofibration that is equivalent to an iso-comma span over $f$ and $g$ is again an iso-comma span over that pair.* □

We now demonstrate that comma squares and iso-comma squares are stable under composition with an iso-comma square on either the left or the top.

3.4.10. **Lemma.** *Consider a diagram*

$$\begin{array}{ccccc} F & \xrightarrow{r} & D & \xrightarrow{p} & B \\ {\scriptstyle s}\downarrow & \cong\psi & {\scriptstyle q}\downarrow & \Leftarrow\phi & \downarrow{\scriptstyle f} \\ E & \xrightarrow[h]{} & C & \xrightarrow[g]{} & A \end{array}$$

*in which $\psi\colon qr\cong hs$ is an iso-comma square.*

*(i) If $\phi\colon fp\Rightarrow gq$ is a comma square, then so is the composite rectangle.*

*(ii) If $\phi\colon fg\cong gq$ is an iso-comma square, then so is the composite rectangle.*

*A dual result holds with the direction of the commas reversed:*

*(iii) If $\phi\colon gq\Rightarrow fp$ is a comma square, then so is the composite rectangle.*

*Proof.* The proofs of (ii) and (iii) are similar to the proof of (i), which we give here. Suppose $\phi$ defines a comma square; the argument for iso-comma squares is parallel. A comma cone over $f$ and $hg$ may be factored through the comma square $\phi$ and then the identity 2-cell $qd = he$ may be factored through the iso-comma square $\psi$. This defines 1-cell induction for the composite rectangle.

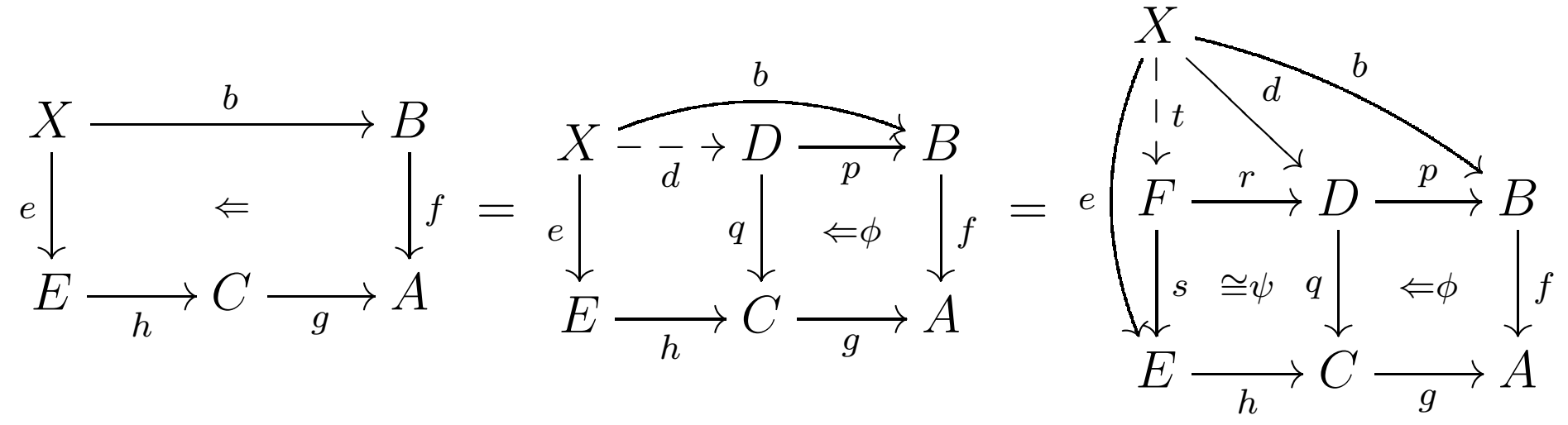

Now consider a parallel pair $t, t' \colon X \to F$ equipped with 2-cells $\beta \colon prt \Rightarrow prt'$ and $\epsilon \colon st \Rightarrow st'$ so that

$$\begin{array}{ccccc} fprt & \xRightarrow{\phi rt} & gqrt & \xRightarrow[\cong]{g\psi t} & ghst \\ {\scriptstyle f\beta}\Downarrow & & & & \Downarrow{\scriptstyle gh\epsilon} \\ fprt' & \xRightarrow[\phi rt']{} & gqrt' & \xRightarrow[g\psi t']{\cong} & ghst' \end{array}$$

commutes. The 2-cell induction property of $\phi$ applies to $\beta \colon prt \Rightarrow prt'$ and the composite $qrt \overset{g\psi t}{\Rightarrow} ghst \overset{gh\epsilon}{\Rightarrow} ghst' \overset{g\psi^{-1}t'}{\Rightarrow} gqrt'$ to induce a 2-cell $\chi \colon rt \Rightarrow rt'$ so that $p\chi = \beta$ and the diagram

$$\begin{array}{ccc} gqrt & \xRightarrow[\cong]{g\psi t} & ghst \\ {\scriptstyle gq\chi}\Downarrow & & \Downarrow{\scriptstyle gh\epsilon} \\ gqrt' & \xRightarrow[g\psi t']{\cong} & ghst' \end{array}$$

commutes. Now the 2-cell induction property of $\psi$ applies to $\chi \colon rt \Rightarrow rt'$ and $\epsilon \colon st \Rightarrow st'$ to induce a 2-cell $\tau \colon t \Rightarrow t'$ so that $s\gamma = \epsilon$ and $pr\gamma = \beta$. This defines 2-cell induction for the composite rectangle.

If $\beta$ and $\epsilon$ are isomorphisms, then so is $\chi$ and hence so is $\gamma$, proving 2-cell conservativity. □

Lemma 3.4.10 has an accompanying cancelation result that can be used to detect iso-comma squares.

3.4.11. **Lemma.** *Consider a diagram*

$$\begin{array}{ccccc} F & \xrightarrow{r} & D & \xrightarrow{p} & B \\ {\scriptstyle s}\downarrow & \cong\psi & {\scriptstyle q}\downarrow & \Leftarrow\phi & \downarrow{\scriptstyle f} \\ E & \xrightarrow[h]{} & C & \xrightarrow[g]{} & A \end{array}$$

*in which* $E \overset{s}{\twoheadleftarrow} F \overset{r}{\twoheadrightarrow} D$ *is a two-sided isofibration.*

*(i) If $\phi\colon fp \Rightarrow gq$ and the composite rectangle are comma squares, then $\psi\colon qr \cong hs$ is an iso-comma square.*
*(ii) If $\phi\colon fg \cong gq$ and the composite rectangle are iso-comma squares, then $\psi\colon qr \cong hs$ is an iso-comma square.*

*Proof.* Applying 1-cell induction to the iso-comma cone $\psi$, there is a map $e\colon F \to E \mathbin{\tilde{\times}_C} D$ in $\underline{\mathrm{Span}}_{\mathcal{C}}(E, D)$ from $F$ to the iso-comma span

$$\begin{array}{ccccc} & & E \mathbin{\tilde{\times}_C} D & & \\ & {\scriptstyle p_1}\swarrow & & \searrow{\scriptstyle p_0} & \\ E & & \overset{\psi'}{\Leftarrow} & & D \\ & {\scriptstyle h}\searrow & & \swarrow{\scriptstyle q} & \\ & & C & & \end{array}$$

By hypothesis (i) or (ii) and the parallel result from Lemma 3.4.10, both the composite rectangles $g\psi \cdot \phi r$ and $g\psi' \cdot \phi r$ are (iso-)commas over the cospan $E \xrightarrow{gh} A \xleftarrow{f} B$. Recollection 3.4.7 implies that the map $e$ is an equivalence in $\underline{\mathrm{Span}}_{\mathcal{C}}(E, B)$. As this 1-cell lies in $\underline{\mathrm{Span}}_{\mathcal{C}}(E, D)$, Lemma 3.4.5 implies that it also defines an equivalence of two-sided isofibrations from $E$ to $D$. It follows from Lemma 3.4.8 that $\psi\colon hs \cong qr$ is then an iso-comma square. □

From these results we obtain the following "pullback stability" result for commas. The connection with pullbacks will be explained in the next section.

3.4.12. **Lemma.** *Suppose that we are given 1-cells $f\colon B \to A$, $g\colon C \to A$, $h\colon B' \to B$, and $k\colon C' \to C$ in an abstract homotopy 2-category $\mathcal{C}$. Then there is a diagram*

$$\begin{array}{ccccc} fh \downarrow gk & \twoheadrightarrow & fh \downarrow g & \xrightarrow{p_0} & B' \\ \downarrow & \cong & \downarrow & \cong & \downarrow{\scriptstyle h} \\ f \downarrow gk & \twoheadrightarrow & f \downarrow g & \xrightarrow{p_0} & B \\ {\scriptstyle p_1}\downarrow & \cong & {\scriptstyle p_1}\downarrow & \Leftarrow\phi & \downarrow{\scriptstyle f} \\ C' & \xrightarrow[k]{} & C & \xrightarrow[g]{} & A \end{array}$$

*in which*
*(i) The lower-right square, the bottom rectangle, the right rectangle, and the outer square are comma squares.*
*(ii) The top-left square, top-right square, the lower-left square, the top rectangle, and the left rectangle are iso-comma squares.*

*Proof.* Form the lower-right comma square and then the three iso-comma squares, ending with the top-left iso-comma square. By Lemma 3.4.10, the bottom and left rectangles are again comma squares. By Corollary 3.4.9 this implies that their summits define the comma objects $f \downarrow gk$ and $fh \downarrow g$. By Lemma 3.4.10, the top rectangle is also an iso-comma square; applying this result to the composite of the top rectangle with the bottom rectangle tells us that the outer composite square is a comma square. By Corollary 3.4.9, this implies that the summit is the comma $fh \downarrow gk$ as claimed. □

3.5. **Iso-commas and pullbacks.** We will frequently work in slices $\mathcal{C}/B$ of an abstract homotopy 2-category over a fixed object $B$. These slice 2-categories are not themselves abstract homotopy 2-categories. However, the ambient abstract homotopy 2-category $\mathcal{C}$ will supply an important operation of pulling back from one slice to another that, while not 2-functorial, will have a number of pleasing properties as we will soon discover.

3.5.1. *Observation* (slice 2-categories). The enriched slice construction of Example 2.1.11 applies in any category theory which is enriched over a cartesian closed category. In particular, it applies to 2-categories since they are neither more nor less than categories enriched in the cartesian closed category of all (small) categories. If $B$ in an object in a 2-category $\mathcal{C}$ then the slice $\mathcal{C}/B$ has objects that are 1-cells $p\colon E \to B$ in $\mathcal{C}$, 1-cells that are commuting triangles of the form

$$\begin{array}{ccc} E & \xrightarrow{f} & F \\ & {\scriptstyle p}\searrow \quad \swarrow{\scriptstyle q} & \\ & B & \end{array} \tag{3.5.2}$$

in $\mathcal{C}$ and 2-cells those $\alpha\colon f \Rightarrow g\colon E \to F$ of $\mathcal{C}$ with the property that $q\alpha$ is the identity 2-cell on $p$. We will use the notation $\hom_B(p, q)$ to denote the hom-category of $\mathcal{C}/B$ of 1-cells and 2-cells between objects $p\colon E \to B$ and $q\colon F \to B$.

When working with abstract homotopy 2-categories $\mathcal{C}$ our default position will be to restrict our attention to the full sub-2-category of the slice 2-category spanned by those objects $p\colon E \twoheadrightarrow B$ that are isofibrations. To avoid cluttering our notation in that context we simply write $\mathcal{C}/B$ to denote this restricted slice.

We often use the term *fibred* to refer to concepts in a slice $\mathcal{C}/B$ of an abstract homotopy 2-category. For instance, specializing Lemma 3.4.5 to the case where one of the objects is terminal (or absent), we have:

3.5.3. **Lemma.** *A 1-cell of $\mathcal{C}/B$ is a* fibred equivalence *if and only if it is an equivalence in $\mathcal{C}$.* □

We are now ready to describe the pullback operation. For motivation, suppose

$$\begin{array}{ccc} F & \xrightarrow{g} & E \\ {\scriptstyle q}\big\downarrow & \cong\psi & \big\downarrow{\scriptstyle p} \\ A & \xrightarrow[f]{} & B \end{array}$$

is an iso-comma square and $p\colon E \twoheadrightarrow B$ is an isofibration. Then the isomorphism $\psi\colon pg \cong fq$ can be lifted along $p$ to define a new functor $\bar{g}\colon F \to E$ isomorphic to $g$

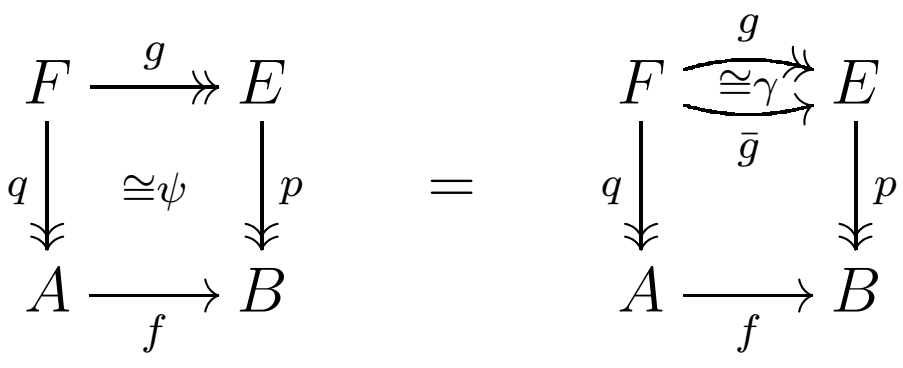

so that $p\gamma = \psi$ and the square $p\bar{g} = fq$ commutes.

The following lemma demonstrates that the commutative square so-obtained functions as a kind of weak pullback. These are more general than the weak pullbacks defined in §I.3.3; see Example 3.5.9.

3.5.4. **Definition** (pullback)**.** A commutative square

$$\begin{array}{ccc} F & \xrightarrow{g} & E \\ {\scriptstyle q}\downarrow & \lrcorner & \downarrow{\scriptstyle p} \\ A & \xrightarrow[f]{} & B \end{array}$$

whose verticals are isofibrations is a *pullback* if

(i) (1-cell induction) Given any commutative square

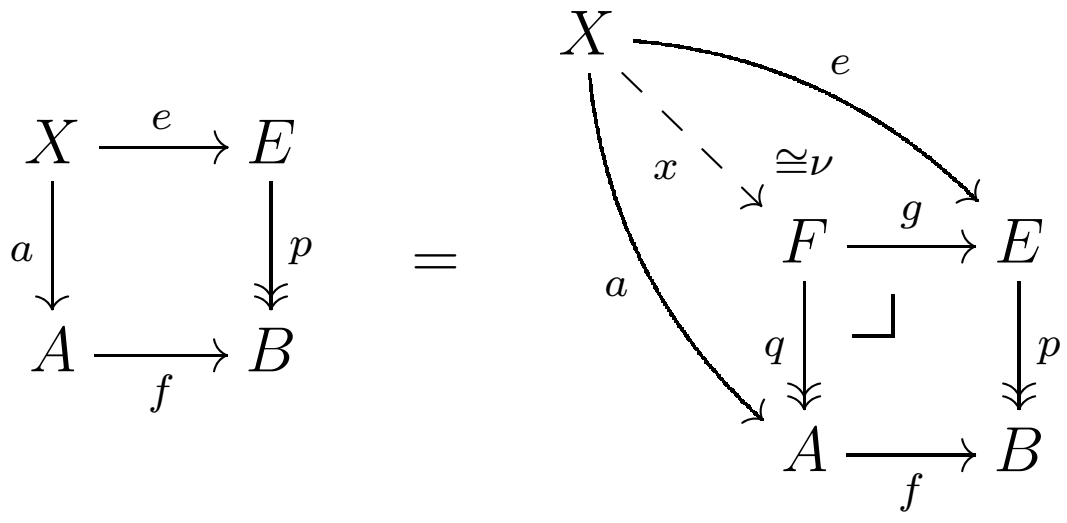

over $f$ and $p$ there exists a morphism $x\colon X \to F$ and an isomorphism $\nu\colon e \cong gx$ so that $qx = a$ and $p\nu$ is an identity.

(ii) (2-cell induction) Suppose we are given 1-cells $x, x'\colon X \to F$. Then for any pair of 2-cells $\epsilon\colon gx \Rightarrow gx'$ and $\alpha\colon qx \Rightarrow qx'$ such that $p\epsilon = f\alpha$, there exists a 2-cell $\tau\colon x \Rightarrow x'$ with $g\tau = \epsilon$ and $q\tau = \alpha$.

(iii) (conservativity) Any 2-cell $\tau\colon x \Rightarrow x'\colon X \to F$ with the property that the whiskered 2-cells $q\tau$ and $g\tau$ are both isomorphisms is also an isomorphism.

3.5.5. *Remark.* The argument of Observation 3.4.6 shows that induced 1-cells $x\colon X \to F$ are unique up to an invertible 2-cell projecting along $q\colon F \twoheadrightarrow A$ to an identity.

In fact, because the right-hand vertical map in a pullback is required to be an isofibration, pullbacks satisfy a more general 1-cell induction property for iso squares. As we will not make explicit use of this here, we leave the details to the reader.

3.5.6. **Lemma.** *An abstract homotopy 2-category $\mathcal{C}$ has pullbacks of isofibrations constructed by forming the iso-comma and then lifting away the isomorphism:*

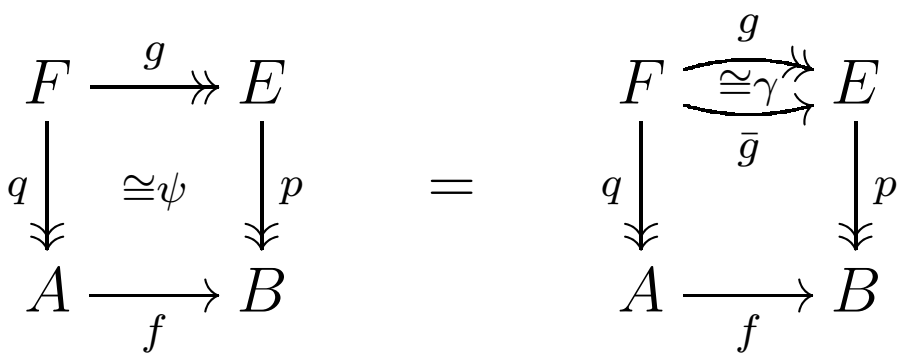

*Proof.* It remains to show that the square $p\bar{g} = fq$ is a pullback. For 1-cell induction, a commutative square over $f$ and $p$ is a special case of an iso-comma cone, so there exists a factorization which gives the desired map $x$ and isomorphism $\gamma x\colon e = gx \cong \bar{g}x$ so that $p\gamma x = \psi x$ is an identity.

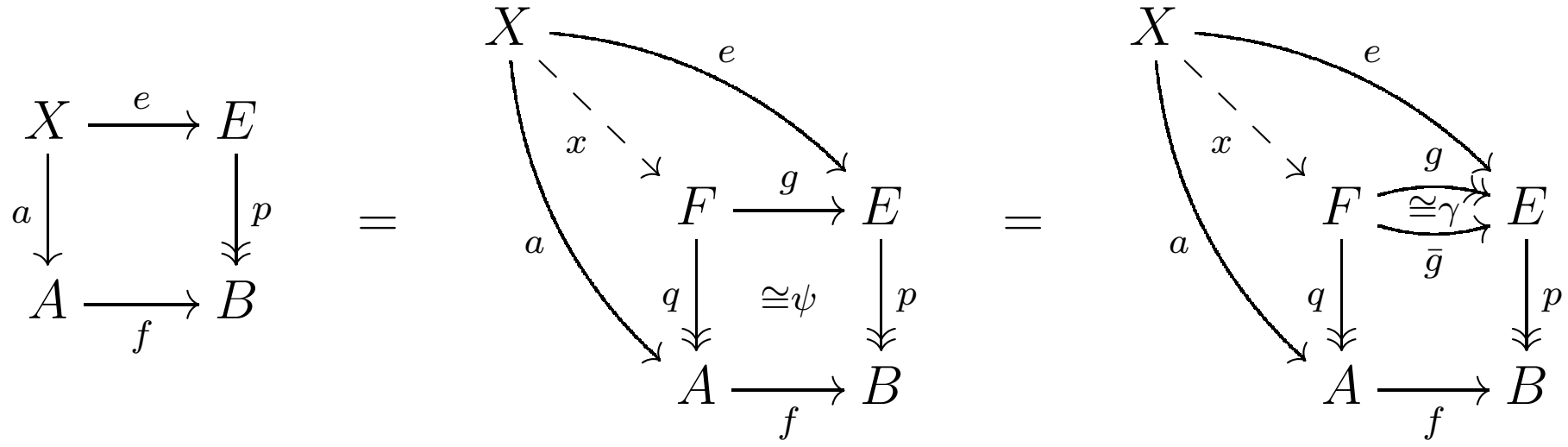

Now fix a pair of functors $x, x'\colon X \to F$. Given a pair of 2-cells $\epsilon\colon \bar{g}x \Rightarrow \bar{g}x'$ and $\alpha\colon qx \Rightarrow qx'$ so that $p\epsilon = f\alpha$, the pair

$$qx \xRightarrow{\alpha} qx' \qquad gx \xRightarrow{\gamma x} \bar{g}x \xRightarrow{\epsilon} \bar{g}x' \xRightarrow{\gamma^{-1}x'} gx'$$

satisfy the compatibility condition necessary to induce a 2-cell $\tau\colon x \Rightarrow x'$ so that $q\tau = \alpha$ and $g\tau\colon gx \Rightarrow gx'$ is the displayed composite. By middle-four interchange,

$$\bar{g}\tau = \gamma x' \cdot g\tau \cdot \gamma^{-1}x = \gamma x' \cdot \gamma^{-1}x' \cdot \epsilon \cdot \gamma x \cdot \gamma^{-1}x = \epsilon,$$

as desired. If $\alpha$ and $\epsilon$ are isomorphisms, then so is the inducing pair, and thus $\tau$ is an isomorphism, by 2-cell conservativity for the iso-comma $\psi$. □

3.5.7. **Definition.** We say that an isofibration $q\colon F \twoheadrightarrow A$ is a *pullback* of an isofibration $p\colon E \twoheadrightarrow B$ along a functor $f\colon A \to B$ if there exists a pullback square

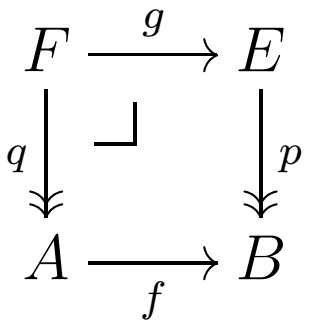

3.5.8. **Example.** If there exists an iso-comma square

$$\begin{array}{ccc} F & \xrightarrow{g} & E \\ {\scriptstyle q}\downarrow & \cong\psi & \downarrow{\scriptstyle p} \\ A & \xrightarrow[f]{} & B \end{array}$$

then $q\colon F \twoheadrightarrow A$ is a pullback of $p\colon E \twoheadrightarrow B$ along $f\colon A \to B$. For instance, the iso-commas of Lemma 3.4.12 exhibit a number of pullback squares involving the functors appearing in comma spans, except that the top horizontals will be replaced by isomorphic functors that are not necessarily isofibrations.

3.5.9. **Example.** If

$$\begin{array}{ccc} F & \xrightarrow{g} & E \\ {\scriptstyle q}\downarrow & \lrcorner & \downarrow{\scriptstyle p} \\ A & \xrightarrow[f]{} & B \end{array}$$

is a commutative square so that for all $X$ the induced functor

$$\hom(X,F) \to \hom(X,E) \underset{\hom(X,A)}{\times} \hom(X,B)$$

is smothering, then $F$ is a pullback in the sense of Definition 3.5.4. On account of the smothering functor, in this case the induced 1-cells into $F$ can be chosen so that both triangles commute strictly.

An argument very similar to the proof of Lemma 3.4.10 proves the following result:

3.5.10. **Lemma** (composition of pullbacks). *Suppose that we are given a pair of squares*

$$\begin{array}{ccccc} E'' & \xrightarrow{g'} & E' & \xrightarrow{g} & E \\ {\scriptstyle p''}\downarrow & \lrcorner & {\scriptstyle p'}\downarrow & \lrcorner & \downarrow{\scriptstyle p} \\ B'' & \xrightarrow[f']{} & B' & \xrightarrow[f]{} & B \end{array} \qquad \textit{respectively} \qquad \begin{array}{ccc} F' & \xrightarrow{h} & F \\ {\scriptstyle q'}\downarrow & \lrcorner & \downarrow{\scriptstyle q} \\ E' & \xrightarrow{g} & E \\ {\scriptstyle p'}\downarrow & \lrcorner & \downarrow{\scriptstyle p} \\ B' & \xrightarrow[f]{} & B \end{array}$$

*in a 2-category $\mathcal{C}$. Then if both squares in the diagram are pullbacks in the sense of 3.5.4 then so is the composite outer rectangle.*

3.5.11. *Remark.* Pullbacks also cancel in the sense of Lemma 3.4.11. The argument, which is somewhat more subtle, is omitted because we will not require this result here.

An easy argument along the lines of Recollection 3.4.7 demonstrates that pullbacks are well-defined up to equivalence. In fact more is true: equivalent isofibrations over $B$ pull back along $f\colon A \to B$ to equivalent isofibrations over $A$, as we shall discover in Corollary 3.6.7. This result and another, equally of interest, will follow from general considerations, which we now discuss.

3.6. **Smothering 2-functors and adjunctions.** An adjunction in an $\infty$-cosmos $\mathcal{K}$ is simply defined to be an adjunction in the associated homotopy 2-category $\mathcal{K}_2$. In the sequel we will need to be able to pull adjunctions in a slice $\mathcal{C}/B$ back along an arbitrary 1-cell $f\colon A \to B$ to give adjunctions in $\mathcal{C}/A$. However, the weak 2-universal property of weak pullbacks in an abstract homotopy 2-category is not strong enough to ensure that we can define a pullback 2-functor $f^*\colon \mathcal{C}/B \to \mathcal{C}/A$. As in §I.3.4, we are able to circumvent this lack of 2-functoriality by making use of a suitable smothering 2-functor.

3.6.1. **Definition** (smothering 2-functors). A *smothering 2-functor* is a 2-functor $F\colon \mathcal{C} \to \mathcal{D}$ that is surjective on objects and *locally smothering*: each of its actions $F\colon \mathcal{C}(A,B) \to \mathcal{D}(FA, FB)$ on hom-categories is a smothering functor.

Smothering 2-functors are *conservative on 1-cells* — a 1-cell in $F$ is an equivalence if and only if its image in $\mathcal{C}$ is an equivalence — and also reflect equivalence between objects. Moreover:

I.4.5.2. **Lemma.** *Suppose $F\colon \mathcal{C} \to \mathcal{D}$ is a smothering 2-functor. Then any adjunction in $\mathcal{D}$ can be lifted to an adjunction in $\mathcal{C}$. Furthermore, if we have previously specified a lift of the objects, 1-cells, and either the unit or counit of the adjunction in $\mathcal{D}$, then there is a lift of the remaining 2-cell that combines with the previously specified data to define an adjunction in $\mathcal{C}$.* □

3.6.2. **Definition.** Given an abstract homotopy 2-category $\mathcal{C}$ and a functor $f\colon A \to B$, define a 2-category $\mathrm{pbk}(\mathcal{C}, f)$ whose:

- objects are pullback squares

$$\begin{array}{ccc} F & \xrightarrow{g} & E \\ {\scriptstyle q}\downarrow & \lrcorner & \downarrow{\scriptstyle p} \\ A & \xrightarrow[f]{} & B \end{array}$$

- 1-cells consist of triples $(k, h, \gamma)$ as depicted in the following commutative diagram:

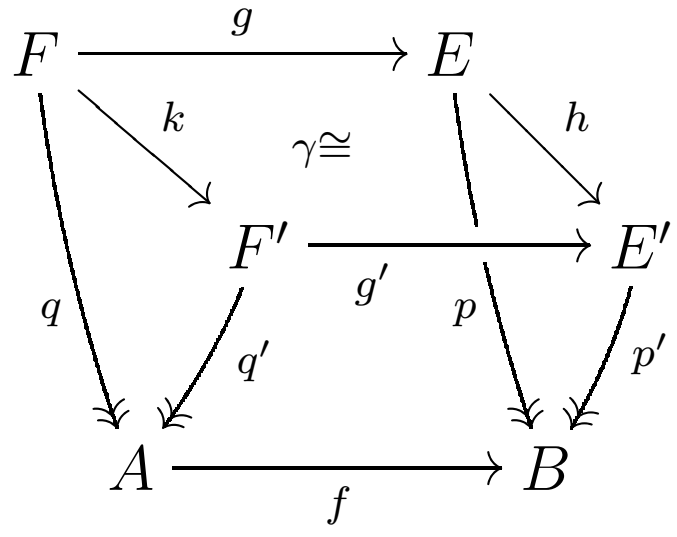

(3.6.3)

  so that $k$ and $h$ are cells in the slices over $A$ and $B$ and $\gamma\colon hg \cong g'k$ is such that $p'\gamma$ is an identity.

- 2-cells consist of pairs of 2-cells $(\beta, \alpha)$ as depicted in the following commutative diagram:

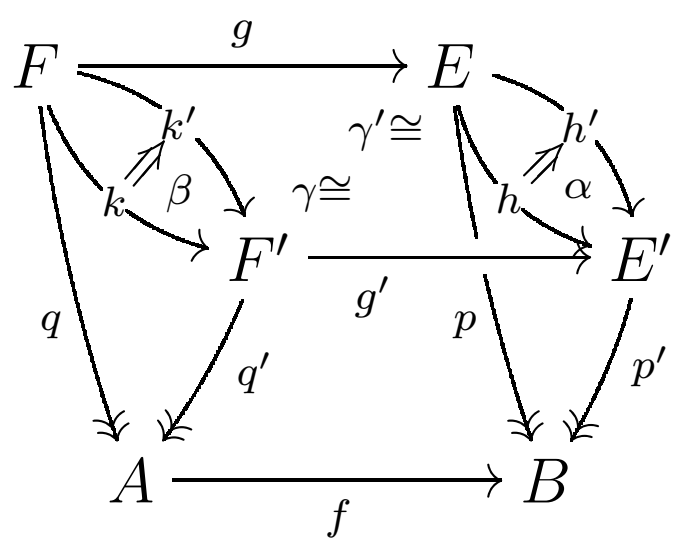

(3.6.4)

Explicitly this means that $\beta$ and $\alpha$ satisfy the equalities required of cells in the slices over $A$ and $B$ and the diagram

$$\begin{array}{ccc} g'k & \xRightarrow{g'\beta} & g'k' \\ \gamma \Downarrow \cong & & \cong \Downarrow \gamma' \\ hg & \xRightarrow[\alpha g]{} & h'g \end{array}$$

commutes.

These cells compose pointwise to make $\mathrm{pbk}(\mathcal{C}, f)$ into a 2-category which admits obvious 2-functorial projections $P_0 \colon \mathrm{pbk}(\mathcal{C}, f) \to \mathcal{C}/A$ and $P_1 \colon \mathrm{pbk}(\mathcal{C}, f) \to \mathcal{C}/B$.

3.6.5. **Lemma.** *The projection* $P_1 \colon \mathrm{pbk}(\mathcal{C}, f) \to \mathcal{C}/B$ *is a smothering 2-functor.*

*Proof.* By Lemma 3.5.6, $P_1 \colon \mathrm{pbk}(\mathcal{C}, f) \to \mathcal{C}/B$ is surjective on objects. From 1-cell induction, it follows that $P_1 \colon \mathrm{pbk}(\mathcal{C}, f) \to \mathcal{C}/B$ is locally surjective on 1-cells: given commutative squares $pg = fq$ and $p'g' = fq'$ as in (3.6.3), the condition $p'h = p$ implies that the former defines a cone over the latter. The 1-cell $k \colon F \to F'$ and isomorphism $\gamma \colon hg \cong g'k$ so that $p'\gamma$ is an identity are induced by the front pullback, defining the desired lift in (3.6.3).

Similarly, using 2-cell induction, it follows that $P_1 \colon \mathrm{pbk}(\mathcal{C}, f) \to \mathcal{C}/B$ is locally full. Given pairs of parallel 1-cells $(k, h, \gamma)$ and $(k', h', \gamma')$ as in (3.6.4) and a 2-cell $\alpha \colon h \Rightarrow h'$ in $\mathcal{C}/B$ so that $p'\alpha = \mathrm{id}_p$, this 2-cell together with $\mathrm{id}_q$ induces the desired lift $\beta \colon k \Rightarrow k'$ displayed in (3.6.4). If $\alpha$ is an isomorphism, then 2-cell conservativity for the pullback $F'$ implies immediately that $\beta$ is also, completing the proof. □

3.6.6. **Corollary.** *Suppose given an adjunction in* $\mathcal{C}/B$ *and a functor* $f \colon A \to B$*. Then there is an adjunction in* $\mathcal{C}/A$ *given by pulling the adjunction in* $\mathcal{C}/B$ *along* $f$*.*

*Proof.* Combining Lemmas 3.6.5 and I.4.5.2, the adjunction in $\mathcal{C}/B$ lifts along the smothering 2-functor $P_1 \colon \mathrm{pbk}(\mathcal{C}, f) \to \mathcal{C}/B$ to give an adjunction in $\mathrm{pbk}(\mathcal{C}, f)$. Applying the other projection $P_0 \colon \mathrm{pbk}(\mathcal{C}, f) \to \mathcal{C}/A$ we get the desired adjunction between pullbacks in $\mathcal{C}/A$. □

3.6.7. **Corollary.** *If* $p \colon E \twoheadrightarrow B$ *and* $p' \colon E' \twoheadrightarrow B'$ *are equivalent in* $\mathcal{C}/B$ *then their pullbacks along any functor* $f \colon A \to B$ *are equivalent in* $\mathcal{C}/A$*.*

*Proof.* The equivalence lifts along the smothering 2-functor $P_1\colon \mathrm{pbk}(\mathcal{C}, f) \to \mathcal{C}/B$ to give an equivalence in $\mathrm{pbk}(\mathcal{C}, f)$. Applying the other projection $P_0\colon \mathrm{pbk}(\mathcal{C}, f) \to \mathcal{C}/A$ we get the desired equivalence between pullbacks in $\mathcal{C}/A$. □

Similar results holds with iso-commas playing the role of pullbacks and without requiring that the initial adjunction is fibred.

3.6.8. **Definition.** Given an abstract homotopy 2-category $\mathcal{C}$ and a functor $f\colon A \to B$, define a 2-category $\mathrm{icom}(\mathcal{C}, f)$ whose:

- objects are iso-comma squares

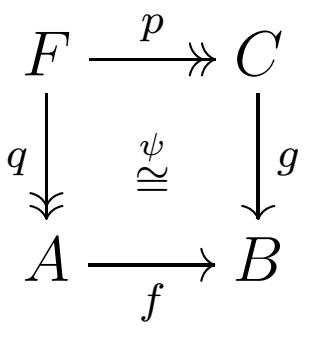

- 1-cells consist of triples $(k, h, \gamma\colon g \cong g'h)$ as depicted in the following commutative diagram:

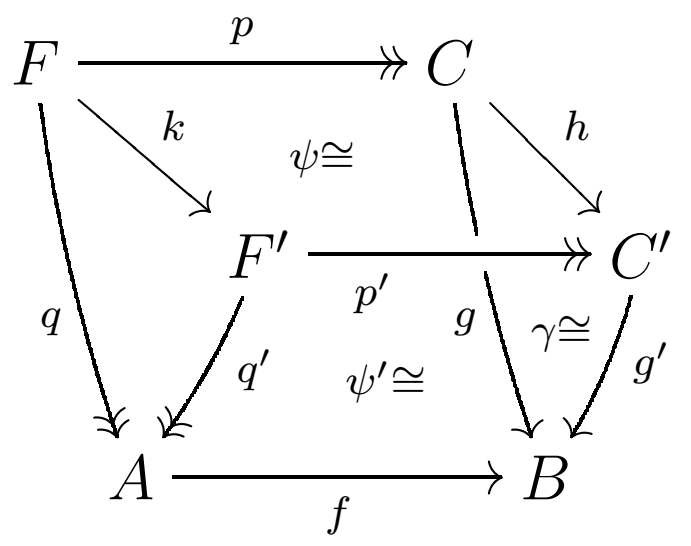

(3.6.9)

so that $k$ is a cell in the slice over $A$, $p'k = hp$, and $\psi' k = g'\psi \cdot \gamma p$.

- 2-cells consist of pairs of 2-cells $(\beta, \alpha)$ as depicted in the following commutative diagram:

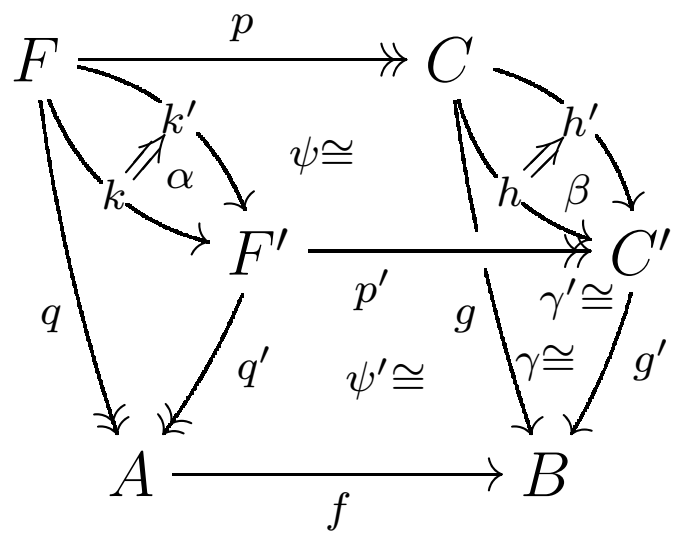

(3.6.10)

so that $\alpha$ is a 2-cell in the slice over $A$, $\gamma' \cdot g'\beta = \gamma$, and $p'\alpha = \beta p$.

These cells compose pointwise to make $\mathrm{icom}(\mathcal{C}, f)$ into a 2-category which admits obvious 2-functorial projections $P_0\colon \mathrm{icom}(\mathcal{C}, f) \to \mathcal{C}/A$ and $P_1\colon \mathrm{icom}(\mathcal{C}, f) \to \mathcal{C}/_{\cong}B$, where $\mathcal{C}/_{\cong}B$ is the *pseudo slice 2-category* over $B$, whose objects are arbitrary maps with codomain $B$, whose morphisms are triangles commuting up to a specified isomorphism, and whose 2-cells are 2-cells between the initial legs of such triangles commuting with the isomorphisms.

3.6.11. **Lemma.** *The projection* $P_1\colon \mathrm{icom}(\mathcal{C}, f) \to \mathcal{C}/_{\cong}B$ *is a smothering 2-functor.*

*Proof.* Existence of iso-commas implies that $P_1$ is surjective on objects. The proof that $P_1$ is locally smothering on the hom category from an isocomma $\psi\colon gp \cong fq$ to an isocomma $\psi'\colon g'p' \cong fq'$ follows easily from 1-cell induction, 2-cell induction, and 2-cell conservativity for $\psi'$. The 1-cell $k$ in (3.6.9) is induced from the iso-comma cone $g'\psi \cdot \gamma p$. The 2-cell $\alpha$ in (3.6.10) is induced from $\mathrm{id}_q$ and $\beta p$, which satisfy the required compatibility condition on account of the equation $\gamma' \cdot g'\beta = \gamma$. 2-cell conservativity is immediate. □

3.6.12. **Corollary.** *Suppose given an adjunction in $\mathcal{C}/_{\cong}B$ and a functor $f\colon A \to B$. Then there is an adjunction in $\mathcal{C}/A$ between the opposing legs of the iso-commas formed from these maps and $f$.*

*Proof.* Combining Lemmas 3.6.11 and I.4.5.2, the adjunction in $\mathcal{C}/_{\cong}B$ lifts along the smothering 2-functor $P_1\colon \mathrm{icom}(\mathcal{C}, f) \to \mathcal{C}/_{\cong}B$ to give an adjunction in $\mathrm{icom}(\mathcal{C}, f)$. Applying the other projection $P_0\colon \mathrm{icom}(\mathcal{C}, f) \to \mathcal{C}/A$ we get the desired adjunction in $\mathcal{C}/A$. □

We conclude by recalling a few basic facts about adjunctions in any 2-category.

3.6.13. *Observation* (adjunctions are representably defined). The adjunction notion is representably defined in any 2-category, in the sense that a 1-cell $u\colon A \to B$ in a 2-category $\mathcal{C}$ admits a left adjoint if and only if for all objects $X$ the functor $\hom(X, u)\colon \hom(X, A) \to \hom(X, B)$ admits a left adjoint in the usual sense. The forward implication simply follows from the fact that adjunctions are preserved by any 2-functor whereas the backward implication is a routine consequence of the *bicategorical Yoneda lemma* [24].

A closely related observation is that a functor $u\colon A \to B$ between cofibrant objects in an $\infty$-cosmos $\mathcal{K}$ has a left adjoint in the homotopy 2-category $\mathcal{K}_2$ if and only if for all cofibrant objects $X$ the functor $\mathrm{map}(X, u)\colon \mathrm{map}(X, A) \to \mathrm{map}(X, B)$ of quasi-categories has a left adjoint. Indeed, this result is an immediate consequence of the observation that the homotopy category construction provides us with a 2-functor $h\colon \underline{\mathrm{qCat}}_2 \to \underline{\mathrm{Cat}}$ which relates the 2-categorical and quasi-categorical representable functors associated with $X$ by the equation $h(\mathrm{map}(X, -)) = \hom(X, -)$. We leave the details to the intrepid reader.

3.6.14. *Observation* (right adjoint right inverse). In any 2-category $\mathcal{C}$, Lemma I.4.1.2 demonstrates that a 1-cell $f\colon B \to A$ has a right adjoint whose counit is an isomorphism if and only if there exists a 1-cell $u\colon A \to B$ for which $fu \cong \mathrm{id}_A$ and and there exists a 2-cell $\eta'\colon \mathrm{id}_B \Rightarrow uf$ for which $\eta' u$ and $f\eta'$ are both isomorphisms. The conclusion of that proof is that we can construct a unit $\eta\colon \mathrm{id}_B \Rightarrow uf$ and a counit $\epsilon\colon fu \cong \mathrm{id}_A$ from the supplied data which demonstrates that $f \dashv u$. Moreover, in the case where $fu = \mathrm{id}_A$ that counit can be taken to be an identity.

If we are given an adjunction $f \dashv u$ in $\mathcal{C}$ whose counit is an isomorphism and in which $f\colon B \to A$ happens to be an isofibration, then we may apply the argument of Example I.4.5.5 and lift the counit $\epsilon\colon fu \cong \mathrm{id}_A$ to give a 1-cell $u' \cong u$ and derive the data for an adjunction $f \dashv u'$ whose counit is an identity. This then may be regarded as being an

adjunction

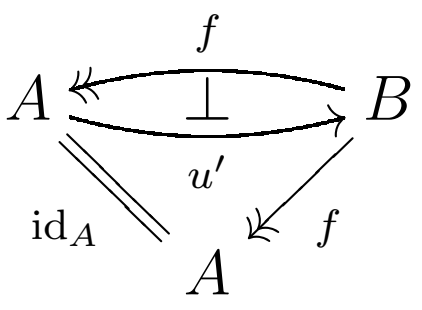

in the slice 2-category $\mathcal{C}/A$.

3.6.15. *Observation.* We conclude with a list of a few other standard 2-categorical results used in the sequel, whose proofs we leave to the reader:

- An adjunction $f \dashv u \colon A \to B$ in $\mathcal{C}$ has an counit which is an isomorphism if and only if $u$ is *representably fully faithful*, in the sense that the actual functor of hom-categories $\hom(X, u)\colon \hom(X, A) \to \hom(X, B)$ is fully faithful for all objects $X \in \mathcal{C}$.
- Suppose that we are given a pair of adjunctions $l \dashv u \dashv r$ then the counit of $l \dashv u$ is an isomorphism if and only if the unit of $u \dashv r$ is an isomorphism. Indeed, this result is a direct consequence of the last result and its dual in $\mathcal{K}^{\mathrm{co}}$.

## 4. Cartesian fibrations

Our purpose in §4.1 to build a 2-categorical theory of *cartesian fibrations* between $\infty$-categories, followed in §4.2 by a 2-category theory of *groupoidal cartesian fibrations*, which is an easier special case. Specialising to the case of quasi-categories, the structures thus defined are equivalent to those introduced in the work of Lurie [12] and Joyal [9] respectively. A companion paper applies this 2-categorical theory in certain slice categories to obtain a notion of *two-sided groupoidal cartesian fibrations* upon which the calculus of modules (profunctors) between $\infty$-categories will be founded [20].

For this section we shall assume that we are working in an abstract homotopy 2-category $\mathcal{C}$. Of course, we shall generally apply these results when $\mathcal{K}$ is an $\infty$-cosmos and $\mathcal{C}$ is its homotopy 2-category $\mathcal{K}_2$, but nothing we say here will depend on that being the case.

### 4.1. **Cartesian fibrations.**

4.1.1. **Definition** (cartesian 2-cells)**.** Suppose that $p\colon E \twoheadrightarrow B$ is an isofibration in the 2-category $\mathcal{C}$. We say that a 2-cell $\chi\colon e' \Rightarrow e\colon A \to E$ is *cartesian* for $p$ if and only if

(i) (induction) for any pair of 2-cells $\tau\colon e'' \Rightarrow e$ and $\gamma\colon pe'' \Rightarrow pe'$ with $p\tau = p\chi \cdot \gamma$ there is some $\overline{\gamma}\colon e'' \Rightarrow e'$ with $p\overline{\gamma} = \gamma$ ($\bar{\gamma}$ *lies over* $\gamma$) and the property that $\tau = \chi \cdot \bar{\gamma}$.

(ii) (conservativity) for any 2-cell $\gamma\colon e' \Rightarrow e'$ if $\chi \cdot \gamma = \chi$ and $p\gamma$ is an identity then $\gamma$ is an isomorphism.

A 2-cell $\rho\colon e \Rightarrow e'\colon A \to E$ in $\mathcal{C}$ is said to be *cocartesian* for $p$ if and only if it is cartesian for $p$ in $\mathcal{C}^{\mathrm{co}}$.

In classical 2-category theory, a "cartesian 2-cell" for $p$ would be a 2-cell for which the induction property holds strongly, in the sense that the induced 2-cell $\bar{\gamma}$ is unique. This however is a notion which we would expect to be of little use in our context, for much the same reason that we find that our homotopy 2-categories only possess weak, but not strong, limits of certain kinds.

4.1.2. *Observation* (isomorphism stability). The class of cartesian 2-cells for an isofibration $p\colon E \twoheadrightarrow B$ is closed under pre-composition and post-composition by arbitrary isomorphisms. That is if $\chi\colon e' \Rightarrow e$ is cartesian for $p$ and $\alpha\colon e \cong \bar{e}$ and $\beta\colon e' \cong \bar{e}'$ are arbitrary invertible 2-cells then $\alpha \cdot \chi\colon \bar{e}' \Rightarrow e$ and $\chi \cdot \beta\colon e' \Rightarrow \bar{e}$ are both cartesian for $p$.

Other closure properties of cartesian 2-cells under composition and left cancelation will be demonstrated later in Lemmas 5.1.8 and 5.1.9.

4.1.3. *Observation* (more conservativity). Suppose that we are given a pair of cartesian 2-cells $\chi\colon e' \Rightarrow e$ and $\chi'\colon e'' \Rightarrow e$ and a third 2-cell $\gamma\colon e'' \Rightarrow e'$ which satisfy the equation $\chi \cdot \gamma = \chi'$ and for which $p\gamma$ is an isomorphism. Then applying the induction property of the cartesian arrow $\chi'$ we may obtain a 2-cell $\gamma'\colon e' \Rightarrow e''$ in the opposite direction with $\chi' \cdot \gamma' = \chi$ and $p\gamma' = (p\gamma)^{-1}$. Now we know that $\chi \cdot (\gamma \cdot \gamma') = \chi$ and $\chi' \cdot (\gamma' \cdot \gamma) = \chi'$ and that $p(\gamma \cdot \gamma') = (p\gamma) \cdot (p\gamma)^{-1} = \mathrm{id}_{pe'}$ and $p(\gamma' \cdot \gamma) = (p\gamma)^{-1} \cdot (p\gamma) = \mathrm{id}_{pe''}$. So we may apply clause (ii) of the cartesian properties of $\chi$ and $\chi'$ to show that $\gamma \cdot \gamma'$ and $\gamma' \cdot \gamma$ are isomorphisms and thus that $\gamma$ and $\gamma'$ are both isomorphisms.

As a special case of this result, we know that if $\chi\colon e' \Rightarrow e$ and $\chi'\colon e'' \Rightarrow e$ are cartesian 2-cells with $p\chi = p\chi'$ then, using the induction property of $\chi$, there exists an induced 2-cell $\gamma\colon e'' \Rightarrow e'$ with $\chi \cdot \gamma = \chi'$ and $p\gamma = \mathrm{id}_{pe'} = \mathrm{id}_{pe''}$ and that this $\gamma$ is an isomorphism by the argument of the last paragraph.

4.1.4. **Definition** (cartesian fibration)**.** We say that an isofibration $p\colon E \twoheadrightarrow B$ in $\mathcal{C}$ is a *cartesian fibration* if and only if:

(i) For every 2-cell on the left of the following diagram

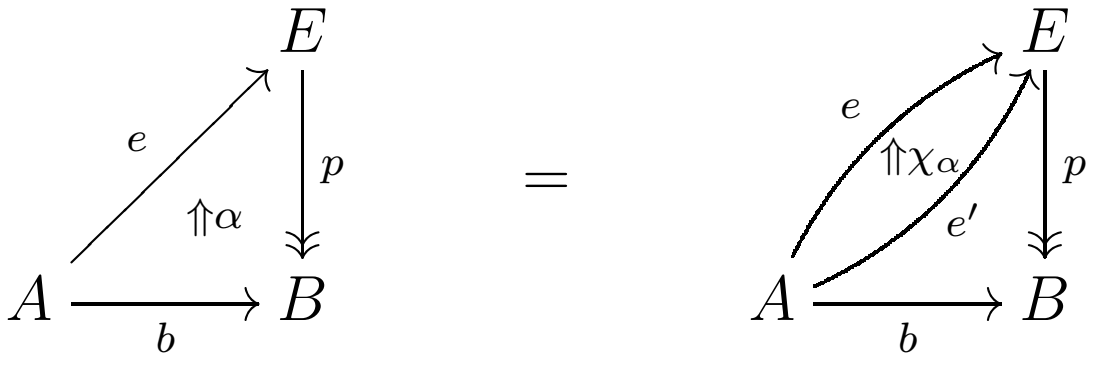

(4.1.5)

there exists a 2-cell $\chi_\alpha\colon e' \Rightarrow e$, as depicted on the right, which is cartesian for $p$ and which lies over $\alpha$, in the sense that $p\chi_\alpha = \alpha$. We call this a *cartesian lift* of $\alpha$ along $p$.

(ii) The class of cartesian 2-cells for $p$ is closed under pre-composition by all 1-cells: that is to say if $\chi\colon e' \Rightarrow e\colon A \to E$ is cartesian for $p$ and $f\colon A' \to A$ is any 1-cell then $\chi f\colon e'f \Rightarrow ef$ is again cartesian for $p$.

In line with traditional usage, we occasionally write $\alpha^*(e)$ to denote the domain of the cartesian lift $\chi_\alpha$ of the 2-cell $\alpha$ in (4.1.5).

4.1.6. *Observation* (uniqueness of cartesian lifts). The argument at the end of Observation 4.1.3 tells us that any pair of cartesian lifts $\chi\colon e' \Rightarrow e$ and $\chi'\colon e'' \Rightarrow e$ of the 2-cell in (4.1.5) are essentially unique up to an invertible 2-cell $\alpha\colon e'' \cong e'$ of their domain whose composite $p\alpha$ is an identity. It follows, in particular, that any $p$-cartesian lift of an identity is an isomorphism.

Applying this result, and Observation 4.1.2, we obtain the following completely trivial fact, which nevertheless will frequently be used. To demonstrate the stability of the class of cartesian 2-cells of $p$ under pre-composition by all 1-cells it is enough to show that for all 2-cells $\alpha$ as in (4.1.5) and all 1-cells $f\colon A' \to A$ there exists *some* cartesian lift $\chi_\alpha\colon e' \Rightarrow e$ of $\alpha$ along $p$ such that $\chi_\alpha f$ is cartesian for $p$.

4.1.7. **Proposition** (composites of cartesian fibrations). *If $q\colon B \twoheadrightarrow A$ and $p\colon E \twoheadrightarrow B$ are both cartesian fibrations, then so is their composite $qp\colon E \twoheadrightarrow A$.*

*Proof.* Note first that the class of isofibrations is closed under composition, so $qp$ is an isofibration. Now suppose that $\chi\colon e' \Rightarrow e\colon X \to E$ is cartesian for $p$ and that $p\chi\colon pe' \Rightarrow pe$ is cartesian for $q$. Then we claim that $\chi$ is also cartesian for $qp$. To prove this fact suppose first that we are given 2-cells $\tau\colon e'' \Rightarrow e$ and $\gamma\colon qpe'' \Rightarrow qpe'$ such that $qp\tau = qp\chi \cdot \gamma$. We may use the fact that $p\chi$ is cartesian for $q$ to infer that there exists some $\bar{\gamma}\colon pe'' \Rightarrow pe'$ with $q\bar{\gamma} = \gamma$ and $p\tau = p\chi \cdot \bar{\gamma}$ and then we may use the fact that $\chi$ is cartesian for $p$ to show that there exists some $\hat{\gamma}\colon e'' \Rightarrow e'$ with $p\hat{\gamma} = q\bar{\gamma} = \gamma$ and $\tau = \chi \cdot \hat{\gamma}$, which verifies the induction property of $\chi$ relative to the composite $qp$. Similarly, if $\gamma\colon e' \Rightarrow e'$ is a 2-cell with $qp\gamma$ an identity and $\chi \cdot \gamma = \chi$ then $p\chi \cdot p\gamma = p\chi$ so we may use the fact that $p\chi$ is cartesian for $q$ to show that $p\gamma$ is an isomorphism. But now we can apply the conservativity property, as discussed in Observation 4.1.3, of the cartesian 2-cell $\chi$ of $p$ to show that $\gamma$ is an isomorphism as required.

With this result in hand, it is now easily observed that we can lift a 2-cell $\alpha\colon a \Rightarrow qpe$ to a cartesian 2-cell for $qp$ by first lifting it to a cartesian 2-cell $\chi_\alpha\colon b \Rightarrow pe$ for $q$ and then lifting that to a cartesian 2-cell $\chi_{\chi_\alpha}\colon e' \Rightarrow e$ for $p$. It is also clear that the stability of the cartesian 2-cells of $qp$ under pre-composition by 1-cells follows directly from this construction of cartesian lifts and the corresponding properties of $p$ and $q$. □

In classical category theory there exists a couple of different ways of characterising cartesian fibrations in terms of certain adjunctions between comma categories. For a 2-categorical account of these equivalent descriptions the reader may wish to consult any of the papers of Street on the topic of fibrations in 2-categories and bicategories [23, 24, 25]. These equivalent characterisations also hold in our context, as we shall demonstrate in Theorem 4.1.10 below. However the 2-categorical arguments used to demonstrate the equivalence of these notions are a little more delicate, precisely because we can only rely upon the weak 2-universal properties of the comma objects used in these characterisations.

4.1.8. **Notation.** Fixing an isofibration $p\colon E \twoheadrightarrow B$ we introduce the following notation for the weakly 2-universal cones, the comma squares of Definition 3.2.1, which display the comma objects $E^2 \cong E \downarrow E$ and $B \downarrow p$:

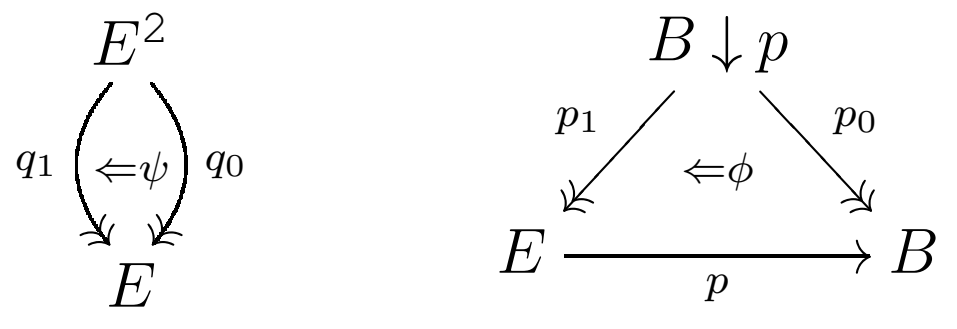

We also define three comparison 1-cells $j\colon E \to E^2$, $k\colon E^2 \to B \downarrow p$, and $i\colon E \to B \downarrow p$ using the 1-cell induction properties of these commas as follows

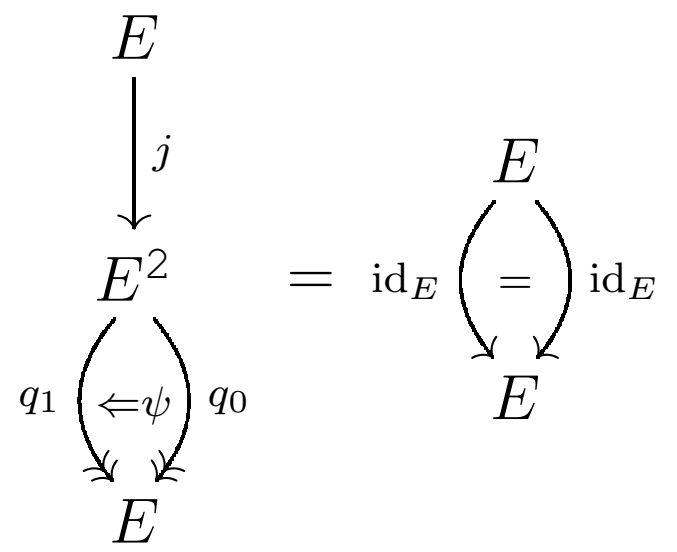

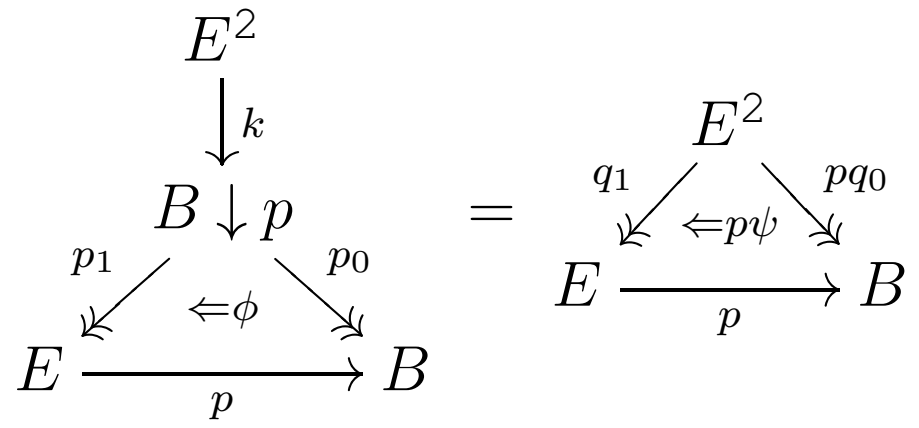

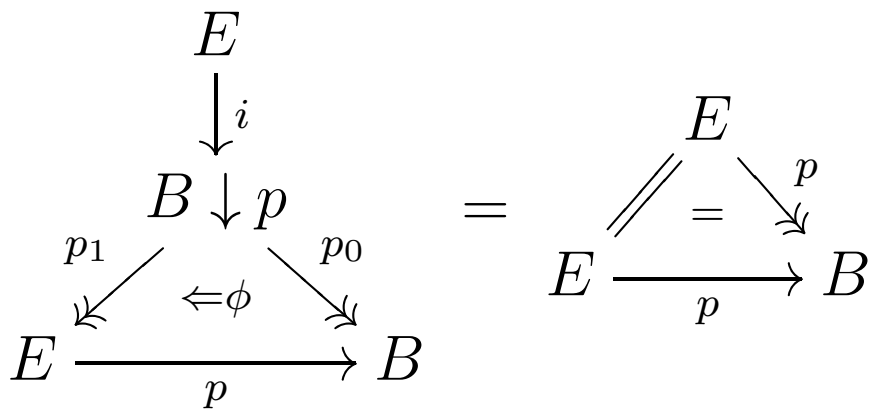

For example, $k\colon E^2 \to B{\downarrow}p$ is an induced 1-cell with the defining properties that $p_0k = pq_0$, $p_1k = q_1$, and $\phi k = p\psi$. We know, by the discussion in Observation 3.4.6, that this induced 1-cell is essentially the unique such. In other words, if $\bar{k}\colon E^2 \to B \downarrow p$ is a second 1-cell satisfying the properties $p_0\bar{k} = pq_0$, $p_1\bar{k} = q_1$ and $\phi\bar{k} = p\psi$ then there exists an invertible 2-cell $\gamma\colon k \cong \bar{k}$ which lies in the 2-category of spans from $E$ to $B$, in the sense that $q_0\gamma$ and $q_1\gamma$ are both identities.

Notice that the defining equations for $k$ and $j$ provide the following computations $p_0kj = pq_0j = p$, $p_1kj = q_1j = \mathrm{id}_E$, and $\phi kj = p\psi j = \mathrm{id}_p$. In other words, $kj$ can also be regarded as being a 1-cell induced by the 1-cell induction property of $B \downarrow p$ under the same defining equations as $i$. Hence, there exists an induced isomorphism $\mu\colon i \cong kj$ with $p_0\mu$ and $p_1\mu$ both identities.

4.1.9. **Lemma.** *The 1-cells $i$ and $j$ feature in adjunctions with the various projections from $B \downarrow p$ and $E^2$ to $E$ as follows:*

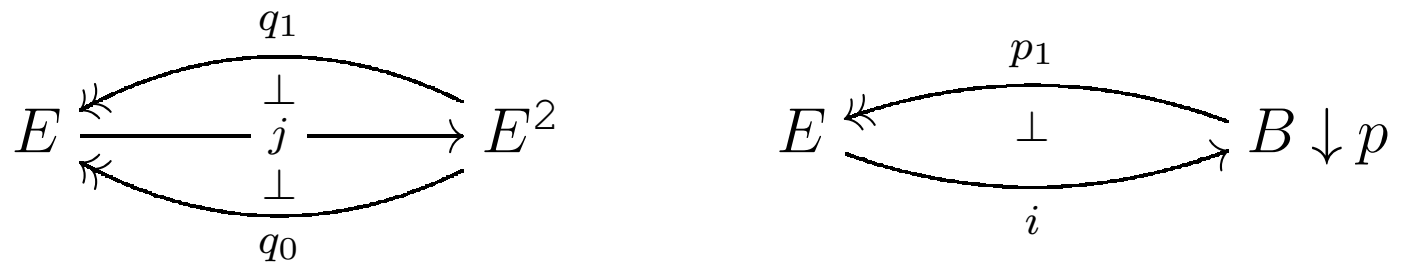

*Here the counits of the adjunctions $q_1 \dashv j$ and $p_1 \dashv i$ are both identities, as is the unit of the adjunction $j \dashv q_0$.*

*Proof.* See Lemma I.4.1.6. The general strategy for proofs of this kind is discussed in Observation 3.6.14. □

We are now ready to state our theorem characterising cartesian fibrations in terms of certain adjunctions in the 2-category $\mathcal{C}$ and its slice $\mathcal{C}/B$:

4.1.10. **Theorem.** *If $p\colon E \twoheadrightarrow B$ is an isofibration then the following are equivalent:*

*(i) The isofibration $p$ is a cartesian fibration,*

*(ii) The 1-cell $i\colon E \to B\downarrow p$ admits a right adjoint which is fibred over $B$. In other words, this condition states that there exists an adjunction*

$$\begin{array}{ccc} B\downarrow p & \underset{r}{\overset{i}{\underset{\bot}{\leftrightarrows}}} & E \\ & {\scriptstyle p_0}\searrow \quad \swarrow{\scriptstyle p} & \\ & B & \end{array} \tag{4.1.11}$$

*in the slice 2-category $\mathcal{C}/B$.*

*(iii) The 1-cell $k\colon E^2 \to B\downarrow p$ is the left adjoint part of an adjunction*

$$B\downarrow p \underset{\bar{r}}{\overset{k}{\underset{\bot}{\leftrightarrows}}} E^2 \tag{4.1.12}$$

*in the 2-category $\mathcal{C}$ whose counit is an isomorphism.*

4.1.13. *Observation.* Because we know from Lemma 4.1.9 that there is an adjunction $p_1 \dashv i$ whose counit is an identity. So it follows by Observation 3.6.15 that whenever the adjunction $i \dashv r$ exists, as in (ii) above, then its unit must be an isomorphism.

*Proof sketch, Theorem 4.1.10.* We prove these equivalences in the order (i) $\Rightarrow$ (ii) $\Rightarrow$ (iii) $\Rightarrow$ (i). A detailed proof is deferred until the appendix; for now, we content ourselves with the following sketch overview:

**(i) ⇒ (ii):** Given (i), the functor $r\colon B\downarrow p \to E$ of (ii) is defined to be the domain of a cartesian 2-cell $\chi_\phi$ lifting $\phi\colon p_0 \Rightarrow pp_1$ of 4.1.8. Pre-composing $\chi_\phi$ by $i\colon E \to B\downarrow p$ we obtain a 2-cell which we show provides an isomorphism $\chi_\phi i\colon ri \cong \mathrm{id}_E$, whose inverse we take as a candidate for the unit $\eta$ of the desired adjunction $i \dashv r$. The candidate counit $\epsilon\colon ir \Rightarrow \mathrm{id}_{B\downarrow p}$ is defined by the 2-cell induction property of $B\downarrow p$ with defining properties $p_0\epsilon = \mathrm{id}_{p_0}$ and $p_1\epsilon = \chi_\phi$. Now we establish the conditions of Observation 3.6.14 to verify that this does indeed display the desired adjunction.

**(ii) ⇒ (iii):** Given (ii), the functor $\bar{r}\colon B\downarrow p \to E^2$ of (iii) is defined by 1-cell induction for $E^2$ with defining properties $q_0\bar{r} = r$, $q_1\bar{r} = p_1$ and $\psi\bar{r} = p_1\epsilon$ (where $\epsilon$ is the counit of $i \dashv r$). Now we observe that $p_0k\bar{r} = p_0$, $p_1k\bar{r} = p_1$, and $\phi k\bar{r} = \phi$ from which it follows, by the essential uniqueness of induced 1-cells into the comma object $B\downarrow p$, that there exists an invertible 2-cell $\bar{\epsilon}\colon k\bar{r} \Rightarrow \mathrm{id}_{B\downarrow p}$. The construction of a candidate unit $\bar{\eta}\colon \mathrm{id}_{E^2} \Rightarrow \bar{r}k$ is a little more involved, yet essentially routine, and we again apply the conditions of Observation 3.6.14 to verify the desired adjunction.

**(iii) ⇒ (i):** Given (iii) and a 2-cell $\alpha\colon b \Rightarrow pe\colon A \to B$ as in (4.1.5), we observe that $\alpha$ induces a 1-cell $a_\alpha\colon A \to B\downarrow p$ and that the 2-cell $\chi_\alpha := \psi\bar{r}a_\alpha$ is a lift of $\alpha$ with codomain $e$. Finally, we use the adjunction $k \dashv \bar{r}$ to demonstrate that this choice of $\chi_\alpha$ is indeed a cartesian 2-cell for $p$ as required. □

4.1.14. *Observation.* The sketch of the proof of (iii) $\Rightarrow$(i) above tells us how to extract cartesian lifts from the data provided in an adjunction of the form given in 4.1.10(iii). Combining this with the proof of implication (ii) $\Rightarrow$ (iii), we can also extract a formula describing how to build cartesian lifts from an adjunction of the form given in 4.1.10(ii). The right adjoint $\bar{r}\colon B \downarrow p \to E^2$ is constructed as a 1-cell induced by the 2-cell $\chi := p_1\epsilon\colon r \Rightarrow p_1$, where $\epsilon$ is the counit of the adjunction $i \dashv r$. So it follows that $\chi_\alpha = \psi\bar{r}a_\alpha = p_1\epsilon a_\alpha$ expresses our cartesian lift of $\alpha$ in terms of the adjunction $i \dashv r$.

An important corollary of Theorem 4.1.10 is that our notion of cartesian fibration is representably defined.

4.1.15. **Corollary.** *Let $p\colon E \twoheadrightarrow B$ be an isofibration in $\mathcal{K}$. Then $p$ is a cartesian fibration if and only if for every cofibrant object $X \in \mathcal{K}$, the isofibration* $\mathrm{map}(X,p)\colon \mathrm{map}(X,E) \twoheadrightarrow \mathrm{map}(X,B)$ *is a cartesian fibration of quasi-categories.*

*Proof.* Recall Proposition 2.1.10, which says that for cofibrant $X$, $\mathrm{map}(X,-)\colon \mathcal{K} \to \underline{\mathrm{qCat}}$ is a functor of $\infty$-cosmoi. By Observation 3.6.13, $k\colon E^2 \to B \downarrow p$ has a right adjoint with isomorphic counit if and only if for each $X$ the functor $\mathrm{map}(X,k)\colon \mathrm{map}(X,E^2) \to \mathrm{map}(X,B\downarrow p)$ of quasi-categories has a right adjoint with isomorphic counit. However we know that $\mathrm{map}(X,-)$ preserves the simplicial limits used to define $E^2$ and $B \downarrow p$ in $h_*\mathcal{K}$ and the trivial fibrations used to define the cofibrant replacements in $\mathcal{K}_2$. It follows that there is a commutative diagram

$$\begin{array}{ccc} \mathrm{map}(X,E^2) & \xrightarrow{\mathrm{map}(X,k)} & \mathrm{map}(X,B\downarrow p) \\ \wr\big\downarrow & & \big\downarrow\wr \\ \mathrm{map}(X,E)^2 & \xrightarrow[k]{} & \mathrm{map}(X,B)\downarrow \mathrm{map}(X,p) \end{array}$$

in which the verticals are trivial fibrations and in particular equivalences of quasi-categories. It follows that $\mathrm{map}(X,k)$ admits a right adjoint with isomorphic counit if and only if $k$ admits a right adjoint with isomorphic counit. Consequently, $p$ is a cartesian fibration in $\mathcal{K}$ if and only if $\mathrm{map}(X,p)$ is a cartesian fibration of quasi-categories for all cofibrant objects $X$. Combining this with Observation 3.6.13 on the representable nature of adjunctions in $\mathcal{K}_2$ completes our proof. $\square$

4.1.16. **Example.** We will use Theorem 4.1.10 to show that the domain projection $p_0\colon E^2 \twoheadrightarrow E$ is a cartesian fibration by constructing an adjunction

$$\begin{array}{ccc} E\downarrow p_0 & \underset{r}{\overset{i}{\underset{\longrightarrow}{\overset{\longleftarrow}{\bot}}}} & E^2 \\ & {\scriptstyle p_0}\searrow \quad \swarrow{\scriptstyle p_0} & \\ & E & \end{array}$$

in $\mathcal{C}/E$. The functors $i$ and $r$ are defined by 1-cell induction

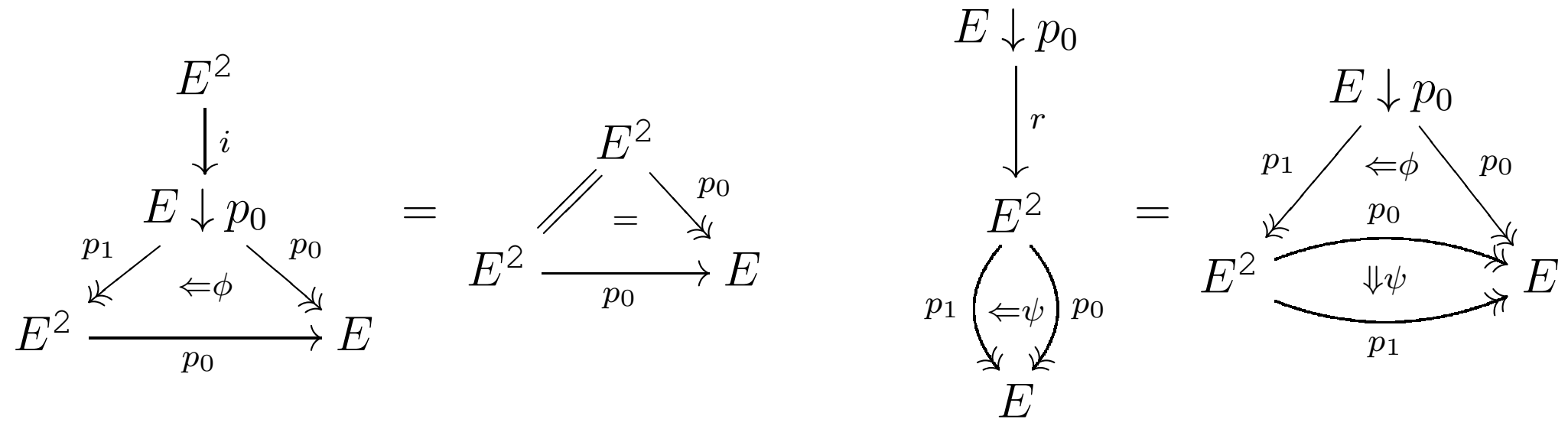

Observe that $i$ and $r$ both lie in $\mathcal{C}/E$. Observation 3.4.6 allows us to induce an isomorphism $\eta\colon 1_{E^2} \Rightarrow ri$ with $p_1\eta$ and $p_0\eta$ both identities.

By 2-cell induction for $E^2$, there is a 2-cell $\gamma\colon r \Rightarrow p_1\colon E \downarrow p_0 \to E^2$ defined by $p_0\gamma = \phi$ and $p_1\gamma = \mathrm{id}$. We use this 2-cell to define the counit $\epsilon\colon ir \Rightarrow 1_{E\downarrow p_0}$ by 2-cell induction for $E \downarrow p_0$. Its defining identities are $p_0\epsilon = \mathrm{id}$ and $p_1\epsilon = \gamma$.

The fact that $p_0\epsilon$ and $p_0\eta$ are identities ensures that this pair of 2-cells lies in $\mathcal{C}/E$. We conclude by applying Observation 3.6.14 to establish the desired adjunction $i \dashv r$. As $\eta$ is an isomorphism, it remains only to show that $\epsilon i$ and $r\epsilon$ are isomorphisms. Both follow immediately from 2-cell conservativity: $p_1\epsilon i = \gamma i$ is an isomorphism by 2-cell conservativity because $p_0\gamma i = \theta i$ is an identity. Thus $\epsilon i$ is an isomorphism. Similarly, $p_0 r\epsilon = p_0\epsilon$ is an identity and $p_1 r\epsilon = p_1 p_1 \epsilon = p_1\gamma$ is an identity, which implies that $r\epsilon$ is an isomorphism.

4.1.17. *Observation* (cartesian lifts for domain projections). We can use Observation 4.1.14 to construct the cartesian lifts for $p_0\colon E^2 \to E$. As described there, a 2-cell is represented by its induced 1-cell $\hat{\alpha}\colon A \to E\downarrow p_0$. Its cartesian lift $\chi_\alpha\colon A \to E^2$ may then be taken to be the restriction of $p_1\epsilon = \gamma\colon r \Rightarrow p_1\colon E \downarrow p_0 \to E^2$ along $\hat{\alpha}$.

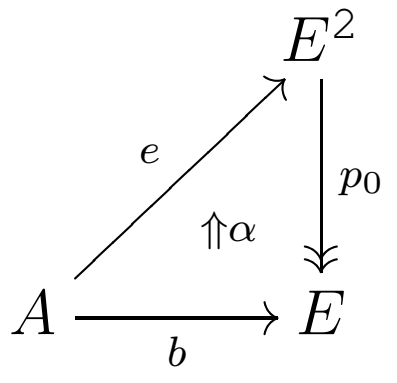

Unpacking the definition, $\chi_\alpha$ is a 2-cell induced by the weak universal property of $E^2$ with defining property $p_0\chi_\alpha = \alpha$ and $p_1\chi_\alpha = \mathrm{id}_{p_1 e}$. By Observation 4.1.6, any $p_0$-cartesian 2-cell is isomorphic to some $\chi_\alpha$, from which we conclude that the $p_0$-cartesian 2-cells $\gamma\colon e' \Rightarrow e\colon X \to E^2$ are precisely those for which $p_1\gamma$ is an isomorphism.

4.1.18. **Example.** For any functor $f\colon B \to A$, the proof given in Example 4.1.16 utilizing Theorem 4.1.10(iii), also demonstrates that the domain projection functor $p_0\colon A \downarrow f \to A$ is a cartesian fibration. For the reasons described in Observation 4.1.17, the $p_0$-cartesian 2-cells are again precisely those 2-cells whose image under $p_1\colon A\downarrow f \to B$ is an isomorphism.

We shall see in section 5 that the domain projection functors $p_0\colon f \downarrow A \to B$ and most generally, for any $g\colon C \to A$, $p_0\colon f\downarrow g \to B$ are also cartesian fibrations. This follows easily from the two cases already demonstrated by Proposition 5.2.1, which proves that cartesian fibrations are stable under the pullback construction of Definition 3.5.7.

4.1.19. **Example.** If $E$ admits pullbacks, in the sense defined in I.5.2.9, then the codomain projection $p_1 \colon E^2 \twoheadrightarrow E$ is a cartesian fibration. Applying Theorem 4.1.10.(iii) it suffices to find a right adjoint right inverse to the restriction functor $E^{2\times 2} \to E^{\lrcorner}$, where $\lrcorner$ denotes the pushout that glues two copies of $\Delta^1$ along their terminal vertex. This was done in Corollary I.5.2.20 in the special case discussed in I.5.2.22.

We have shown that if $E$ admits pullbacks then the codomain projection $p_1 \colon E^2 \twoheadrightarrow E$ is both a cartesian and a cocartesian fibration. We call such a functor a *bifibration*. Our particular interest in this notion is principally derived from the following result:

4.1.20. **Proposition.** *Let $p \colon E \twoheadrightarrow B$ be a bifibration in an abstract homotopy 2-category $\mathcal{C}$. Then any arrow $\alpha \colon a \Rightarrow b \colon X \to B$ in the base induces an adjunction $\Sigma_\alpha \dashv \alpha^* \colon E_a \to E_b$ in $\mathcal{C}/X$ between the fibres of $p$ over $a$ and $b$.*

*Proof.* The fibres over $a$ and $b$ are defined to be the pullbacks of $p \colon E \twoheadrightarrow B$ along the functors $a, b \colon X \to B$.

$$\begin{array}{ccc} E_a & \xrightarrow{i_a} & E \\ {\scriptstyle p_a}\downarrow & \lrcorner & \downarrow{\scriptstyle p} \\ X & \xrightarrow[a]{} & B \end{array} \qquad\qquad \begin{array}{ccc} E_b & \xrightarrow{i_b} & E \\ {\scriptstyle p_b}\downarrow & \lrcorner & \downarrow{\scriptstyle p} \\ X & \xrightarrow[b]{} & B \end{array}$$

By Lemma 3.5.6, the pulled-back isofibrations are formed as the parallel legs of iso-comma squares

$$\begin{array}{ccc} E_a & \xrightarrow{q_a} & E \\ {\scriptstyle p_a}\downarrow & \overset{\psi_a}{\cong} & \downarrow{\scriptstyle p} \\ X & \xrightarrow[a]{} & B \end{array} \qquad\qquad \begin{array}{ccc} E_b & \xrightarrow{q_b} & E \\ {\scriptstyle p_b}\downarrow & \overset{\psi_b}{\cong} & \downarrow{\scriptstyle p} \\ X & \xrightarrow[b]{} & B \end{array}$$

Because $p \colon E \twoheadrightarrow B$ is a cartesian fibration, by Theorem 4.1.10, $k \colon E^2 \to B \downarrow p$ admits a right adjoint with isomorphic counit. This data defines an adjunction

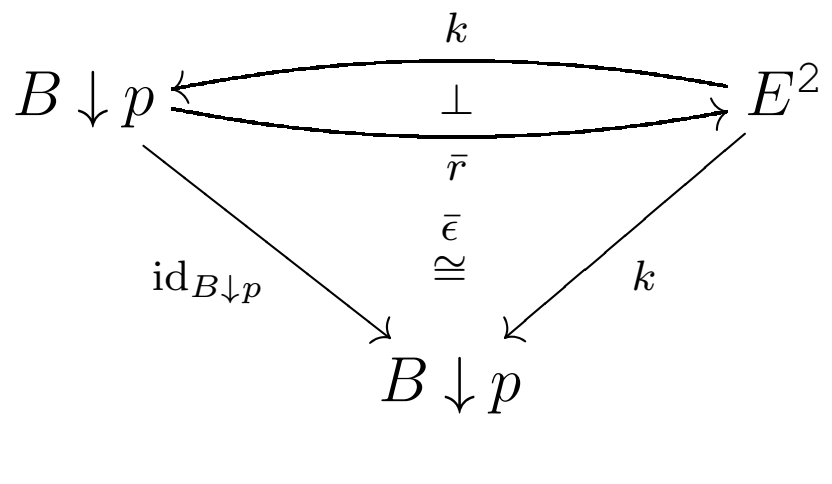

(4.1.21)

in the pseudo slice 2-category $\mathcal{C}/_{\cong}(B \downarrow p)$ introduced in Definition 3.6.8. Lemma 3.4.12 provides an iso-comma square

$$\begin{array}{ccc} B \downarrow p & \twoheadrightarrow & E \\ {\scriptstyle \ell}\downarrow & \cong & \downarrow{\scriptstyle p} \\ B^2 & \underset{p_1}{\twoheadrightarrow} & B \end{array}$$

Composition with $\ell$ defines a 2-functor $\mathcal{C}/_{\cong}(B \downarrow p) \to \mathcal{C}/_{\cong}B^2$ which carries the adjunction (4.1.21) to an adjunction in the pseudo slice 2-category over $B^2$. Note that by construction, the composite $\ell k$ is isomorphic to the map $p^2 \colon E^2 \to B^2$. Forming the iso-commas of $\ell \colon B \downarrow p \to B^2$ and $p^2 \colon E^2 \to B^2$ with $x \colon X \to B^2$, Corollary 3.6.12, via Lemma 3.6.11 allows us to transport the adjunction from the pseudo slice 2-category $\mathcal{C}/_{\cong}B^2$ to an adjunction in the usual slice 2-category $\mathcal{C}/X$, i.e., to an adjunction fibred over $X$.

To analyse this latter adjunction further, consider the following diagram:

$$\begin{array}{ccccc} E_b & \twoheadrightarrow & B \downarrow p & \twoheadrightarrow & E \\ {\scriptstyle p_b}\downarrow & \cong & {\scriptstyle l}\downarrow & \cong & \downarrow{\scriptstyle p} \\ X & \xrightarrow{x} & B^2 & \xrightarrow{p_1} & B \end{array} \qquad b = p_1 x \colon X \twoheadrightarrow B$$

By Lemma 3.4.10, the composite of the displayed iso-commas defines an iso-comma rectangle. By Corollary 3.4.9, the left-hand vertical isofibration is equivalent to $p_b \colon E_b \twoheadrightarrow X$ over $X$. This identifies one object in the adjunction fibred over $X$ as $p_b \colon E_b \twoheadrightarrow X$.

Notice also that the other object $q \colon F \twoheadrightarrow X$ on the right of this adjunction over $X$ is obtained by forming the iso-comma:

$$\begin{array}{ccc} F & \twoheadrightarrow & E^2 \\ & & {\scriptstyle k}\downarrow \;\cong \\ {\scriptstyle q}\downarrow & \cong & B \downarrow p \qquad {\scriptstyle p^2} \\ & & {\scriptstyle l}\downarrow \\ X & \underset{x}{\longrightarrow} & B^2 \end{array}$$

So we may depict the resulting adjunction in the slice $\mathcal{C}/X$ as

$$\begin{array}{ccc} E_b & \underset{\perp}{\rightleftarrows} & F \\ & {\scriptstyle p_b}\searrow \quad \swarrow{\scriptstyle q} & \\ & X & \end{array}$$

and its counit is an isomorphism

Now we also know that $p \colon E \twoheadrightarrow B$ is a co-cartesian fibration, so applying a dual argument we obtain a second adjunction fibred over $X$. Most notably its left hand object is also formed as an iso-comma between $p^2\colon E^2 \twoheadrightarrow B^2$ and $x\colon X \to B^2$, so this time we obtain an adjunction

$$\begin{array}{ccc} F & \underset{\bot}{\rightleftarrows} & E_a \\ & {}_{q}\searrow \qquad \swarrow_{p_a} & \\ & X & \end{array}$$

in $\mathcal{C}/X$ whose unit is an isomorphism. Composing these two we obtain the adjunction postulated in the statement. □

4.1.22. *Observation.* In the theory of quasi-categories there exists an important converse to the last proposition. Specifically, if we are given an adjunction $f \dashv u\colon A \to B$ of quasi-categories we may form a bifibration $p\colon E \twoheadrightarrow 2$ whose fibres are equivalent to $A$ and $B$ respectively. Indeed, this observation lies at the very heart of Lurie's presentation of the theory of adjunctions of quasi-categories [12, 5.2.2.1].

Our final result in this section specializes to the $\infty$-cosmos $\underline{\mathrm{qCat}}$ and demonstrates that the class of cartesian fibrations characterized by Definition 4.1.4 coincides precisely with the class of cartesian fibrations between quasi-categories defined in [12, 2.4.2.1]. Our proof uses the characterization provided by Theorem 4.1.10.(iii). In the homotopy 2-category of quasi-categories, the induced functor $k\colon E^2 \twoheadrightarrow B \downarrow p$ can be modeled by an isofibration, namely the Leibniz cotensor of $p\colon E \twoheadrightarrow B$ with $d^0\colon \Delta^0 \to \Delta^1$. Lifting the isomorphic counit along the isofibration $k\colon E^2 \twoheadrightarrow B \downarrow p$ as in Observation 3.6.14, we might as well assume that $k$ admits a *right adjoint right inverse* (with counit an identity). Such isofibrations of quasi-categories can be characterized by a lifting property, which we now recall.

I.4.4.12. **Lemma** (right adjoint right inverse as a lifting property)**.** *A isofibration $f\colon A \twoheadrightarrow B$ of quasi-categories admits a right adjoint right inverse in $\underline{\mathrm{qCat}}_2$ if and only if for all $b \in B_0$ there exists $ub \in A_0$ with $fub = b$ and so that any lifting problem with $n \geq 1$*

$$\begin{array}{ccccc} \Delta^0 & \xrightarrow{\{n\}} & \partial\Delta^n & \longrightarrow & A \\ & & \downarrow & \nearrow & \downarrow f \\ & & \Delta^n & \longrightarrow & B \end{array} \qquad (\text{top composite } ub) \tag{4.1.24}$$

*has a solution.*

4.1.24. **Corollary.** *An isofibration $p\colon E \twoheadrightarrow B$ of quasi-categories is a cartesian fibration if and only if any $\alpha\colon b \to pe \in B_1$ admits a lift $\chi\colon e' \to e \in E_1$ so that any lifting problem for $n \geq 2$*

$$\begin{array}{ccccc} \Delta^1 & \xrightarrow{\{n-1,n\}} & \Lambda^{n,n} & \longrightarrow & E \\ & & \downarrow & \nearrow & \downarrow p \\ & & \Delta^n & \longrightarrow & B \end{array} \qquad (\text{top composite } \chi) \tag{4.1.25}$$

*has a solution.*

Definition [12, 2.4.2.1] only requires that $p$ is an inner fibration with the lifting property (4.1.25), but it follows easily that any such $p$ must be an isofibration.

*Proof.* By Theorem 4.1.10, $p$ is a cartesian fibration if and only if $k\colon E^2 \twoheadrightarrow B\downarrow p$ admits a right adjoint right inverse. On recalling that $k$ is the Leibniz hom of $\delta^0\colon \Delta^0 \to \Delta^1$ and $p$ and transposing the lifting property of Lemma I.4.4.12, we see that this is the case if and only if any $\alpha\colon b \to pe \in (B\downarrow p)_0$ admits a lift $\chi\colon e' \to e \in (E^2)_0$ along $k$ so that any lifting problem

$$\begin{array}{ccccc} \Delta^1 & \xrightarrow{\{n\}\times\Delta^1} & \partial\Delta^n \times \Delta^1 \cup_{\partial\Delta^n\times\Delta^0} \Delta^n \times \Delta^0 & \longrightarrow & E \\ & & \downarrow & \nearrow & \downarrow p \\ & & \Delta^n \times \Delta^1 & \longrightarrow & B \end{array} \qquad (\Delta^1 \xrightarrow{\chi} E)$$

has a solution.

To solve this lifting problem, we must find lifts for each of the $n+1$ shuffles of $\Delta^n \times \Delta^1$. We number these shuffles $0, \ldots, n$ starting from the closed end of the cylinder. Proceeding inductively for $k < n$, we choose a lift for the $k$th shuffle by filling a $\Lambda^{n+1,k+1}$ horn. To lift the $n$th shuffle, we're required to fill a $\Lambda^{n+1,n+1}$ horn whose $\{n, n+1\}$ edge is $\chi$. If such horns can be filled, then we can complete our construction of the lift $\Delta^n \times \Delta^1 \to E$ proving that $k\colon E^2 \twoheadrightarrow B\downarrow p$ admits a right adjoint right inverse.

Conversely, the right adjoint right inverse adjunction $k \dashv \bar{r}$ constructed in Theorem 4.1.10 is fibered over $E$, so we can pull it back along a vertex $e\colon \Delta^0 \to E$ to obtain a right adjoint right inverse to $E\downarrow e \twoheadrightarrow B\downarrow pe$. A special case of Proposition I.2.4.13 proves that this isofibration is equivalent to $E_{/e} \twoheadrightarrow B_{/pe}$, so we again have a right adjoint right inverse. Applying Lemma I.4.4.12, this says that any $\alpha\colon b \to pe \in B_1$ admits a lift $\chi\colon e' \to e \in E_1$ so that any lifting problem

$$\begin{array}{ccccc} \Delta^0 & \xrightarrow{\{n\}} & \partial\Delta^n & \longrightarrow & E_{/e} \\ & & \downarrow & \nearrow & \downarrow f \\ & & \Delta^n & \longrightarrow & B_{/pe} \end{array} \qquad (\Delta^0 \xrightarrow{\chi} E_{/e})$$

has a solution. The adjoint form of this lifting property is exactly (4.1.25). □

## 4.2. **Groupoidal cartesian fibrations.**

We now turn our attention to groupoidal cartesian fibrations, a special case of the cartesian fibrations defined in 4.1.4. Let $\mathcal{C}$ be an abstract homotopy 2-category.

4.2.1. **Definition** (groupoidal objects)**.** We say that an object $A$ in the 2-category $\mathcal{C}$ is *groupoidal* if and only if every 2-cell $\gamma\colon a \Rightarrow a'\colon X \to A$ with codomain $A$ is an isomorphism.

4.2.2. **Lemma.** *An isofibration $p\colon E \twoheadrightarrow B$ in $\mathcal{C}$ is a groupoidal object of $\mathcal{C}/B$ if and only if it is conservative.*

A 1-cell in a 2-category is *conservative* if and only if it is representably conservative, i.e., reflects invertible 2-cells.

*Proof.* Unwinding Definition 4.2.1, an object $f\colon A \to B$ is groupoidal in $\mathcal{C}/B$ if and only if, as a 1-cell of $\mathcal{C}$, it has the property that whenever $\gamma$ is a 2-cell with $f\gamma$ an *identity* then $\gamma$ is an isomorphism. It is clear from this that all conservative 1-cells $f\colon A \to B$ of $\mathcal{C}$ are, in particular, groupoidal when considered as 0-cells of $\mathcal{C}/B$.

Conversely, suppose that $p\colon E \twoheadrightarrow B$ is both an isofibration and a groupoidal object of $\mathcal{C}/B$. Suppose now that $\gamma\colon e \Rightarrow e'$ is a 2-cell with $p\gamma$ an isomorphism. Then we can lift $p\gamma$ to an isomorphism $\bar{\gamma}\colon e'' \Rightarrow e'$ with $p\bar{\gamma} = p\gamma$ and observe that the 2-cell $\bar{\gamma}^{-1} \cdot \gamma\colon e \Rightarrow e''$ has $p(\bar{\gamma}^{-1} \cdot \gamma) = (p\gamma)^{-1} \cdot p\gamma = \mathrm{id}_{pe}$. It follows, since $p\colon E \twoheadrightarrow B$ is groupoidal in $\mathcal{C}/B$, that we may infer that $\bar{\gamma}^{-1} \cdot \gamma$ is an isomorphism and thus that $\gamma$ itself is an isomorphism as required. □

4.2.3. **Definition** (groupoidal cartesian fibrations)**.** An isofibration $p\colon E \twoheadrightarrow B$ of $\mathcal{C}$ is a *groupoidal cartesian fibration* if and only it is a cartesian fibration and it is groupoidal as an object of the slice $\mathcal{C}/B$.

As Lemma 4.2.2 reveals, groupoidalness of an isofibration is also a representably defined notion. Combining this with Corollary 4.1.15, the following result is immediate.

4.2.4. **Corollary.** *Let $p\colon E \twoheadrightarrow B$ be an isofibration in $\mathcal{K}$. Then $p$ is a groupoidal cartesian fibration if and only if for every cofibrant object $X \in \mathcal{K}$, the induced isofibration* $\mathrm{map}(X, p)\colon \mathrm{map}(X, E) \twoheadrightarrow \mathrm{map}(X, B)$ *is a groupoidal cartesian fibration of quasi-categories.* □

4.2.5. **Proposition.** *An isofibration $p\colon E \twoheadrightarrow B$ is a groupoidal cartesian fibration if and only if every 2-cell $\alpha\colon b \Rightarrow pe\colon X \to B$ has an essentially unique lift $\chi\colon e' \Rightarrow e\colon X \to E$. Hence, if $p\colon E \twoheadrightarrow B$ is a groupoidal cartesian fibration, then every 2-cell with codomain $E$ is $p$-cartesian.*

Here the term "essentially unique" means that if $\chi\colon e' \Rightarrow e$ and $\tau\colon e'' \Rightarrow e$ are two lifts of $\alpha$ with the same codomain then there exists some isomorphism $\gamma\colon e'' \Rightarrow e'$ with $\chi \cdot \gamma = \tau$ and for which $p\gamma$ is an identity 2-cell.

*Proof.* To prove necessity, since $p$ is a cartesian fibration we know that any 2-cell $\alpha\colon b \Rightarrow pe$ has a lift $\chi\colon e' \Rightarrow e$ which is cartesian for $p$. To prove the essential uniqueness of this lift, suppose that $\tau\colon e'' \Rightarrow e$ is any other 2-cell with $p\tau = \alpha$ and observe that it factors through $\chi$, by the 1-cell induction property of that cartesian 2-cell, to give a 2-cell $\gamma\colon e'' \Rightarrow e'$ which lies over the identity 2-cell on $b$. It follows then that $\gamma$ is also an isomorphism, as required, since $p\colon E \twoheadrightarrow B$ is groupoidal.

Conversely, to prove sufficiency start by proving that the stated condition ensures that $p\colon E \twoheadrightarrow B$ is groupoidal. We do this by considering a 2-cell $\chi\colon e' \Rightarrow e$ for which $p\chi$ is an identity, and observing that the identity $\mathrm{id}_e\colon e \Rightarrow e$ is also a lift for $p\chi$. Applying the postulated essential uniqueness of such lifts we can demonstrate that these two lifts are isomorphic and thus that $\chi$ is an isomorphism since $\mathrm{id}_e$ patently is.

To complete our proof it clearly suffices now to show that every 2-cell $\chi\colon e' \Rightarrow e$ with domain $E$ is cartesian for $p$, simply because then the lifts assumed in the statement provide the cartesian lifts we seek to show that $p$ is a cartesian fibration. The required conservativity property of $\chi$ follows immediately from the argument of the last paragraph. It remains to show that if we are given $\tau\colon e'' \Rightarrow e$ and $\gamma\colon pe'' \Rightarrow pe'$ with $p\chi \cdot \gamma = p\tau$ then we can construct a 2-cell $\bar{\gamma}$ with $p\bar{\gamma} = \gamma$ and $\chi \cdot \bar{\gamma} = \tau$. To do this we take the lift $\mu\colon e^* \Rightarrow e'$ of $\gamma\colon pe'' \Rightarrow pe'$ guaranteed by the statement and observe that $p(\chi \cdot \mu) = p\chi \cdot p\mu = p\chi \cdot \gamma = p\tau$. So we now have two lifts $\chi \cdot \mu$ and $\tau$ of the same 2-cell $p\tau$, and so the essential uniqueness of such things provides an invertible 2-cell $\alpha\colon e' \Rightarrow e^*$ with $p\alpha$ an identity and $\chi \cdot \mu \cdot \alpha = \tau$. Finally $p(\mu \cdot \alpha) = p\mu \cdot p\alpha = p\mu \cdot \mathrm{id}_{pe''} = p\mu = \gamma$, so it follows that the 2-cell $\bar{\gamma} := \mu \cdot \alpha$ provides the factorisation of $\tau$ through $\chi$ that we seek. □

4.2.6. **Lemma.** *If $p$ and $pq$ are groupoidal cartesian fibrations, then so is $q$.*

*Proof.* This result follows via routine application of the characterisation of groupoidal cartesian fibrations given in Proposition 4.2.5, the details of which are left to the reader. □

4.2.7. **Proposition.** *An isofibration $p\colon E \twoheadrightarrow B$ is a groupoidal cartesian fibration if and only if the functor $k\colon E^2 \to B \downarrow p$ is an equivalence.*

*Proof.* First assume that $p\colon E \twoheadrightarrow B$ is a groupoidal cartesian fibration. Since $p$ is, in particular, a cartesian fibration, Theorem 4.1.10 tells us that $k\colon E^2 \to B \downarrow p$ has a right adjoint $\bar{r}\colon B \downarrow p \to E^2$ with unit $\bar{\eta}\colon \mathrm{id}_{E^2} \Rightarrow \bar{r}k$ and isomorphic counit $\bar{\epsilon}\colon k\bar{r} \cong \mathrm{id}_{B\downarrow p}$. To show that $k \dashv \bar{r}$ is an adjoint equivalence, it suffices to prove that $\bar{\eta}$ is an isomorphism.

Now observe that $k\bar{\eta}$ is also an isomorphism, since the the counit of this adjunction is an isomorphism, and that the defining properties of $k$ give us $p_0k\bar{\eta} = pq_0\bar{\eta}$ and $p_1k\bar{\eta} = q_1\bar{\eta}$, so in particular $q_1\bar{\eta}$ and $pq_0\bar{\eta}$ are both isomorphisms. However $p\colon E \twoheadrightarrow B$ is assumed groupoidal in $\mathcal{C}/B$ so we know, by Lemma 4.2.2, that it is conservative as a 1-cell in $\mathcal{C}$ and thus that $q_0\bar{\eta}$ is also an isomorphism. Now we can apply the conservativity property of the comma object $E^2$ to show that $\bar{\eta}$ is an isomorphism and it follows that $k \dashv \bar{r}$ is an adjoint equivalence.

Conversely, suppose that $k\colon E^2 \to B \downarrow p$ is an equivalence. Then we know that we can pick an equivalence inverse $\bar{r}\colon B \downarrow p \to E^2$ along with a unit and counit 2-cells in a way which gives us an adjoint equivalence $k \dashv \bar{r}$. The counit of that adjoint equivalence is an isomorphism, so Theorem 4.1.10 tells us that $p\colon E \twoheadrightarrow B$ is a cartesian fibration. It remains to prove that it is also conservative. To that end, suppose that $\gamma\colon e \Rightarrow e'\colon X \to E$ is a 2-cell with $p\gamma$ an isomorphism and apply the 1-cell induction property of the comma object $E^2$ to obtain a 1-cell $g\colon X \to E^2$ with $q_0g = e$, $q_1g = e'$, and $\psi g = \gamma$. Using 2-cell induction property of $E^2$ we can also construct a 2-cell $\bar{\gamma}\colon g \Rightarrow je'$, where $j\colon E \to E^2$ is the 1-cell introduced in 4.1.8, with defining properties $q_0\bar{\gamma} = \gamma$ and $q_1\bar{\gamma} = \mathrm{id}_{e'}$. Now observe, from the construction of $\bar{\gamma}$ and the definition of $k$, that $p_0k\bar{\gamma} = pq_0\bar{\gamma} = p\gamma$, which is an isomorphism by assumption, and $p_1k\bar{\gamma} = q_1\gamma = \mathrm{id}_{e'}$, so we may apply the conservativity property of the comma object $B \downarrow p$ to show that $k\bar{\gamma}$ is an isomorphism. We also assumed that $k$ is an equivalence so it is, in particular, conservative and it follows that $\bar{\gamma}$ is therefore an isomorphism and thus that $\gamma = q_0\bar{\gamma}$ is also such. This completes our proof that $p$ is conservative and thus that it is a groupoidal object in $\mathcal{C}/B$. □

4.2.8. **Corollary.** *Any representable surjective equivalence $p\colon E \xrightarrow{\sim\!\!\!\!\rightarrow} B$ is a groupoidal cartesian fibration.*

Recall that trivial fibrations in an $\infty$-cosmos define representable surjective equivalences by Lemma 3.1.4.

*Proof.* By Recollection 3.1.3, a surjective equivalence $p\colon E \xrightarrow{\sim\!\!\!\!\rightarrow} B$ admits a right inverse $r\colon B \to E$ to $p$ and an isomorphism $\gamma\colon rp \cong \mathrm{id}_E$ with $\gamma r = \mathrm{id}_r$ and $p\gamma = \mathrm{id}_p$. We can use the 1-cell induction property of the comma object $E^2$ to induce a 1-cell $\bar{r}\colon B\downarrow p \to E^2$ with defining properties $q_0\bar{r} = rp_0$, $q_1\bar{r} = p_1$, and $\psi\bar{r} = \gamma p_1 \cdot r\phi$ which we claim is an equivalence inverse to $k\colon E^2 \to B\downarrow p$. To validate this claim, start by observing that the equalities $p_0k\bar{r} = pq_0\bar{r} = prp_0 = p_0$, $p_1k\bar{r} = q_1\bar{r} = p_1$, and $\phi k\bar{r} = p\psi\bar{r} = p(\gamma p_1 \cdot r\phi) = p\gamma p_1 \cdot pr\phi = \mathrm{id}_p p_1 \cdot \phi = \phi$ follow directly from the definitions of $\bar{r}$ and $k$ and the fact that $r$ is right inverse to $p$. Consequently $\mathrm{id}_{B\downarrow p}$ and $k\bar{r}$ both enjoy the same defining properties as induced 1-cells into $B\downarrow p$, so we can apply Observation 3.4.6 to show that $k\bar{r} \cong \mathrm{id}_{B\downarrow p}$. To construct an isomorphism $\mathrm{id}_{E^2} \cong \bar{r}k$ first observe that $q_0\bar{r}k = rp_0k = rpq_0$, $q_1\bar{r}k = p_1k = q_1$, and $\psi\bar{r}k = (\gamma p_1 \cdot r\phi)k = \gamma p_1 k \cdot rp\psi = \gamma q_1 \cdot rp\psi = \psi \cdot \gamma q_0$, where the last step of this computation is a middle four interchange. From these it is easily checked that the pair $\gamma q_0\colon rpq_0 \cong q_0$ and $\mathrm{id}_{q_1}$ satisfy the condition required to induce a 2-cell $\mathrm{id}_{E^2} \Rightarrow \bar{r}k$, by the weak 2-universal property of $E^2$, which is an isomorphism by the conservativity property as required. □

4.2.9. **Example.** In the $\infty$-cosmos $\underline{\mathrm{qCat}}$, Proposition 4.2.7 can be used to show that the class of groupoidal cartesian fibrations coincides exactly with the class of *right fibrations* introduced by Joyal [8]. Observe that in this case $k\colon E^2 \twoheadrightarrow B\downarrow p$ can be modeled by the Leibniz cotensor of $p\colon E \twoheadrightarrow B$ with $d^0\colon \Delta^0 \to \Delta^1$.

Proposition 4.2.7 characterizes the groupoidal right fibrations $p$ as those maps for which $k\colon E^2 \xrightarrow{\sim\!\!\!\!\rightarrow} B\downarrow p$ is a trivial fibration, i.e., for which the Leibniz cotensor of $p$ and $d^0$ has the right lifting property against $\partial\Delta^n \hookrightarrow \Delta^n$ for any $n$. By adjunction, this condition asserts that $p$ lifts against the Leibniz product $(\partial\Delta^n \hookrightarrow \Delta^n) \mathbin{\widehat{\times}} (d^0\colon \Delta^0 \to \Delta^1)$. By a standard combinatorial lemma due to Joyal [12, 2.1.2.6], these Leibniz products generate the same weakly saturated class as the right horn inclusions, from which we conclude that this right lifting property holds precisely when $p$ is a right fibration.

4.2.10. *Observation* (biterminal objects in the homotopy 2-category). By 2.1.1(a), any $\infty$-cosmos $\mathcal{K}$ has a terminal object 1, with the universal property that $\mathrm{map}(X,1) \cong \Delta^0$ and hence that $\hom(X,1) \cong 1$ for any cofibrant object $X$. If 1 is cofibrant, then it defines a 2-terminal object in $\mathcal{K}_2$. If not, then its cofibrant replacement, which we also denote by 1, is a *biterminal object*, meaning that the hom-category $\hom(X,1)$ is equivalent to the terminal category for all $X \in \mathcal{K}_2$. Explicitly, this universal property says that any object has at least one functor $!\colon X \to 1$ and any parallel pair of functors with codomain 1 are connected by a unique isomorphism. A morphism $x\colon 1 \to X$ will be called a *point* of $X$.

4.2.11. **Example.** The domain projection functor $p_0 \colon B \downarrow b \to B$ is groupoidal if $b \colon 1 \to B$ is a point. Example 4.1.18, relying on Proposition 5.2.1, shows that $p_0$ is cartesian. Groupoidalness is immediate from the conservativity property of 2-cell induction, which implies that $p_0$ is conservative, by the fact that the codomain of $p_1 \colon B \downarrow b \to 1$ is biterminal.

Our next objective is to complete the proofs of Examples 4.1.18 and 4.2.11 by demonstrating that the notions of cartesian and groupoidal cartesian fibrations are pullback stable. We turn to this topic, among others, in the next section.

## 5. Cartesian functors and pullbacks of cartesian fibrations

Our aim in this section is to conclude the unfinished business of §4 and show that:

(i) Cartesian 2-cells compose and can be canceled on the left.
(ii) Pullbacks of (groupoidal) cartesian functors are (groupoidal) cartesian functors.

To begin this work, in §5.1 we introduce cartesian functors and prove a relative version of Theorem 4.1.10. This allows us to finish (i). In §5.2, we prove pullback stability and show moreover that pullback squares are cartesian functors in the sense to be introduced.

### 5.1. Cartesian functors.

5.1.1. **Definition.** Let $p \colon E \twoheadrightarrow B$ and $q \colon F \twoheadrightarrow A$ be a pair of cartesian fibrations in $\mathcal{C}$. A commutative square

$$\begin{array}{ccc} E & \xrightarrow{g} & F \\ {\scriptstyle p}\big\downarrow & & \big\downarrow{\scriptstyle q} \\ B & \xrightarrow[f]{} & A \end{array} \tag{5.1.2}$$

defines a *cartesian functor* if and only if $g$ preserves cartesian 2-cells: i.e., if whiskering with $g$ carries $p$-cartesian 2-cells to $q$-cartesian 2-cells.

5.1.3. *Observation* (all functors into a groupoidal fibration are cartesian). If $q$ is a groupoidal cartesian fibration, then Proposition 4.2.5 noted that any 2-cell with codomain $F$ is $q$-cartesian, and, trivially, any commutative square from with codomain $q$ defines a cartesian functor.

We now extend the equivalent characterizations of cartesian fibrations provided by Theorem 4.1.10 to cartesian functors. Note that any functor $g \colon E \to F$ induces a functor $g^2 \colon E^2 \to F^2$, well defined up to an isomorphic 2-cell over $g$, and similarly any commutative square (5.1.2) induces a functor

$$\langle g, f \rangle \colon B \downarrow p \to A \downarrow q$$

over $f$ and $g$.

5.1.4. **Theorem.** *Let $p \colon E \twoheadrightarrow B$ and $q \colon F \twoheadrightarrow A$ be cartesian fibrations in $\mathcal{C}$ that commute with the pair of functors $g \colon E \to F$ and $f \colon B \to A$. The following are equivalent:*

*(i) The pair $(g, f)$ defines a cartesian functor.*

*(ii) The mate of the canonical isomorphism*

$$\begin{array}{ccc} E & \xrightarrow{g} & F \\ {\scriptstyle i}\downarrow & \cong & \downarrow{\scriptstyle i} \\ B \downarrow p & \xrightarrow[\langle g,f\rangle]{} & A \downarrow q \end{array} \tag{5.1.5}$$

*in the diagram of functors over $f\colon B \to A$ is an isomorphism.*

*(iii) The mate of the canonical isomorphism*

$$\begin{array}{ccc} E^2 & \xrightarrow{g^2} & F^2 \\ {\scriptstyle k}\downarrow & \cong & \downarrow{\scriptstyle k} \\ B \downarrow p & \xrightarrow[\langle g,f\rangle]{} & A \downarrow q \end{array} \tag{5.1.6}$$

*in the diagram of induced functors is an isomorphism.*

*Proof sketch.* We will demonstrate that (i) ⇔ (ii) and (i) ⇔ (iii). The idea of each proof is similar. Conditions (ii) and (iii) imply that $g$ preserves the explicitly chosen cartesian lifts described in Observation 4.1.14 up to isomorphism, which by Observation 4.1.2 implies that $f$ preserves all cartesian 2-cells. Conversely, assuming (i), we must show that a whiskered copy of the counit of $i \dashv r$ and of the unit of $k \dashv \bar{r}$ are isomorphisms. The counit of $i \dashv r$ and the unit of $k \dashv \bar{r}$ each encode the data of the factorization of a 2-cell through the cartesian lift of its projection. It follows from (i) that the 2-cells in question are themselves cartesian and the factorizations live over identities. Thus Observation 4.1.3 implies that the desired 2-cells are isomorphisms. □

*Proof.* **(i) ⇒ (ii):** By Observation 4.1.13 the unit of $i \dashv r$ is an isomorphism. Thus, the mate of (5.1.5) is an isomorphism if and only if the image of the counit $\epsilon$ of $i \dashv r$ under $r\langle g, f\rangle\colon B \downarrow p \to A \downarrow q \to F$ is an isomorphism. Recall that $\epsilon\colon ri \Rightarrow \mathrm{id}\colon B \downarrow p \to B \downarrow p$ has defining equations $p_0\epsilon = \mathrm{id}_{p_0}$ and $\pi_1\epsilon = \chi_\phi\colon r \Rightarrow p_1$, where $\phi$ is the canonical 2-cell under the comma object $B \downarrow p$. Whiskering with $B \downarrow p \to A \downarrow q$, we get a 2-cell whose projection along $p_0$ is still an identity and whose projection along $p_1$ is $g\chi_\phi\colon gr \Rightarrow gp_1$.

Applying $r\colon A \downarrow q \to F$, we get a 2-cell $r\langle g, f\rangle\epsilon$ whose domain is the domain of a $q$-cartesian lift of the identity at $qgr$, some isomorphism $q^*(gr) \cong gr$ projecting to identity, and whose codomain is the domain of a $q$-cartesian lift of $f\phi\colon fp_0 \Rightarrow fpp_1 = qgp_1$, the 2-cell inducing $B \downarrow p \to A \downarrow q$. The 2-cell from the former to the latter is defined by factoring $q^*(gr) \cong gr \overset{g\chi_\phi}{\Longrightarrow} gp_1$ through the $q$-cartesian lift of $f\phi$. By (i) and Observation 4.1.2, $q^*(gr) \cong gr \overset{g\chi_\phi}{\Longrightarrow} gp_1$ is also a $q$-cartesian lift of $f\phi$, so Observation 4.1.3 implies that this induced 2-cell, the image of the counit, is an isomorphism, as required.

**(i) ⇒ (iii):** Because the counits of the adjunctions $k \dashv \bar{r}$ are isomorphisms, the mate of (5.1.6) is an isomorphism if and only if the restriction $\bar{\eta}g^2\bar{r}$ of the unit of $k \dashv \bar{r}$ is an isomorphism. We prove this using 2-cell conservativity for $F^2$. The codomain projection of $\bar{\eta}$ is an identity, so it suffices to consider the domain projection. In (7.0.5) this is defined to be the composite of isomorphisms, which we can ignore, and $rk\tau$, the 2-cell $\tau\colon jq_0 \Rightarrow \mathrm{id}$ being the counit of the adjunction of Lemma 4.1.9. The idea is that $\tau$ encodes the canonical commutative square from the identity on the domain of a generic arrow $f_0 \Rightarrow f_1$ to that arrow. Then $rk\tau$ is the factorization of $f_0 \Rightarrow f_1$ through the $q$-cartesian lift of its projection. This factorization necessarily lies over the identity on $qf_0$ and so is an isomorphism if and only if $f_0 \Rightarrow f_1$ is $q$-cartesian.

When we restrict $\bar{\eta}$ along $g^2\bar{r}$, we restrict the generic arrow to one of the form $g\chi\colon ge' \Rightarrow ge$, i.e., to $g\colon E \to F$ applied to a $p$-cartesian lift of some $b \Rightarrow pe$. By (i), this arrow $q$-cartesian, which demonstrates that $rk\tau g^2\bar{k}$ is an isomorphism, as desired.

**(iii) ⇒ (i):** In the presence of a right adjoint $\bar{r}$ to $k$, the image under $g$ of the cartesian lift of a 2-cell $\alpha$ represented by $A \xrightarrow{a_\alpha} B \downarrow p$ is the composite $g\psi\bar{r}a_\alpha$, the composite 2-cell displayed on the top-right below.

$$\begin{array}{ccccccc} A & \xrightarrow{a_\alpha} & B\downarrow p & \xrightarrow{\bar{r}} & E^2 & \overset{\Uparrow\psi}{\rightrightarrows} & E \\ & & {\scriptstyle\langle g,f\rangle}\big\downarrow & \cong & \big\downarrow{\scriptstyle g^2} & & \big\downarrow{\scriptstyle g} \\ & & A\downarrow q & \xrightarrow[\bar{r}]{} & F^2 & \overset{\Uparrow\psi}{\rightrightarrows} & F \end{array}$$

By the definition of $g^2$, this 2-cell equals the composite $\psi g^2\bar{r}a_\alpha$. Assuming (iii), this is isomorphic to the composite $\psi\bar{r}\langle g,f\rangle a_\alpha$ along the bottom-left, which represents the $q$-cartesian lift of the 2-cell $f\alpha$. Thus, $g$ preserves the cartesian lifts specified by the adjunctions $k \dashv \bar{r}$, which by Observation 4.1.6 suffices to prove that $g$ preserves all cartesian 2-cells.

**(ii) ⇒ (i):** Similarly, in the presence of a right adjoint $r$ to $i$, the image under $g$ of the cartesian lift of a 2-cell $\alpha$ represented by $A \xrightarrow{a_\alpha} B\downarrow p$ is the 2-cell displayed in the following commutative diagram.

$$\begin{array}{ccccccc} & & & E & & & \\ & & {\scriptstyle r}\nearrow & \Downarrow\epsilon & \searrow{\scriptstyle i} & & \\ A & \xrightarrow{a_\alpha} & B\downarrow p & \xrightarrow[1]{\qquad\qquad} & B\downarrow p & \xrightarrow{p_1} & E \\ & & {\scriptstyle\langle g,f\rangle}\big\downarrow & & \big\downarrow{\scriptstyle\langle g,f\rangle} & & \big\downarrow{\scriptstyle g} \\ & & A\downarrow q & \xrightarrow[1]{\qquad\qquad} & A\downarrow q & \xrightarrow[p_1]{} & F \end{array}$$

This 2-cell diagram can be rewritten as

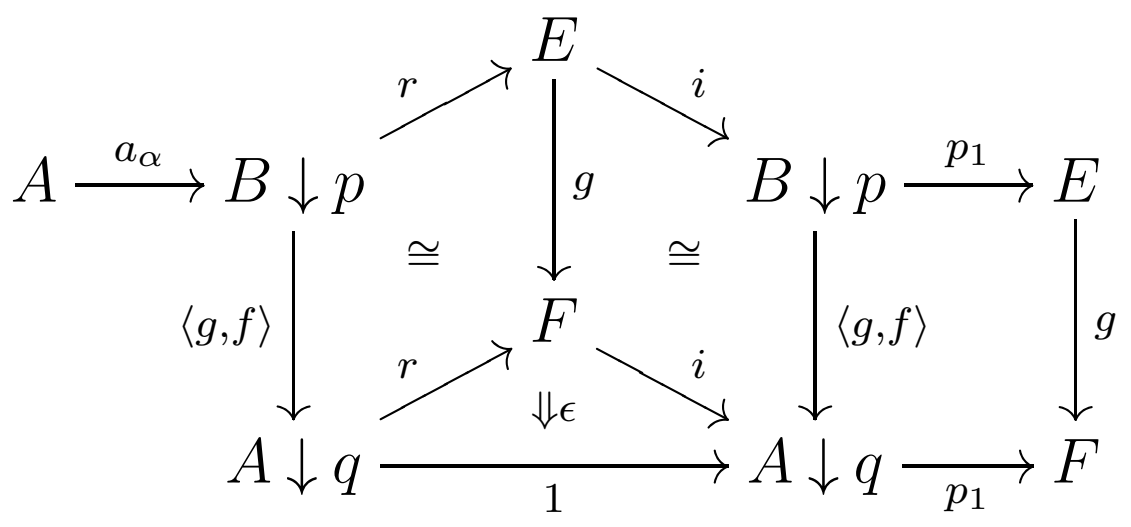

where the 2-cells in the squares are the inverse of (5.1.5) and its mate, which we assume is also an isomorphism. In this way we see that $g$ applied to the $p$-cartesian lift of $\alpha$ is isomorphic to a $q$-cartesian lift of $f\alpha \colon fb \Rightarrow qge$, which by Observation 4.1.6 suffices to prove that $g$ preserves all cartesian 2-cells. ☐

5.1.7. **Corollary.** *If a functor between cartesian fibrations*

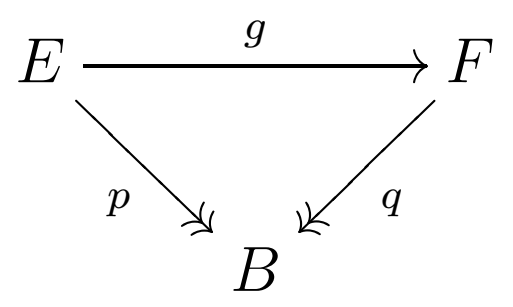

*admits a left adjoint $\ell \colon F \to E$ over $B$, then $g$ is a cartesian functor.*

*Proof.* By Corollary 3.6.6, the adjunction $\ell \dashv g$ over $B$ can be pulled back along $p_1 \colon B^2 \to B$ to define an adjunction

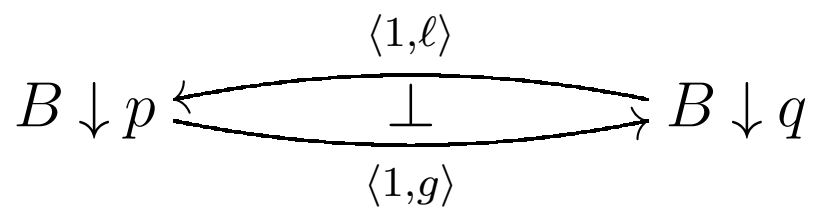

over $B^2$. Here we use Lemma 3.4.12 to identify iso-comma squares

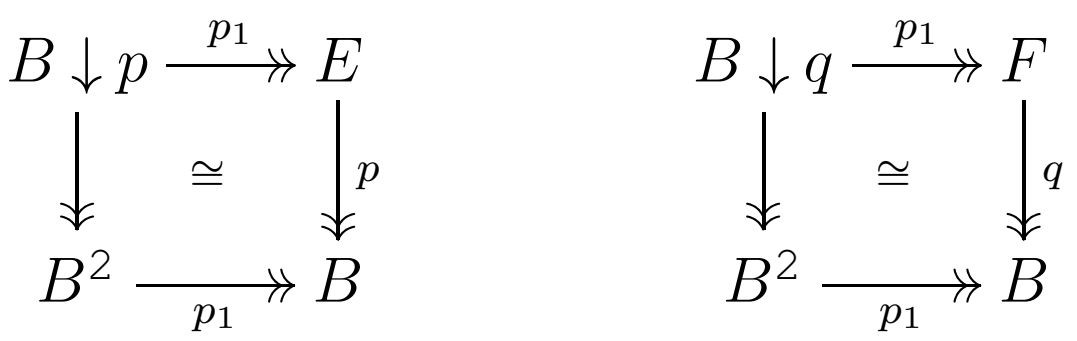

which by Lemma 3.5.6 proves that the pulled back adjunction has the form displayed.

By a standard 2-categorical result, the mate of (5.1.5) with respect to the right adjoints $i \dashv r$ is an isomorphism if and only if the mate with respect to the left adjoints $\ell \dashv g$ and $\langle 1, \ell\rangle \dashv \langle 1, f\rangle$ is an isomorphism. The latter is the case, because the left adjoint $\ell$ lies over $B$. Thus Theorem 5.1.4 implies that the right adjoint $g$ is cartesian. ☐

Using this result, we can prove some further stability properties of cartesian 2-cells.

5.1.8. **Lemma.** *Let $p \colon E \twoheadrightarrow B$ be a cartesian fibration. If $\chi \colon e' \Rightarrow e$ and $\chi' \colon e'' \Rightarrow e'$ are both cartesian for $p$, then so is $\chi \cdot \chi' \colon e'' \Rightarrow e$.*

*Proof.* By 1-cell induction, the 2-cell $p\chi\colon pe' \Rightarrow pe\colon A \to B$ may be regarded as a functor $a_{p\chi}\colon A \to B \downarrow p$. By Observation 4.1.14, its cartesian lift is $\chi_{p\chi} := p_1 \epsilon a_{p\chi}$, where $\epsilon$ is the counit of the adjunction $i \dashv r$ of Theorem 4.1.10. Because $\chi$ is cartesian, Observation 4.1.3 implies that $\chi$ and $\chi_{p\chi}$ are isomorphic. For convenience, we replace $\chi$ by $\chi_{p\chi}$ and redefine $\chi'$ so as to absorb the isomorphism. This modification does not change the 2-cell $\chi \cdot \chi'$, which we desire to show is cartesian.

By 2-cell induction, the diagram

$$\begin{array}{ccc} pe'' & \overset{p\chi'}{\Longrightarrow} & pe' \\ {\scriptstyle p(\chi\cdot\chi')}\Downarrow & & \Downarrow{\scriptstyle p\chi} \\ pe & \underset{\mathrm{id}}{\Longrightarrow} & pe \end{array}$$

defines a 2-cell $a_{p(\chi\cdot\chi')} \Rightarrow a_{p\chi}\colon A \to B \downarrow p$. Observation 4.1.17, applied in the context of Example 4.1.18, tells us that this 2-cell is cartesian for $p_0\colon B{\downarrow}q \to B$, because its codomain component is an identity. By Corollary 5.1.7, the right adjoint $r\colon B \downarrow q \to E$ carries $p_0$-cartesian 2-cells to $p$-cartesian 2-cells. Thus, we obtain a $p$-cartesian 2-cell $\gamma\colon ra_{p(\chi\cdot\chi')} \Rightarrow ra_{p\chi} = e'$ with $p\gamma = p\chi'$. Middle four interchange applied to the 2-cells

$$A \underset{a_{p\chi}}{\overset{a_{p(\chi\cdot\chi')}}{\Downarrow a_{(p\chi',\mathrm{id})}}} B \downarrow p \underset{p_1}{\overset{r}{\Downarrow p_1\epsilon}} E$$

provides us with a commuting diagram of 2-cells

$$\begin{array}{ccc} ra_{p(\chi\cdot\chi')} & \overset{\gamma}{\Longrightarrow} & ra_{p\chi} = e' \\ {\scriptstyle \chi_{p(\chi\cdot\chi')}}\Downarrow & & \Downarrow{\scriptstyle \chi_{p\chi}=\chi} \\ e & \underset{\mathrm{id}}{\Longrightarrow} & e \end{array}$$

whose left-hand side is a cartesian lift of $p(\chi \cdot \chi')$; in particular $\chi \cdot \gamma$ is cartesian. Because $\chi'$ is cartesian, it is isomorphic to $\gamma$. This tells us $\chi' \cdot \chi \cong \chi \cdot \gamma$ is weakly $p$-cartesian. □

5.1.9. **Lemma.** *Let $p\colon E \twoheadrightarrow B$ be a cartesian fibration. If $\chi\colon e' \Rightarrow e$ and $\chi''\colon e'' \Rightarrow e$ are both cartesian for $p$, and if $\chi'' = \chi \cdot \chi'$, then $\chi'\colon e'' \Rightarrow e'$ is also cartesian for $p$.*

*Proof.* Let $\psi\colon \bar{e} \Rightarrow e'$ denote a cartesian lift of $p\chi'$. As $p\psi = p\chi'$, there is some $\theta\colon e'' \Rightarrow \bar{e}$ with $\psi\cdot\theta = \chi'$ and $p\theta = \mathrm{id}$. Composing on the left with $\chi$ we see that $\chi\cdot\psi\cdot\theta = \chi\cdot\chi'\cdot\theta = \chi''\cdot\theta$. By Lemma 5.1.8, $\chi''$ and $\chi \cdot \psi$ define two cartesian lifts of $p\chi''$. Observation 4.1.3 now implies that $\theta$ is an isomorphism, so we conclude by Observation 4.1.2 that $\chi' = \psi \cdot \theta$ is cartesian. □

5.2. **Pullback stability.** Recall Definition 3.5.7, which defines the pullback of an isofibration $p\colon E \twoheadrightarrow B$ along a functor $f\colon A \to B$. Corollary 3.6.7 observes that this notion is well-defined up to equivalence in $\mathcal{C}/A$. We now demonstrate that pullbacks of (co-)cartesian fibrations are again (co-)cartesian. The analogous result for groupoidal (co-)cartesian fibrations will follow as an easy corollary.

5.2.1. **Proposition** (pullbacks of cartesian fibrations). *In an abstract homotopy 2-category $\mathcal{C}$, suppose*

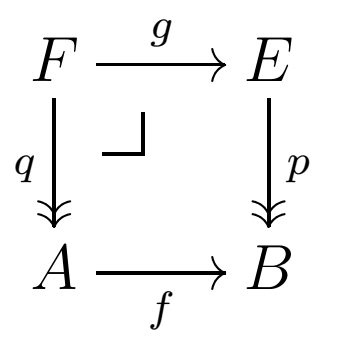

*is a pullback and* $p\colon E \twoheadrightarrow B$ *is a cartesian fibration. Then* $q\colon F \twoheadrightarrow A$ *is a cartesian fibration. Moreover, a 2-cell* $\chi\colon x' \Rightarrow x\colon X \to F$ *is cartesian for* $q$ *if and only if* $g\chi\colon gx' \Rightarrow gx\colon X \to E$ *is cartesian for* $p$*. In particular, the pullback square defines a cartesian functor from* $q$ *to* $p$*.*

*Proof.* Applying Theorem 4.1.10, it will suffice to construct a right adjoint to the induced map $\bar{i}\colon F \to A \downarrow q$ over $A$, which we do using the data of the adjunction

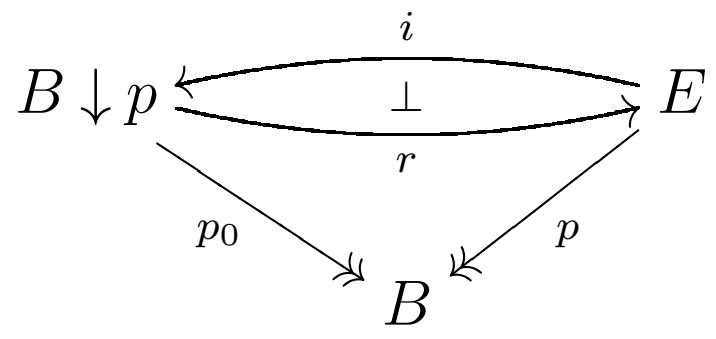

The pair $(f, g)$ induces a functor $\langle g, f\rangle\colon A \downarrow q \to B \downarrow p$ by 1-cell induction applied to the diagram

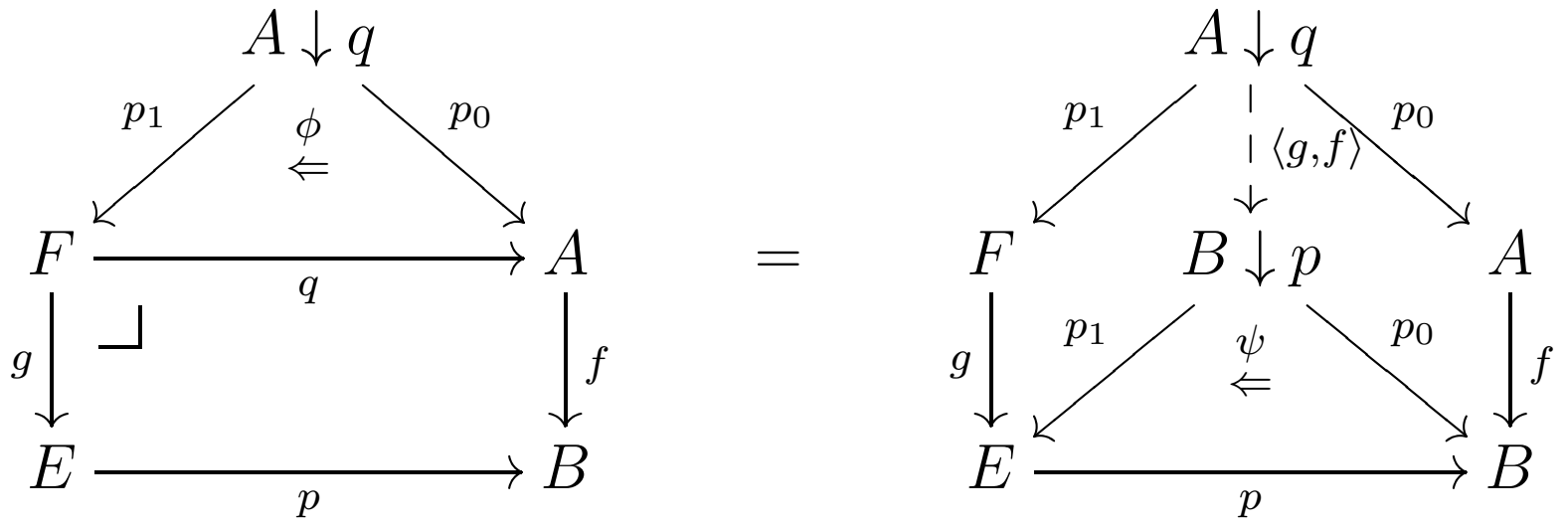

Using this, we define the right adjoint to $\bar{i}$ using 1-cell induction into the pullback $F$ applied to the identity 2-cell in the diagram

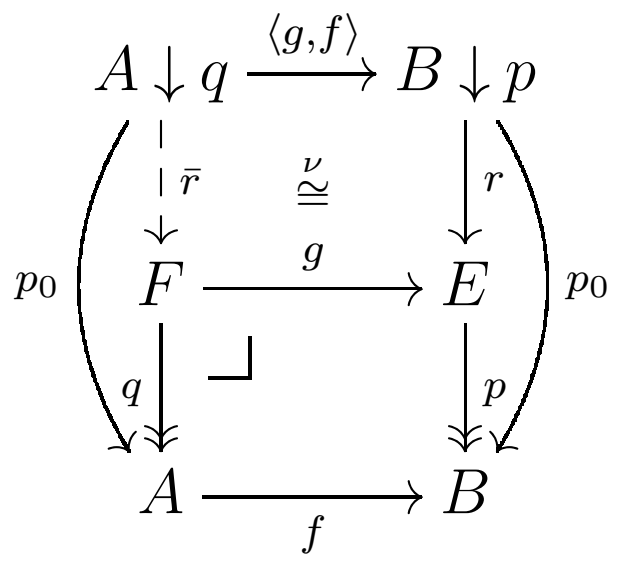

Uniqueness of 1-cell induction into $B \downarrow p$ implies that the square

$$\begin{array}{ccc} F & \xrightarrow{\quad g \quad} & E \\ {\scriptstyle \bar{i}}\big\downarrow & \cong\gamma & \big\downarrow{\scriptstyle i} \\ A \downarrow q & \xrightarrow[\langle g,f\rangle]{} & B \downarrow p \end{array}$$

commutes up to an isomorphic 2-cell over $E$ and $B$.

The unit 2-cell $\bar{\eta}\colon \mathrm{id}_F \Rightarrow \bar{r}\bar{i}$ is defined by 2-cell induction from the pair $q\bar{\eta} := \mathrm{id}_q$ and $g\bar{\eta} := g \overset{\eta}{\Rightarrow} rig \overset{r\gamma}{\Rightarrow} r\langle g, f\rangle\bar{i} \overset{\nu\bar{i}}{\Rightarrow} g\bar{r}\bar{i}$. By construction, $\bar{\eta}$ is fibered over $A$, and 2-cell conservativity implies that it is an isomorphism.

Applying Lemma I.4.1.2, recalled in Observation 3.6.14, in $\mathcal{C}/A$, to prove that $\bar{i} \dashv \bar{r}$ with unit $\bar{\eta}$ it suffices to exhibit a fibered 2-cell $\bar{\epsilon}\colon \bar{i}\bar{r} \Rightarrow \mathrm{id}_{A\downarrow q}$ that whiskers with $\bar{i}$ and with $\bar{r}$ to isomorphisms. We define $\bar{\epsilon}$ using 2-cell induction applied to $\mathrm{id}_{p_0}$, which ensures that the induced 2-cell lies over $A$, and a 2-cell $\tau\colon p_1\bar{i}\bar{r} = \bar{r} \Rightarrow p_1\colon A \downarrow q \to F$ that is itself defined by 2-cell induction. We define $\tau$ so that $q\tau := \phi\colon p_0 \Rightarrow qp_1$ and $g\tau := g\bar{r} \overset{\nu}{\Rightarrow} r\langle g, f\rangle \overset{\chi\langle g,f\rangle}{\Rightarrow} p_1\langle g, f\rangle = gp_1$, where $\chi$ is a $p$-cartesian lift of the $\psi$ used in the proof of Theorem 4.1.10.(i)⇒(ii) to define the functor $r$.

The proof that $\bar{\epsilon}\bar{i}$ is an isomorphism is an easy consequence of 2-cell conservativity: it suffices to show that $\tau\bar{i}$ is an isomorphism and this follows from the fact that $\phi\bar{i} = \mathrm{id}_q$ and $\chi\langle g, f\rangle\bar{i}$ defines a $p$-cartesian lift of $\mathrm{id}_{pg}$ and is thus an isomorphism.

The proof that $\bar{r}\bar{\epsilon}$ is an isomorphism also uses 2-cell conservativity. We have $q\bar{r}\bar{\epsilon} = p_0\bar{\epsilon} = \mathrm{id}_{p_0}$, so it remains only to show that $g\bar{r}\bar{\epsilon}$ is an isomorphism. For this, consider the middle four interchange square

$$\begin{array}{ccc} g\bar{r}\bar{i}\bar{r} & \overset{g\tau\bar{i}\bar{r}}{\Longrightarrow} & gp_1\bar{i}\bar{r} \\ {\scriptstyle g\bar{r}\epsilon}\Big\Downarrow & & \Big\Downarrow{\scriptstyle gp_1\bar{\epsilon} = g\tau} \\ g\bar{r} & \underset{g\tau}{\Longrightarrow} & gp_1 \end{array}$$

Recall $g\tau$ is isomorphic to a $p$-cartesian cell, and hence is itself $p$-cartesian. We have just argued that $\tau\bar{i}$ is an isomorphism, so the commutative square implies that $(g\bar{r}\bar{\epsilon}) \cdot (g\tau\bar{i}\bar{r})^{-1}$ is an automorphism of the $p$-cartesian cell $g\tau$ that maps, upon application of $p$, to an identity. Definition 4.1.1.(ii) now implies that this composite, and hence $g\bar{r}\bar{\epsilon}$ is an isomorphism, as required.

It remains to argue that $g$ preserves and reflects cartesian cells. Preservation follows from the characterization of cartesian functors presented in Theorem 5.1.4.(ii). The definition of $\bar{\eta}$ ensures that the mate of the canonical 2-cell $\gamma$ is an identity.

Now suppose given a generic 2-cell $\lambda\colon x \Rightarrow y\colon X \to F$ and consider the factorization $\lambda = \chi_{q\lambda}\zeta$ of $\lambda$ through the $q$-cartesian lift of $q\lambda$. We have just shown that $g\chi_{q\lambda}$ is $p$-cartesian. If $g\lambda$ is also $p$-cartesian, then Observation 4.1.2 implies that then $g\zeta$ is an isomorphism. By construction $q\zeta$ is an identity, so 2-cell conservativity of pullbacks now tells us that $\zeta$ is an isomorphism, which implies that $\lambda$ is $q$-cartesian. Thus $g$ creates cartesian 2-cells. □

5.2.2. **Corollary.** *Suppose*

$$\begin{array}{ccc} F & \xrightarrow{g} & E \\ {\scriptstyle q}\downarrow & \lrcorner & \downarrow{\scriptstyle p} \\ A & \xrightarrow[f]{} & B \end{array}$$

*is a pullback square and $p$ is a groupoidal cartesian fibration. Then $q$ is a groupoidal cartesian fibration.*

*Proof.* It remains only to prove groupoidalness. Given a 2-cell $\alpha$ with codomain $F$ so that $q\alpha$ is an identity, then $pg\alpha$ is an isomorphism, whence $g\alpha$ is an isomorphism by Lemma 4.2.2, which demonstrates the conservativity of the groupoidal isofibration $p$, whence $\alpha$ is an isomorphism by 2-cell conservativity for iso-commas. □

5.2.3. **Example.** Example 4.1.16 shows that the domain-projection functor $p_0 \colon A^2 \twoheadrightarrow A$ is a cartesian fibration. For any $f \colon B \to A$, Lemma 3.4.12 and Example 3.5.8 provides a pullback square

$$\begin{array}{ccc} f \downarrow A & \longrightarrow & A^2 \\ {\scriptstyle p_0}\downarrow & \lrcorner & \downarrow{\scriptstyle p_0} \\ B & \xrightarrow[f]{} & A \end{array}$$

Thus, we see that $p_0 \colon f \downarrow A \twoheadrightarrow B$ is a pullback of $p_0 \colon A^2 \twoheadrightarrow A$ and so we conclude from Proposition 5.2.1 that the pullback again defines a cartesian fibration.

For any $g \colon C \to A$, Example 4.1.18 demonstrates that $p_0 \colon A \downarrow g \twoheadrightarrow A$ is also cartesian. From this, Lemma 3.4.12, Example 3.5.8 and Proposition 5.2.1 it follows that the pullback

$$\begin{array}{ccc} f \downarrow g & \longrightarrow & A \downarrow g \\ {\scriptstyle p_0}\downarrow & \lrcorner & \downarrow{\scriptstyle p_0} \\ B & \xrightarrow[f]{} & A \end{array}$$

is also cartesian. This completes the proof of Example 4.1.18. A dual argument shows that the codomain projection functors from commas $p_1 \colon f \downarrow g \twoheadrightarrow C$ are cocartesian fibrations.

## 6. The Yoneda lemma

Let $\mathcal{K}$ denote an $\infty$-cosmos and consider a cartesian fibration $p \colon E \twoheadrightarrow B$ and a point $b \colon 1 \to B$, whose domain is the biterminal object described in Observation 4.2.10. Write $\mathrm{map}_B$ for the hom quasi-category in $\mathcal{K}/B$ defined in Example 2.1.11; as usual, we allow the domain to be an object such as $b \colon 1 \to B$ that is not necessarily an isofibration.

The point

$$\begin{array}{ccc} 1 & \xrightarrow{t} & B \downarrow b \\ & {\scriptstyle b}\searrow \quad \swarrow{\scriptstyle p_0} & \\ & B & \end{array}$$

induced in the comma object $B{\downarrow}b$ is terminal (see Lemma I.4.1.6). Evaluation at the terminal object $t$ induces a functor $\mathrm{map}_B(p_0, p) \to \mathrm{map}_B(b, p)$. Let $\mathrm{map}_B^{\mathrm{cart}}(p_0, p) \subset \mathrm{map}_B(p_0, p)$ denote the full sub quasi-category whose vertices are the cartesian functors.

6.0.1. **Theorem** (Yoneda lemma)**.** *Restriction along* $t\colon 1 \to B \downarrow b$ *induces an equivalence of quasi-categories*

$$\mathrm{map}_B^{\mathrm{cart}}(p_0\colon B\downarrow b \twoheadrightarrow B, p\colon E \twoheadrightarrow B) \simeq \mathrm{map}_B(b\colon 1 \to B, p\colon E \twoheadrightarrow B).$$

We prove this result along the following lines. The terminal object $t\colon 1 \to B \downarrow b$ defines an adjunction

$$1 \underset{t}{\overset{!}{\underset{\perp}{\rightleftarrows}}} B \downarrow b \tag{6.0.2}$$

not over $B$ but in the *lax comma 2-category* $\mathcal{K}_2/\!\!/B$, introduced below. We use this data to construct a right adjoint

$$\mathrm{map}_B(b, p) \underset{R}{\overset{\mathrm{map}_B(t,p)}{\underset{\perp}{\rightleftarrows}}} \mathrm{map}_B(p_0, p)$$

that lands in the sub quasi-category $\mathrm{map}_B^{\mathrm{cart}}(p_0, p)$ and show that this pair of functors restrict to an adjoint equivalence of quasi-categories.

Because there is no 2-functor $\mathrm{map}_B(-, p)\colon (\mathcal{K}_2/\!\!/B)^{\mathrm{op}} \to \underline{\mathrm{qCat}}_2$, in order to derive the second adjunction from the first, we instead lift along an appropriate smothering 2-functor — the same strategy employed to prove Corollary 3.6.6. We introduce this smothering 2-functor in §6.1, which gives a general characterization of the weak 2-functoriality of the construction of the pullback of a cartesian fibration. In §6.2, we introduce the lax slice 2-category and prove Theorem 6.0.1 along the lines just sketched.

6.1. **Functoriality of pullbacks of cartesian fibrations.** The aim in this section is to establish a sort of weak 2-functoriality satisfied by the operation of pulling back a cartesian fibration, encoded by Proposition 6.1.6 below.

Let $\mathcal{C}$ be an abstract homotopy 2-category.

6.1.1. **Definition.** Let $\mathcal{C}^{\lrcorner}$ denote the 2-category whose objects are cospans whose right map is a cartesian fibration, whose 1-cells are diagrams of the form

$$\begin{array}{ccccc} A' & \xrightarrow{f'} & B' & \overset{p'}{\twoheadleftarrow} & E' \\ {\scriptstyle a}\downarrow & \Uparrow\phi & \downarrow{\scriptstyle b} & & \downarrow{\scriptstyle e} \\ A & \xrightarrow[f]{} & B & \underset{p}{\twoheadleftarrow} & E \end{array} \tag{6.1.2}$$

and whose 2-cells consist of a triple of 2-cells $\alpha\colon a \Rightarrow \bar{a}$, $\beta\colon b \Rightarrow \bar{b}$, $\epsilon\colon e \Rightarrow \bar{e}$ between the verticals so that $p\epsilon = \beta p'$ and $\bar{\phi} \cdot f\alpha = \beta f' \cdot \phi$.

Given any object $A \xrightarrow{f} B \overset{p}{\twoheadleftarrow} E$ in $\mathcal{C}^{\lrcorner}$, we may form a pullback

$$\begin{array}{ccc} E \times_B A & \xrightarrow{g} & E \\ {\scriptstyle q}\downarrow & \lrcorner & \downarrow{\scriptstyle p} \\ A & \xrightarrow[f]{} & B \end{array} \tag{6.1.3}$$

Proposition 5.2.1 demonstrates that the functor $q$ is again a cartesian fibration, which implies that the 2-functor $\mathcal{C}^{\square} \to \mathcal{C}^{\lrcorner}$ that we presently introduce is surjective on objects.

6.1.4. **Definition.** Let $\mathcal{C}^{\square}$ denote the 2-category whose objects are pullback squares (6.1.3) whose verticals are cartesian fibrations and whose 1-cells are cubes

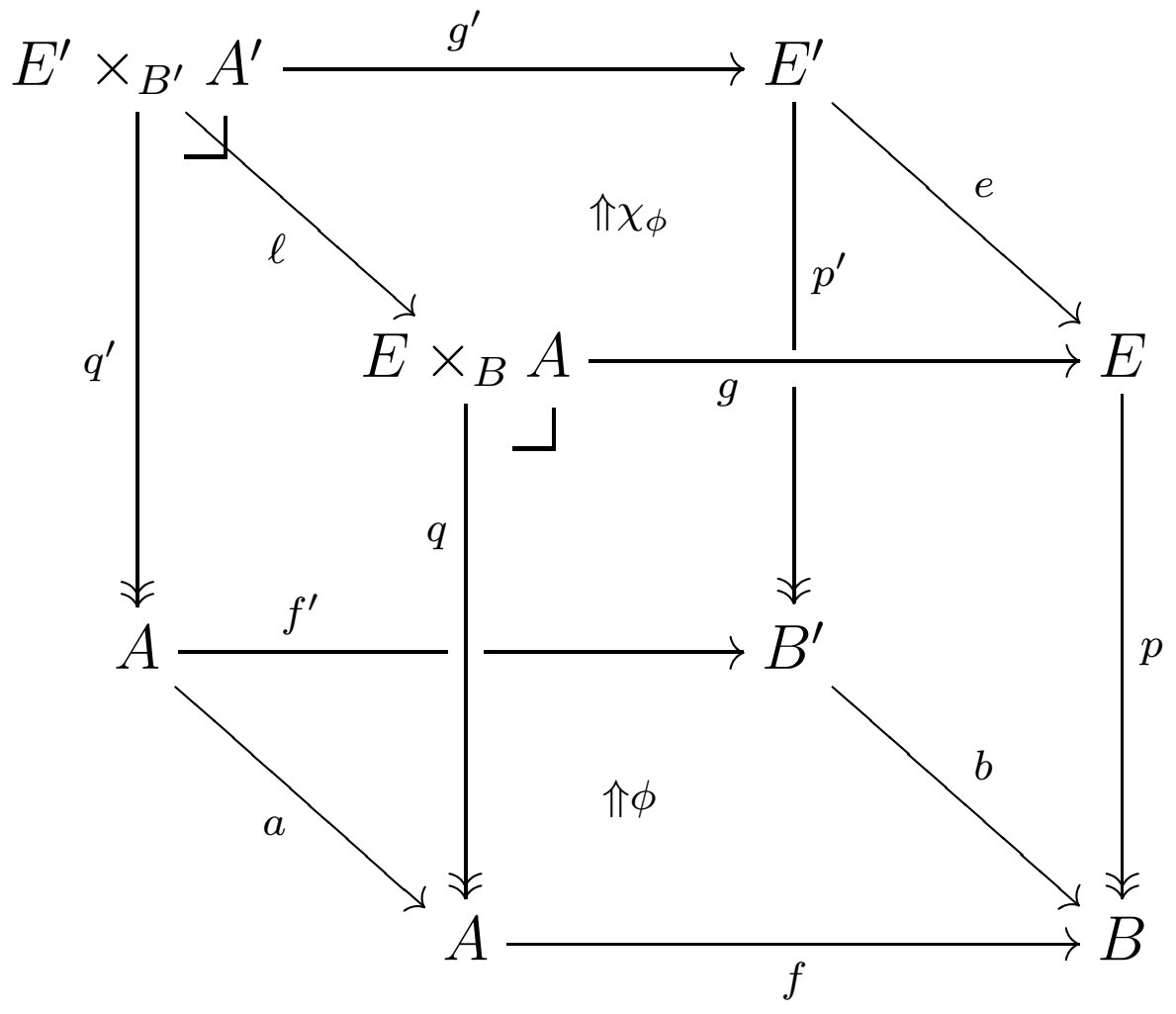

(6.1.5)

whose vertical faces commute and with $\chi_\phi \colon g\ell \Rightarrow eg'$ a $p$-cartesian lift of $\phi$. A 2-cell consists of a quadruple $\alpha\colon a \Rightarrow \bar{a}$, $\beta\colon b \Rightarrow \bar{b}$, $\epsilon\colon e \Rightarrow \bar{e}$, $\gamma\colon \ell \Rightarrow \bar{\ell}$ in which $\epsilon$ and $\gamma$ are, respectively, lifts of $\beta p'$ and $\alpha q'$ and so that $\bar{\phi} \cdot f\alpha = \beta f' \cdot \phi$ and $\chi_{\bar{\phi}} \cdot g\gamma = \epsilon g' \cdot \chi_\phi$.

We will make use of the pair of forgetful 2-functors $\mathcal{C} \leftarrow \mathcal{C}^{\square} \to \mathcal{C}^{\lrcorner}$ that project onto the pullback and the cospan, respectively.

6.1.6. **Proposition.** *$\mathcal{C}^{\square} \to \mathcal{C}^{\lrcorner}$ is a smothering 2-functor.*

*Proof.* Recall a smothering 2-functor is surjective on objects and locally smothering. Proposition 5.2.1 implies the former, so it remains to show that $\mathcal{C}^{\square} \to \mathcal{C}^{\lrcorner}$ is locally surjective on 1-cells, locally full, and conservative for 2-cells.

For local surjectivity on 1-cells, we must produce a diagram (6.1.5) expanding the data of (6.1.2) and a pair of chosen pullbacks. To begin, take $\chi_\phi$ to be any $p$-cartesian lift of $\phi q'\colon faq' \Rightarrow peg'$. By 1-cell induction, its domain factors, up to an isomorphism projecting along $p$ to an identity, as $g\ell$ for some functor $\ell$. By Observation 4.1.2, we can absorb this isomorphism into the $p$-cartesian cell $\chi_\phi$, so that it's domain equals $g\ell$ and the new $\chi_\phi\colon g\ell \Rightarrow eg'$ remains a $p$-cartesian lift of $\phi q'\colon faq' \Rightarrow peg'$. This completes the construction of (6.1.5).

For local fullness, suppose given a triple of 2-cells $\alpha\colon a \Rightarrow \bar{a}$, $\beta\colon b \Rightarrow \bar{b}$, $\epsilon\colon e \Rightarrow \bar{e}$ as in Definition 6.1.1. Define a 2-cell $\psi\colon g\ell \Rightarrow g\bar{\ell}$ by using the induction property of the cartesian 2-cell $\chi_{\bar{\phi}}\colon g\bar{\ell} \Rightarrow \bar{e}g'$ applied to the 2-cell $\epsilon g' \cdot \chi_\phi$ and the factorization of its whiskered composite with $p$ as $\bar{\phi}q' \cdot f\alpha q'$. By construction $p\psi = f\alpha q'$ so this pair induces a 2-cell $\gamma\colon \ell \Rightarrow \bar{\ell}$ projecting to $\alpha q'$ and $\psi$. The first condition tells us that $\gamma$ is a lift along $q$ of $\alpha q'$, and the second tells us that the pasted composite of $\chi_{\bar{\phi}}$ with $\gamma$ is the pasted composite of $\epsilon$ with $\chi_\phi$. In summary, $(\alpha, \beta, \epsilon, \gamma)$ is a 2-cell in $\mathcal{C}^\square$.

For 2-cell conservativity, suppose now that $\alpha$, $\beta$, and $\epsilon$ are isomorphisms. By the conservativity property for pullbacks, $\gamma$ is an isomorphism if both $q\gamma = \alpha q'$ and $g\gamma = \psi$ are. Using the inverses of $\alpha$, $\beta$, and $\gamma$, and arguing as in the previous paragraph, we produce a 2-cell $\bar{\psi}\colon g\bar{\ell} \Rightarrow g\ell$ lifting $f\alpha^{-1}q'$. The defining equations of these cartesian lifts tell us that $\bar{\psi} \cdot \psi$ is an automorphism of $\chi_\phi$ projecting to an identity and $\psi \cdot \bar{\psi}$ is an automorphism of $\chi'_\phi$ projecting to an identity. Now conservativity of cartesian 2-cells tells us that both composites are isomorphisms, so $\psi$ is as well. By 2-cell conservativity of weak pullbacks, we conclude that $\gamma$ is an isomorphism. □

6.2. **Proof of the Yoneda lemma.** We now specialize the results of the previous section to the pullback squares (2.1.12) defining the hom quasi-categories $\mathrm{map}_B$ appearing in the statement of the Yoneda lemma. Example 3.5.9 and Proposition I.3.3.14 demonstrate that they are pullbacks in the homotopy 2-category $\underline{\mathrm{qCat}}_2$, and thus fall under the purview of Proposition 6.1.6.

6.2.1. **Definition** (lax slice 2-category)**.** Given a 2-category $\mathcal{C}$ and an object $B$, write $\mathcal{C}/\!\!/B$ for the *lax slice 2-category*, with objects $f\colon X \to B$, 1-cells

$$\begin{array}{ccc} X & \xrightarrow{\quad k \quad} & Y \\ & {\scriptstyle f}\searrow \overset{\alpha}{\Rightarrow} \swarrow{\scriptstyle g} & \\ & B & \end{array} \tag{6.2.2}$$

and 2-cells

$$\begin{array}{ccc} X & \overset{k'}{\underset{k}{\overset{\Uparrow}{\rightrightarrows}}} & Y \\ & {\scriptstyle f}\searrow \Rightarrow\alpha \swarrow{\scriptstyle g} & \\ & B & \end{array} \quad = \quad \begin{array}{ccc} X & \xrightarrow{\quad k' \quad} & Y \\ & {\scriptstyle f}\searrow \overset{\alpha'}{\Rightarrow} \swarrow{\scriptstyle g} & \\ & B & \end{array}$$

For the remainder of this section, fix an $\infty$-cosmos $\mathcal{K}$.

6.2.3. *Observation.* The data of (6.0.2) defines an adjunction in $\mathcal{K}_2/\!\!/B$. The right adjoint $t$ is fibered over $B$. The left adjoint 1-cell is the universal comma cone:

$$\begin{array}{ccc} B \downarrow b & \xrightarrow{\quad ! \quad} & 1 \\ & {\scriptstyle p_0}\searrow \overset{\phi}{\Rightarrow} \swarrow{\scriptstyle b} & \\ & B & \end{array} \tag{6.2.4}$$

The counit is an identity and the unit is defined by 2-cell induction from $p_0\eta = \phi$ and $p_1\eta = \mathrm{id}$, as is required to ensure that this construction defines a 2-cell in $\mathcal{K}_2/\!\!/B$. The proof that this data defines an adjunction is a special case of Lemma I.4.1.6.

6.2.5. *Observation.* Fixing a cartesian fibration $p\colon E \twoheadrightarrow B$ in $\mathcal{K}$, we define a 2-functor $(\mathcal{K}_2/\!\!/B)^{\mathrm{op}} \to \underline{\mathrm{qCat}}_2^{\lrcorner}$ that caries a 1-cell (6.2.2) to

$$\begin{array}{ccccc} \Delta^0 & \xrightarrow{g} & \mathrm{map}(Y,B) & \xleftarrow{\mathrm{map}(Y,p)} & \mathrm{map}(Y,E) \\ \| & \Uparrow\alpha & \downarrow{\scriptstyle \mathrm{map}(k,E)} & & \downarrow{\scriptstyle \mathrm{map}(k,E)} \\ \Delta^0 & \xrightarrow[f]{} & \mathrm{map}(X,B) & \xleftarrow[\mathrm{map}(X,p)]{} & \mathrm{map}(X,E) \end{array}$$

By Corollary 4.1.15 and the fact that all objects in $\mathcal{K}_2$ are cofibrant, the functor of quasi-categories $\mathrm{map}(X,p)\colon \mathrm{map}(X,E) \twoheadrightarrow \mathrm{map}(X,B)$ is a cartesian fibration. Note that the pullbacks of the top and bottom cospans are the quasi-categories $\mathrm{map}_B(g,p)$ and $\mathrm{map}_B(f,p)$.

6.2.6. **Lemma.** *Let $p\colon E \twoheadrightarrow B$ be a cartesian fibration in $\mathcal{K}$. A 2-cell $\chi\colon x \Rightarrow y\colon Q \to \mathrm{map}(X,E)$ is $\mathrm{map}(X,p)$-cartesian if and only if each of its components*

$$\Delta^0 \xrightarrow{q} Q \underset{x}{\overset{y}{\underset{\Uparrow\chi}{\rightrightarrows}}} \mathrm{map}(X,E)$$

*define $p$-cartesian 2-cells $\chi q\colon xq \Rightarrow yq\colon X \to E$ in $\mathcal{K}_2$.*

*Proof.* If $\chi$ is $\mathrm{map}(X,p)$-cartesian then so is $\chi q$, by stability of cartesian cells under restriction. Definition 4.1.1, applied to the $\mathrm{map}(X,p)$-cartesian 2-cell $\chi q$, is expressed entirely in reference to the functor

$$\hom(\Delta^0, \mathrm{map}(X,E)) \xrightarrow{\mathrm{map}(X,p)\circ -} \hom(\Delta^0, \mathrm{map}(X,B))$$

between hom-categories in the homotopy 2-category $\underline{\mathrm{qCat}}_2$. This functor is isomorphic to the functor

$$\hom(X,E) \xrightarrow{p\circ -} \hom(X,B)$$

between hom-categories in the homotopy 2-category $\mathcal{K}_2$, from which we conclude that $\chi q\colon xq \Rightarrow yq\colon X \to E$ is $p$-cartesian.

Conversely, the argument just given tells us that if $\chi q\colon xq \Rightarrow yq\colon X \to E$ is $p$-cartesian then $\chi q\colon xq \Rightarrow yq\colon \Delta^0 \to \mathrm{map}(X,E)$ is $\mathrm{map}(X,p)$-cartesian. To show that $\chi$ is $\mathrm{map}(X,p)$-cartesian, consider its factorization $\chi = \chi' \cdot \lambda$ through a $\mathrm{map}(X,p)$-cartesian lift $\chi'$ of its projection. By Observation 4.1.3, we know that each of the components $\lambda q$ are isomorphisms. Observation I.3.2.3, which tells us that pointwise isomorphisms in $\underline{\mathrm{qCat}}_2$ are isomorphisms, then tells us that $\lambda$ is an isomorphism. Hence $\chi$ is $\mathrm{map}(X,p)$-cartesian. □

The 2-functor of Observation 6.2.5 carries the adjunction described in Observation 6.2.3 to an adjunction in $\underline{\mathrm{qCat}}_2^{\lrcorner}$; this 2-functor transforms the left adjoint 1-cell (6.2.4) into the

right adjoint 1-cell

$$\begin{array}{ccccc} \Delta^0 & \xrightarrow{b} & \mathrm{map}(1,B) & \xleftarrow{\mathrm{map}(1,p)} & \mathrm{map}(1,E) \\ \| & \Uparrow\phi & \downarrow{\scriptstyle \mathrm{map}(!,B)} & & \downarrow{\scriptstyle \mathrm{map}(!,E)} \\ \Delta^0 & \xrightarrow[p_0]{} & \mathrm{map}(B\downarrow b,B) & \xleftarrow[\mathrm{map}(B\downarrow b,p)]{} & \mathrm{map}(B\downarrow b,E) \end{array} \tag{6.2.7}$$

By Proposition 6.1.6, the projection $\underline{\mathrm{qCat}}_2^{\square} \to \underline{\mathrm{qCat}}_2^{\lrcorner}$ is a smothering 2-functor, so we may apply Lemma I.4.5.2 to lift this data to an adjunction in $\underline{\mathrm{qCat}}_2^{\square}$ between the simplicial pullbacks (2.1.12) defining the sliced mapping quasi-categories. Projecting along the 2-functor $\underline{\mathrm{qCat}}_2^{\square} \to \underline{\mathrm{qCat}}_2$ that evaluates at the upper-left-hand vertex gives us an adjunction

$$\mathrm{map}_B(b,p) \underset{R}{\overset{\mathrm{map}_B(t,p)}{\rightleftarrows}} \mathrm{map}_B(p_0,p) \quad (\bot) \tag{6.2.8}$$

with the right adjoint 1-cell (6.2.7) lifting and then projecting to the functor $R$.

6.2.9. **Lemma.** *For each vertex*

$$\begin{array}{ccc} 1 & \xrightarrow{e} & E \\ & {\scriptstyle b}\searrow \quad \swarrow{\scriptstyle p} & \\ & B & \end{array}$$

*in* $\mathrm{map}_B(b,p)$, $Re\colon p_0 \to p$ *defines a cartesian functor in* $\mathcal{K}/B$.

*Proof.* The proof of Proposition 6.1.6 reveals that $R$ is defined to be the domain component of a $\mathrm{map}(B\downarrow b,p)$-cartesian 2-cell lifting $\phi$.

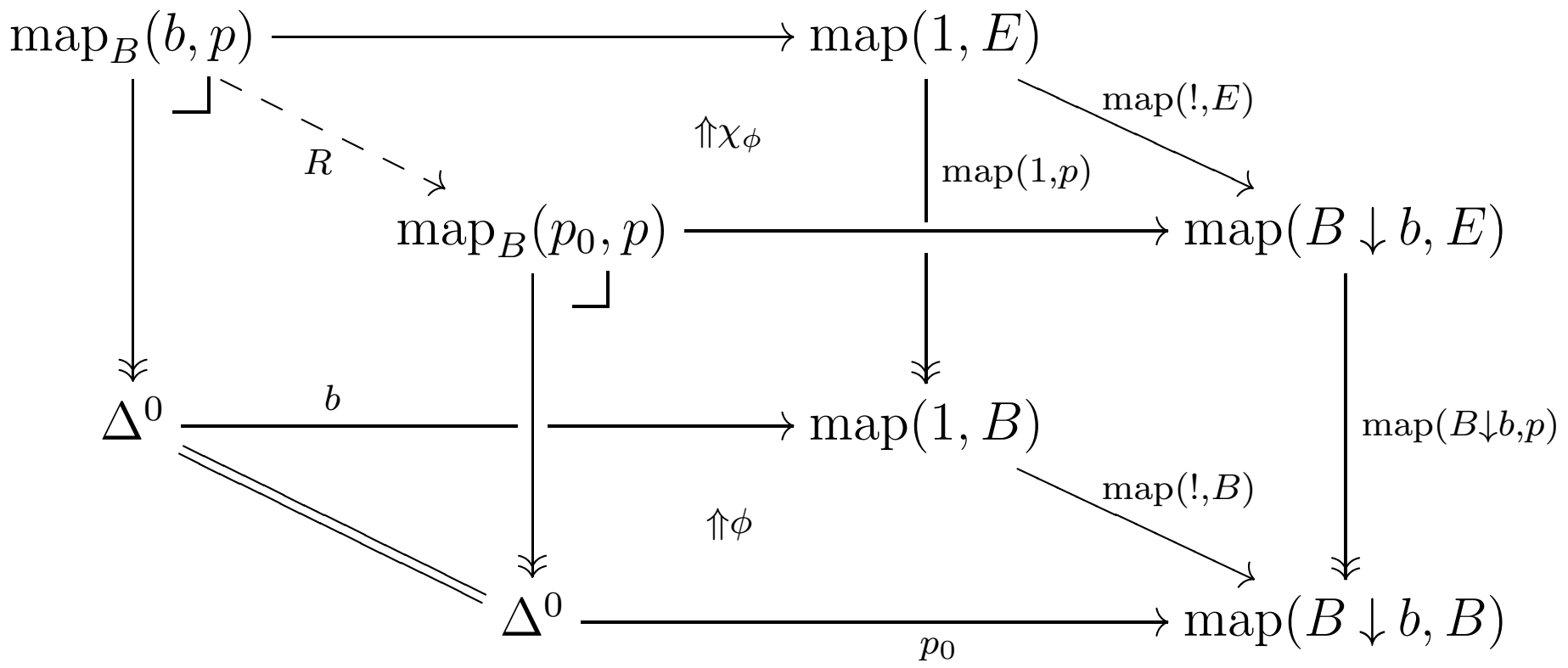

(6.2.10)

Applying Lemma 6.2.6, we see that the image $Re \in \mathrm{map}_B(p_0, p)_0$ of an object $e \in \mathrm{map}_B(b, p)_0$ is the domain component of a $p$-cartesian lift of the composite 2-cell:

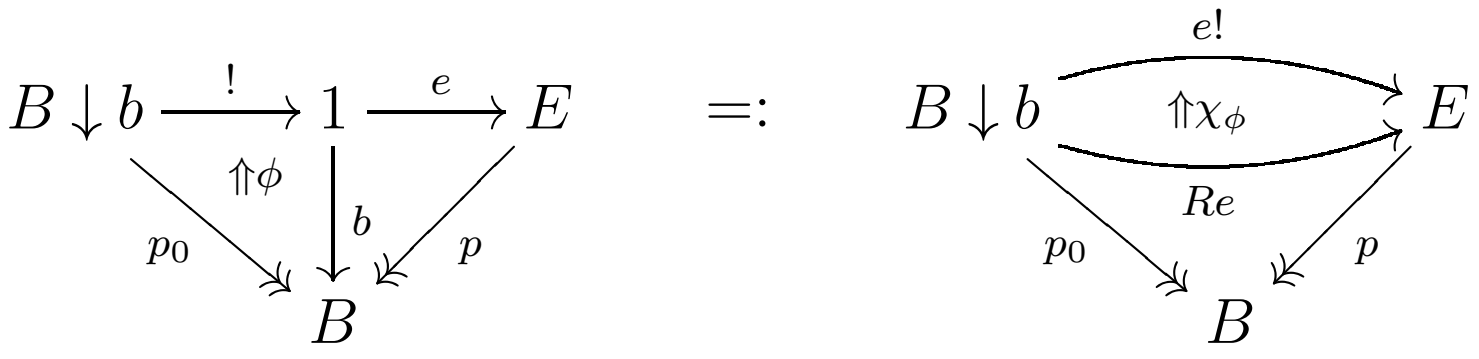

By Example 4.2.11 $p_0$ is a groupoidal cartesian fibration, and so Proposition 4.2.5 proves that every 2-cell $\tau\colon a \Rightarrow b\colon X \to B \downarrow b$ is $p_0$-cartesian. To show that $Re$ is a cartesian functor, we must show that $Re\tau$ is $p$-cartesian. By middle-four interchange, we have $\chi_\phi b \cdot Re\tau = e!\tau \cdot \chi_\phi a$ and the right-hand side, as the composite of a $p$-cartesian cell and an isomorphism, is $p$-cartesian. Thus, $Re\tau$ is $p$-cartesian by Lemma 5.1.9. □

Lemma 6.2.9 tells us that (6.2.8) restricts to an adjunction

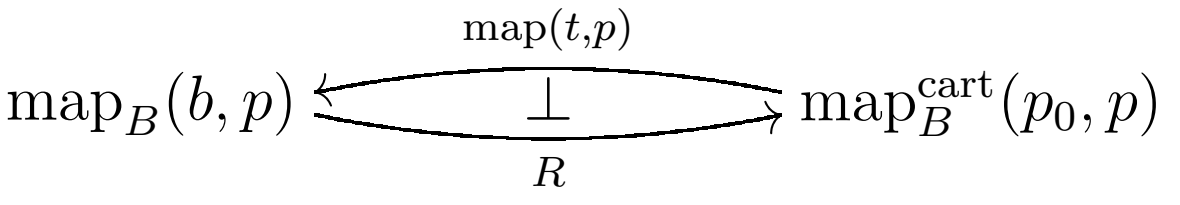

(6.2.11)

Moreover, by 2-cell conservativity of the smothering 2-functor $\underline{\mathrm{qCat}}_2^{\square} \to \underline{\mathrm{qCat}}_2^{\lrcorner}$, the counit is an isomorphism because the counit of the adjunction of 6.2.3 was. To prove that (6.2.11) is an adjoint equivalence, it remains only to show that the unit is an isomorphism.

6.2.12. **Lemma.** *If*

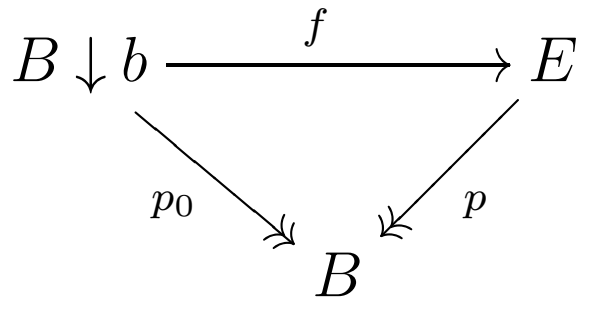

*is a cartesian functor then the component of the unit of the adjunction* (6.2.11) *at this object is an isomorphism.*

*Proof.* Let $\eta\colon \mathrm{id} \Rightarrow t!\colon B \downarrow b \to B \downarrow b$ denote the 2-cell in $\mathcal{K}_2$ defined in Observation 6.2.3 that serves as the unit of the adjunction in the lax slice 2-category $\mathcal{K}_2/\!\!/B$. Unpacking Observation 6.2.5 and the definition of $\underline{\mathrm{qCat}}_2^{\square}$, the unit $\hat{\eta}$ of (6.2.8) is defined to be a factorization of the restriction along $\mathrm{map}_B(p_0, p) \to \mathrm{map}(B\downarrow b, E)$ of the unit 2-cell $\mathrm{map}(\eta, E)$ of $\mathrm{map}(t, E) \dashv \mathrm{map}(!, E)$ through the restriction along $\mathrm{map}_B(t, p)$ of $\chi_\phi$, the $\mathrm{map}(B\downarrow b, p)$-cartesian lift (6.2.10).

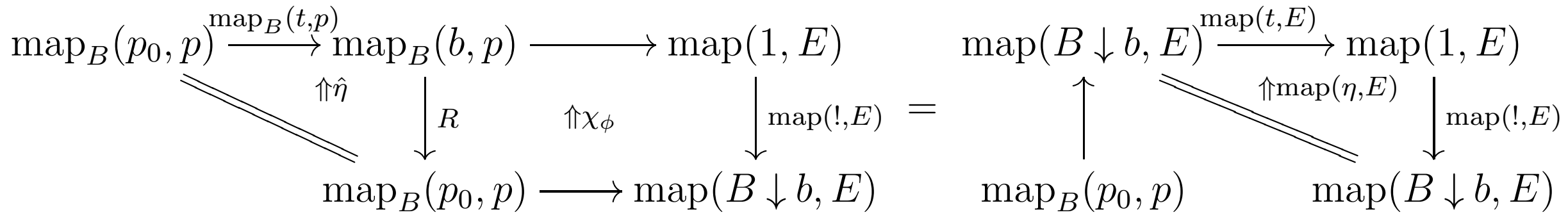

The condition $p_0\eta = \phi$ tells us that $\hat{\eta}$ projects to an identity.

The component of $\mathrm{map}(\eta, E)$ at a vertex $f \in \mathrm{map}_B(p_0, p)$ is $f\eta$. As $p_0 \colon B \downarrow b \twoheadrightarrow B$ is groupoidal, any 2-cell, such as $\eta$, with codomain $B \downarrow b$ is $p_0$-cartesian, so the hypothesis that $f$ is a cartesian functor implies that $f\eta$ is a $p$-cartesian 2-cell. By Lemma 6.2.6, the components of a $\mathrm{map}(B \downarrow b, p)$-cartesian 2-cell, obtained by evaluating at a vertex of $\mathrm{map}(B\downarrow b, E)$, are $p$-cartesian 2-cells in $\mathcal{K}_2$. So Observation 4.1.3 tells us that the component at $f$ of the factorization $\hat{\eta}$ of $f\eta$ through the $p$-cartesian lift of $\phi$ is an isomorphism, as claimed by the statement. □

Observation I.3.2.3 demonstrates that pointwise isomorphisms in $\underline{\mathrm{qCat}}_2$ are isomorphisms, so Lemma 6.2.12 tells us that $\mathrm{map}_B(t, p) \colon \mathrm{map}_B^{\mathrm{cart}}(p_0, p) \to \mathrm{map}_B(b, p)$ is the left adjoint part of an adjoint equivalence. This completes the proof of Theorem 6.0.1.

Recall that any functor whose target is a groupoidal cartesian fibration is a cartesian functor. This gives rise to a simplified statement of the Yoneda lemma in the special case of maps into a groupoidal cartesian fibration.

6.2.13. **Corollary.** *Let $p \colon E \twoheadrightarrow B$ be a groupoidal cartesian fibration in $\mathcal{K}$ and let $b \colon 1 \to B$ be a point. Then restriction along $t \colon 1 \to B \downarrow b$ induces an equivalence of quasi-categories*

$$\mathrm{map}_B(p_0 \colon B \downarrow b \twoheadrightarrow B, p \colon E \twoheadrightarrow B) \simeq \mathrm{map}_B(b \colon 1 \to B, p \colon E \twoheadrightarrow B).$$

## 7. Appendix: proof of Theorem 4.1.10

In this section we give a detailed proof of Theorem 4.1.10, whose statement we now recall.

4.1.10. **Theorem.** *If $p \colon E \twoheadrightarrow B$ is an isofibration in the 2-category $\mathcal{C}$ then the following are equivalent:*

*(i) The isofibration $p$ is a cartesian fibration,*

*(ii) The 1-cell $i \colon E \to B \downarrow p$ admits a right adjoint which is fibred over $B$. In other words, this condition states that there exists an adjunction*

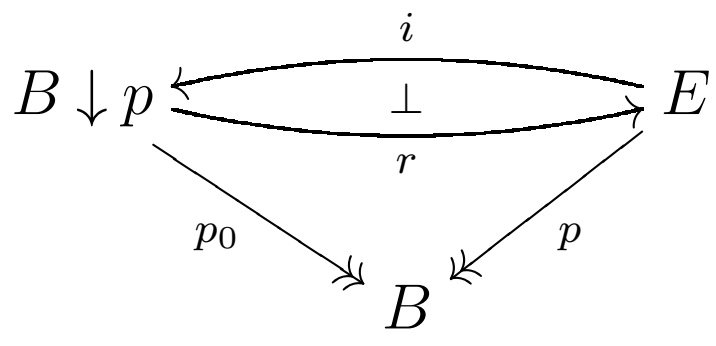

(7.0.1)

*in the slice 2-category $\mathcal{C}/B$.*

*(iii) The 1-cell $k \colon E^2 \to B \downarrow p$ is the left adjoint part of an adjunction*

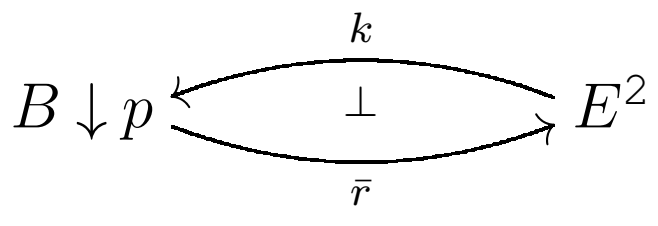

(7.0.2)

*in the 2-category $\mathcal{C}$ whose counit is an isomorphism.*

7.0.3. **Notation.** Recall Lemma 4.1.9, which establishes various adjunctions between the induced 1-cells defined in 4.1.8 and the projection isofibrations:

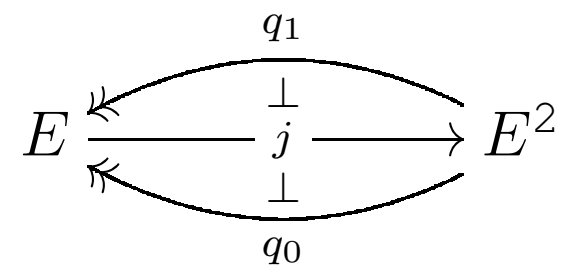
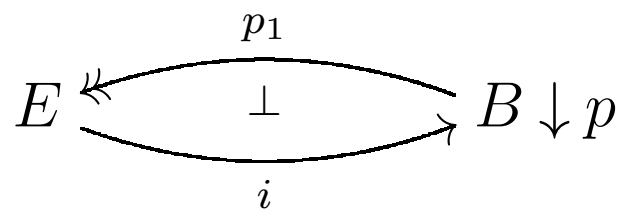

Here the counits of the adjunctions $q_1 \dashv j$ and $p_1 \dashv i$ are both identities, as is the unit of the adjunction $j \dashv q_0$.

We adopt the notation $\lambda\colon \mathrm{id}_{E^2} \Rightarrow jq_1$ and $\kappa\colon \mathrm{id}_{B\downarrow p} \Rightarrow ip_1$ for the units of the adjunctions $q_1 \dashv j$ and $p_1 \dashv i$ respectively, and observe that these enjoy the defining properties that $q_0\lambda = \psi$ and $q_1\lambda = \mathrm{id}_{q_1}$ and $p_0\kappa = \phi$ and $p_1\kappa = \mathrm{id}_{p_1}$. We also adopt the notation $\tau\colon jq_0 \Rightarrow \mathrm{id}_{E^2}$ for the counit of the adjunction $j \dashv q_0$, defined by 2-cell induction via the conditions $q_0\tau = \mathrm{id}_{q_0}$ and $q_1\tau = \psi$.

*Proof.* **(i) ⇒ (ii):** Assume that $p\colon E \twoheadrightarrow B$ is a cartesian fibration and start by taking a cartesian lift $\chi_\phi\colon \phi^*(p_1) \Rightarrow p_1$ of the defining 2-cell associated with the comma object $B\downarrow p$ and take the 1-cell $r\colon B\downarrow p \to E$ to be the domain $r := \phi^*(p_1)$ of that lift. Precomposing the cartesian 2-cell $\chi_\phi\colon r \Rightarrow p_1$ by the 1-cell $i\colon E \to B\downarrow p$ we obtain a 2-cell $\chi_\phi i\colon ri \Rightarrow p_1 i = \mathrm{id}_E$ which, by pre-composition stability, is again a cartesian 2-cell for $p$. Now $p\chi_\phi i = \phi i = \mathrm{id}_p$ so $\chi_\phi i$ is a cartesian lift of the identity 2-cell $\mathrm{id}_p\colon p \Rightarrow p$, and we know, from Observation 4.1.2, that this also has the cartesian lift $\mathrm{id}_{\mathrm{id}_E}\colon \mathrm{id}_E \Rightarrow \mathrm{id}_E$. So, applying Observation 4.1.6, we obtain an isomorphism $\eta\colon \mathrm{id}_E \cong ri$ for which $p\eta$ is an identity, that is to say that $\eta$ is actually an isomorphism between these 1-cells in $\mathcal{C}/B$. Notice that this argument also tells us that $\chi_\phi i$ is also itself an isomorphism (the inverse to $\eta$).

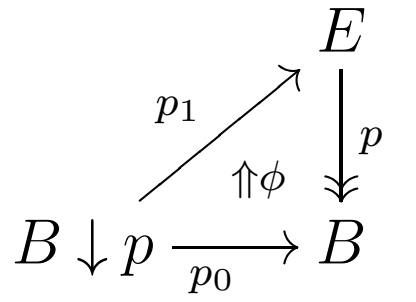

We can apply the 2-cell induction property of $B\downarrow p$ to induce a candidate counit $\epsilon\colon ir \Rightarrow \mathrm{id}_{B\downarrow p}$ with the defining properties $p_0\epsilon = \mathrm{id}_{p_0}$ and $p_1\epsilon = \chi_\phi$ on account of the equalities $p_0 ir = pr = p_0$, $p_1 ir = r$ and $\phi ir = \mathrm{id}_p r = \mathrm{id}_{p_0}$. In particular, the first of these defining equations tells us that this $\epsilon$ is actually a 2-cell in $\mathcal{C}/B$. To complete our proof, we verify that $\epsilon i$ and $r\epsilon$ are both isomorphisms and apply Observation 3.6.14, in $\mathcal{C}/B$, to complete the construction of the adjunction.

To that end, observe that $p_0\epsilon i = \mathrm{id}_{p_0} i = \mathrm{id}_p$, and $p_1\epsilon i = \chi_\phi i$, the second of which we have already seen to be an isomorphism, so we may apply the conservativity clause of the weak universal property of $B\downarrow p$ to infer that $\epsilon i$ is an isomorphism. Finally, to show that

$r\epsilon$ is also an isomorphism consider the naturality (middle four interchange) square

$$\begin{array}{ccc} rir & \xRightarrow{\chi_\phi ir} & p_1 ir \\ {\scriptstyle r\epsilon}\Downarrow & & \Downarrow{\scriptstyle p_1\epsilon} \\ r & \xRightarrow[\chi_\phi]{} & p_1 \end{array}$$

in which we know that that the upper horizontal 2-cell $\chi_\phi ir$ is an isomorphism which lies over the identity on $p_0$, because $\chi_\phi i$ is an isomorphism with $p\chi_\phi i = \mathrm{id}_p$, and we also know that the right hand vertical $p_1\epsilon = \chi_\phi$. Consequently this commutative square tells us that $\chi_\phi \cdot (r\epsilon \cdot (\chi_\phi ir)^{-1}) = \chi_\phi$ and, furthermore, $pr\epsilon = p_0\epsilon = \mathrm{id}_{p_0}$ so the 2-cell $r\epsilon \cdot (\chi_\phi ir)^{-1}$ lies over the identity on $p$. So we may apply the conservativity property of cartesian 2-cells, as discussed in Observation 4.1.3, to show that $r\epsilon \cdot (\chi_\phi ir)^{-1}$ is an isomorphism and thus that $r\epsilon$ is an isomorphism as required.

**(ii) ⇒ (iii):** Suppose that we have the fibred adjunction in $\mathcal{C}/B$ as depicted in (7.0.1) and adopt the notation $\eta\colon \mathrm{id}_E \cong ri$ and $\epsilon\colon ir \Rightarrow \mathrm{id}_{B\downarrow p}$ for its unit and counit respectively. Also define a 2-cell $\chi := p_1\epsilon\colon r = p_1 ir \Rightarrow p_1$ and construct a functor $\bar{r}\colon B\downarrow p \twoheadrightarrow E^2$ using the 1-cell induction property of $E^2$ and the defining equations $q_0\bar{r} = r$, $q_1\bar{r} = p_1$, and $\psi\bar{r} = \chi$:

$$\begin{array}{c} B\downarrow p \\ \downarrow{\scriptstyle \bar{r}} \\ E^2 \\ {\scriptstyle q_1}\left(\Leftarrow\psi\right){\scriptstyle q_0} \\ E \end{array} \quad = \quad \begin{array}{c} B\downarrow p \\ {\scriptstyle p_1}\left(\Leftarrow\chi\right){\scriptstyle r} \\ E \end{array}$$

Observe also that we have a naturality (middle four interchange) square

$$\begin{array}{ccc} p_0 ir & \xRightarrow{p_0\epsilon} & p_0 \\ {\scriptstyle \phi ir}\Downarrow & & \Downarrow{\scriptstyle \phi} \\ pp_1 ir & \xRightarrow[pp_1\epsilon]{} & pp_1 \end{array}$$

in which $p_0\epsilon = \mathrm{id}_{p_0}$, because $\epsilon$ is a 2-cell in $\mathcal{C}/B$, and $\phi ir = \mathrm{id}_r$, from the definition of the 1-cell $i$. It follows from the commutativity of this square that $p\chi = pp_1\epsilon = \phi$.

Now consider the composite $k\bar{r}\colon B\downarrow p \to B\downarrow p$ and observe that $p_0 k\bar{r} = pq_0\bar{r} = pr = p_0$, $p_1 k\bar{r} = q_1\bar{r} = p_1$, and $\phi k\bar{r} = p\psi\bar{r} = p\chi = \phi$. Of course, a second endo-1-cell on $B\downarrow p$ for which these particular equalities hold is the identity $\mathrm{id}_{B\downarrow p}$, so it follows by Lemma I.3.3.27, the essential uniqueness of induced 1-cells into comma objects, that there exists an isomorphism $\bar{\epsilon}\colon k\bar{r} \cong \mathrm{id}_{B\downarrow k}$ with $p_0\bar{\epsilon} = \mathrm{id}_{p_0}$ and $p_1\bar{\epsilon} = \mathrm{id}_{p_1}$.

To construct a 2-cell $\bar{\eta}\colon \mathrm{id}_{E^2} \Rightarrow \bar{r}k$ first observe that $q_0\bar{r}k = rk$, $q_1\bar{r}k = p_1 k = q_1$ and $\psi\bar{r}k = \chi k$ so we may use the 2-cell induction property of $E^2$ to construct a 2-cell $\bar{\eta}\colon \mathrm{id}_{E^2} \Rightarrow \bar{r}k$ from a pair of 2-cells $\hat{\mu}\colon q_0 \Rightarrow rk$ and $\mathrm{id}_{q_1}\colon q_1 \Rightarrow q_1$ which make the following

square commute:

$$\begin{array}{ccc} q_0 & \overset{\hat{\mu}}{\Longrightarrow} & rk \\ {\scriptstyle\psi}\Downarrow & & \Downarrow{\scriptstyle\chi k} \\ q_1 & = & p_1k \end{array} \tag{7.0.4}$$

To construct this 2-cell $\hat{\mu}$ observe that we may transpose the isomorphism $\mu\colon i \cong kj$, as discussed in 4.1.9, under the adjunction $j \dashv q_0$ to give a 2-cell $\tilde{\mu}\colon iq_0 \Rightarrow k$ and then transpose that under the adjunction $i \dashv r$ to give $\hat{\mu}\colon q_0 \Rightarrow rk$. More concretely, this may be expressed as a composite

$$q_0 \underset{\cong}{\overset{\eta q_0}{\Longrightarrow}} riq_0 \underset{\cong}{\overset{r\mu q_0}{\Longrightarrow}} rkjq_0 \overset{rk\tau}{\Longrightarrow} rk \tag{7.0.5}$$

where $\tau$ is the counit of $j \dashv q_0$ and $\eta$ is the unit of $i \dashv r$. With this definition the commutativity of the square in (7.0.4) reduces to the following computation:

| | |
|---|---|
| $\chi k \cdot \hat{\mu} = \chi k \cdot rk\tau \cdot r\mu q_0 \cdot \eta q_0$ | by definition of $\hat{\mu}$ |
| $= p_1k\tau \cdot \chi kjq_0 \cdot r\mu q_0 \cdot \eta q_0$ | middle four interchange left composite |
| $= p_1k\tau \cdot p_1\mu q_0 \cdot \chi iq_0 \cdot \eta q_0$ | middle four interchange centre composite |
| $= q_1\tau \cdot p_1\mu q_0 \cdot p_1\epsilon iq_0 \cdot p_1i\eta q_0$ | $p_1k = q_1$, $\chi = p_1\epsilon$, and $p_1i = \mathrm{id}_E$ |
| $= q_1\tau \cdot p_1\mu q_0 \cdot p_1(\epsilon i \cdot i\eta)q_0$ | factor $p_1$ and $q_0$ out of right hand composite |
| $= q_1\tau \cdot p_1\mu q_0$ | apply triangle identity for $i \dashv r$ |
| $= \psi$ | from 4.1.9 we have $q_1\tau = \psi$ and $p_1\mu$ an identity |

Consequently there is an induced 2-cell $\bar{\eta}$ as advertised with the defining properties that $q_0\bar{\eta} = \hat{\mu}$ and $q_1\bar{\eta} = \mathrm{id}_{q_1}$ so, by Observation 3.6.14, all that remains for us to prove in order to demonstrate the desired adjunction $k \dashv \bar{r}$ is that both of the whiskered 2-cells $k\bar{\eta}$ and $\bar{\eta}\bar{r}$ are isomorphisms.

The first of these is easy because $p_0k\bar{\eta} = pq_0\bar{\eta} = p\hat{\mu}$ and $p_1k\bar{\eta} = q_1\bar{\eta} = \mathrm{id}_{q_0}$ by the definitions of $k$ and $\bar{\eta}$. Furthermore, from the definition of $\hat{\mu}$ it is easily verified that $p\hat{\mu} = \mathrm{id}_{pq_0}$ using the facts that $r$, $\eta$ and $\mu$ are cells fibred over $B$, that is $pr = p_0$ and the 2-cells $p\eta$ and $p_0\mu$ are both identities, and that $\tau$ has defining property $q_0\tau = \mathrm{id}_{q_0}$. So it follows that we may use the conservativity clause of the weak 2-universal property of $B\downarrow p$ to conclude that $k\bar{\eta}$ is an isomorphism.

The proof of the second of these isomorphism properties is only slightly more involved, from the definition of $\bar{\eta}$ we have $q_0\bar{\eta}\bar{r} = \hat{\mu}\bar{r}$ and $q_1\bar{\eta}\bar{r} = \mathrm{id}_{q_1}\bar{r} = \mathrm{id}_{p_1}$. So to apply the conservativity clause of the weak 2-universal property of $E^2$ to prove that $\bar{\eta}\bar{r}$ is an isomorphism, as required, then all we need do is show that $\hat{\mu}\bar{r}$ is an isomorphism. However, consulting the definition of $\hat{\mu}$ in (7.0.5) we see that it is a composite of a pair of isomorphisms and a 2-cell $rk\tau\colon rkjq_0 \Rightarrow rk$, so to show that $\hat{\mu}\bar{r}$ is an isomorphism it is sufficient to verify that $rk\tau\bar{r}$ is an isomorphism.

To do so first define a map $\epsilon'\colon ir \Rightarrow \mathrm{id}_{B\downarrow p}$ as the following composite

$$ir \underset{\cong}{\overset{\mu r}{\Longrightarrow}} kjq_0\bar{r} \overset{k\tau\bar{r}}{\Longrightarrow} k\bar{r} \underset{\cong}{\overset{\bar{\epsilon}}{\Longrightarrow}} \mathrm{id}_{B\downarrow p}$$

and note that since the 2-cells at either end of this composite are isomorphisms it follows that $rk\tau\bar{r}$ is an isomorphism, as required, if and only if $r\epsilon'$ is an isomorphism. Now we have a naturality (middle four interchange) square

$$\begin{array}{ccc} irir & \overset{ir\epsilon}{\Longrightarrow} & ir \\ {\scriptstyle \epsilon' ir}\Downarrow & & \Downarrow{\scriptstyle \epsilon'} \\ ir & \underset{\epsilon}{\Longrightarrow} & \mathrm{id}_{B\downarrow p} \end{array}$$

and observe that $r\epsilon$ and $\epsilon i$ are both isomorphisms, since $i \dashv r$ is an adjunction with isomorphic unit, so in particular the horizontal map at the top of this square is an isomorphism. Furthermore we have $p_0 k\tau\bar{r}i = pq_0\tau\bar{r}i = p\,\mathrm{id}_{q_0}\bar{r}i = \mathrm{id}_{pq_0\bar{r}i} = \mathrm{id}_p$, where the first equality comes from the definition of $k$ and the second from the definition of $\tau$, and we have $p_1 k\tau\bar{r}i = q_1\tau\bar{r}i = \psi\bar{r}i = \chi i = p_1\epsilon i$ which is an isomorphism. So, by the conservativity clause of the weak 2-universal property of $B \downarrow p$, it follows that the 2-cell $k\tau\bar{r}i$ is an isomorphism and, on consulting the definition of $\epsilon'$, it is clear then that $\epsilon' i$ is also an isomorphism. In other words, we have found that the left hand vertical of the square above is an isomorphism and, consequently, that the other two 2-cells $\epsilon$ and $\epsilon'$ in that square are related by composition with the isomorphism $\epsilon' ir \cdot (ir\epsilon)^{-1}$. It follows, therefore, that since $r\epsilon$ is an isomorphism we can infer that $r\epsilon'$ is also an isomorphism and thus complete our demonstration of the adjunction $k \dashv \bar{r}$ as required.

**(iii) ⇒ (i):** Suppose now that we are given an adjunction as in (7.0.2) and adopt the notation $\bar{\epsilon}\colon k\bar{r} \cong \mathrm{id}_{B\downarrow p}$ and $\bar{\eta}\colon \mathrm{id}_{E^2} \Rightarrow \bar{r}k$ for its (isomorphic) counit and unit respectively. As a first step we make this adjunction into a fibred adjunction in $\mathcal{C}/(E \times B)$. The defining equations of $k$ tell us that $(p_1, p_0)k = (q_1, pq_0)$ so it is already a 1-cell in this slice as depicted, we need however to pick $\bar{r}$ and the unit and counit to also be cells in the slice. To do this consider the isomorphism $(p_1, p_0)\bar{\epsilon}\colon (q_1, pq_0)\bar{r} = (p_1, p_0)k\bar{r} \cong (p_1, p_0)$ which we may lift along the isofibration $(q_1, pq_0)$ to give an isomorphism $\kappa\colon \bar{r} \cong \bar{r}'$ with $(q_1, pq_0)\kappa = (p_1, p_0)\bar{\epsilon}$ and $(q_1, pq_0)\bar{r}' = (p_1, p_0)$, so in particular $\bar{r}'$ is a 1-cell in the slice 2-category. Now it is easily checked that $\bar{\epsilon}' := \bar{\epsilon} \cdot k\kappa^{-1}$ and $\bar{\eta}' := \kappa k \cdot \bar{\eta}$ satisfy the triangle identities and that they are both 2-cells in the slice 2-category, thus providing the fibred adjunction we seek. We drop the primes and simply assume that our original cells were selected to be fibred over $E \times B$.

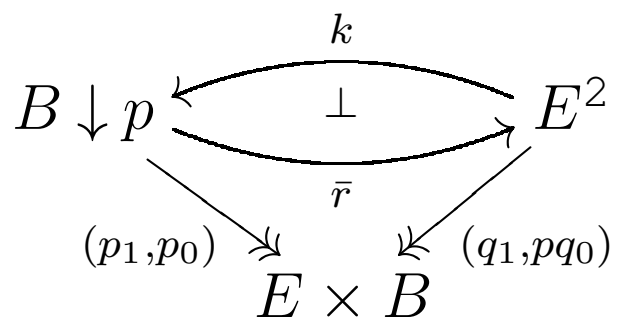

Now, suppose that we are given a 2-cell $\alpha\colon b \Rightarrow pe\colon A \to B$, as depicted in the diagram on the left of (4.1.5), which we are to lift to a cartesian 2-cell. To do this, start by observing that this data provides us with a comma cone over the arrow $p$ and so it induces a 1-cell $a_\alpha\colon A \to B \downarrow p$ with the defining properties $p_0 a_\alpha = b$, $p_1 a_\alpha = e$, and $\phi a_\alpha = \alpha$. Using this 1-cell we may define a 2-cell $\chi_\alpha := \psi\bar{r}a_\alpha$ whose codomain is $q_1\bar{r}a_\alpha = p_1 a_\alpha = e$, because $\bar{r}$ is a 1-cell in the slice over $E \times B$, and for which $p\chi_\alpha = p\psi\bar{r}a_\alpha = \phi k\bar{r}a_\alpha$ by the defining property of $k$. Now the 2-cell $\phi k\bar{r}$ features in the following middle four interchange square

$$\begin{array}{ccc} p_0k\bar{r} & \xRightarrow{p_0\bar{\epsilon}} & p_0 \\ {\scriptstyle \phi k\bar{r}}\Downarrow & & \Downarrow{\scriptstyle \phi} \\ pp_1k\bar{r} & \xRightarrow[pp_1\bar{\epsilon}]{} & pp_1 \end{array}$$

whose upper and lower horizontals are actually identities, because $\bar{\epsilon}$ is a 2-cell in the slice over $E \times B$, so it follows that $\phi k\bar{r} = \phi$ and thus that $p\chi_\alpha = \phi a_\alpha = \alpha$. In other words, we have demonstrated that $\chi_\alpha$ is indeed a lift of $\alpha$. We now show that it is cartesian for $p$. So let $e' \colon A \to E$ denote the domain of $\chi_\alpha$, which has $pe' = b$, and suppose that we are given a 2-cell $\gamma \colon e'' \Rightarrow e$ and a 2-cell $\beta \colon pe'' \Rightarrow b$ with $p\gamma = \alpha \cdot \beta$. Applying the 1-cell induction property of $E^2$ the 2-cell $\gamma$ induces a 1-cell $g \colon A \to E^2$ with defining properties $q_0g = e''$, $q_1g = e$, and $\psi g = \gamma$. It follows therefore that the composite $kg \colon A \to B \downarrow p$ has $p_0kg = pq_0g = pe''$, $p_1kg = q_1g = e$, and $\phi kg = p\psi g = p\gamma$, so we can regard $kg$ as being induced by the 1-cell induction property of $B \downarrow p$ with these defining properties. Of course, the 1-cell $a_\alpha \colon A \to B \downarrow p$ is induced in this way by $\alpha \colon b \Rightarrow pe$, so we may use 2-cell induction to construct a $\bar{\beta} \colon kg \Rightarrow a_\alpha$ with the defining properties $p_0\bar{\beta} = \beta$ and $p_1\bar{\beta} = \mathrm{id}_e$, simply because the required compatibility condition for this to happen reduces to the assumed equation $p\gamma = \alpha \cdot \beta$. Taking the adjoint transpose of $\bar{\beta}$ under the adjunction $k \dashv \bar{r}$ we obtain a corresponding 2-cell $\hat{\beta} \colon g \Rightarrow \bar{r}a_\alpha$ and these are related by the equation $\bar{\beta} = \bar{\epsilon}a_\alpha \cdot k\hat{\beta}$. In particular, since we chose $\bar{\epsilon}$ to be a 2-cell in the slice over $E \times B$ it follows that $p_0\bar{\beta} = p_0\bar{\epsilon}a_\alpha \cdot p_0k\hat{\beta} = pq_0\hat{\beta}$ and $p_1\bar{\beta} = p_1\bar{\epsilon}a_\alpha \cdot p_1k\hat{\beta} = q_1\hat{\beta}$. Now by horizontally composing $\hat{\beta}$ with $\psi$ we get the middle four interchange square

$$\begin{array}{ccc} q_0g & \xRightarrow{q_0\hat{\beta}} & q_0\bar{r}a_\alpha \\ {\scriptstyle \psi g}\Downarrow & & \Downarrow{\scriptstyle \psi\bar{r}a_\alpha} \\ q_1g & \xRightarrow[q_1\hat{\beta}]{} & q_1\bar{r}a_\alpha \end{array}$$

in which our various definitions give $\psi g = \gamma$, $q_1\hat{\beta} = p_1\bar{\beta} = \mathrm{id}_e$, and $\psi\bar{r}a_\alpha = \chi_\alpha$. So that square reduces to the equation $\chi_\alpha \cdot q_0\hat{\beta} = \gamma$ and, furthermore, we have also seen that $pq_0\hat{\beta} = p_0\bar{\beta} = \beta$. In other words we have shown that $q_0\hat{\beta}$ is the factorisation we seek to verify the induction part of the cartesian property of $\chi_\alpha$.

To check the required conservativity property of $\chi_\alpha$, suppose now that $\gamma \colon e' \Rightarrow e'$ is a 2-cell for which $p\gamma = \mathrm{id}_b$ and $\chi_\alpha \cdot \gamma = \chi_\alpha$. The equation $\chi_\alpha = \psi\bar{r}a_\gamma$ tells us that $\bar{r}a_\gamma$ is a 1-cell induced by $\chi_\alpha$ under the 1-cell induction property of $E^2$, so we may use $\gamma$ to induce a 2-cell $\bar{\gamma} \colon \bar{r}a_\alpha \Rightarrow \bar{r}a_\alpha$ with defining properties $q_0\bar{\gamma} = \gamma$ and $q_1\bar{\gamma} = \mathrm{id}_e$. We are assured of the existence of an induced 2-cell with these defining properties because the required compatibility condition reduces to the assumed equation $\chi_\alpha \cdot \gamma = \chi_\alpha$. Now horizontally composing $\bar{\gamma}$ and $\bar{\eta}$ we get the following middle four interchange square

$$\begin{array}{ccc} \bar{r}a_\alpha & \xRightarrow{\bar{\eta}\bar{r}a_\alpha} & \bar{r}k\bar{r}a_\alpha \\ {\scriptstyle\bar{\gamma}}\Downarrow & & \Downarrow{\scriptstyle\bar{r}k\bar{\gamma}} \\ \bar{r}a_\alpha & \xRightarrow[\bar{\eta}\bar{r}a_\alpha]{} & \bar{r}k\bar{r}a_\alpha \end{array}$$

in which we observe that $k\bar{\gamma}$ is an isomorphism, by the conservativity property of $B \downarrow p$ because $p_0 k\bar{\gamma} = pq_0\bar{\gamma} = p\gamma = \mathrm{id}_b$ and $p_1 k\bar{\gamma} = q_1\bar{\gamma} = \mathrm{id}_e$, as is $\bar{\eta}\bar{r}$, since the counit of $k \dashv \bar{r}$ is assumed to be an isomorphism. So this square expresses $\bar{\gamma}$ as a composite of three isomorphisms, so it follows that both it and $\gamma = q_0\bar{\gamma}$ are isomorphisms as hoped.

Finally all that remains is to show that the cartesian 2-cells for $p$ are preserved by precomposition by arbitrary 1-cells. To that end, suppose that $f\colon A' \to A$ is a 1-cell and consider the cartesian lifts $\chi_\alpha$ and $\chi_{\alpha f}$. Observe that the 1-cell $a_\alpha f$ satisfies the equations $p_0 a_\alpha f = bf$, $p_1 a_\alpha f = ef$ and $\phi a_\alpha f = \alpha f$, so as an induced 1-cell into $B{\downarrow}p$ it enjoys the same defining properties as $a_{\alpha f}$ and it follows that there exists an isomorphism $\gamma\colon a_{\alpha f} \cong a_\alpha f$ in the slice over $E \times B$. Horizontal composition of that isomorphism with $\psi\bar{r}$ we obtain the following middle four interchange square

$$\begin{array}{ccc} q_0\bar{r}a_{\alpha f} & \overset{q_0\bar{r}\gamma}{\cong} & q_0\bar{r}a_\alpha f \\ {\scriptstyle\psi\bar{r}a_{\alpha f}}\Downarrow & & \Downarrow{\scriptstyle\psi\bar{r}a_\alpha f} \\ q_1\bar{r}a_{\alpha f} & \underset{q_1\bar{r}\gamma}{\cong} & q_1\bar{r}a_\alpha f \end{array}$$

whose verticals are, by definition, the 2-cells $\chi_{\alpha f}$ and $\chi_\alpha f$ respectively. However, we have shown that $\chi_{\alpha f}$ is cartesian for $p$ so it follows that the 2-cell $\chi_\alpha f$ is also cartesian for $p$. Finally the comment at the end of Observation 4.1.6 applies to complete our proof. □

## References

[1] C. Barwick and C. Schommer-Pries. On the unicity of the homotopy theory of higher categories. `https://arxiv.org/abs/1112.0040`, 2011.

[2] T. Beke. Sheafifiable homotopy model categories. *Mathematical Proceedings of the Cambridge Philosophical Society*, 129(03):447–475, 2001.

[3] K. S. Brown. Abstract homotopy theory and generalized sheaf cohomology. *Transactions of the American Mathematical Society*, 186:419–458, 1973.

[4] D. Dugger. Replacing model categories with simplicial ones. *Transactions of the American Mathematical Society*, 353(12):5003–5028, 2001.

[5] D. Dugger and D. I. Spivak. Mapping spaces in quasi-categories. *Algebr. Geom. Topol.*, 11(1):263–325, 2011.

[6] W. G. Dwyer, P. S. Hirschhorn, D. M. Kan, and J. H. Smith. *Homotopy limit functors on model categories and homotopical categories*, volume 113 of *Mathematical Surveys and Monographs*. American Mathematical Society, 2004.

[7] P. Hirschhorn. *Model Categories and their Localizations*. Number 99 in Mathematical Surveys and Monographs. American Mathematical Society, 2003.

[8] A. Joyal. Quasi-categories and Kan complexes. *Journal of Pure and Applied Algebra*, 175:207–222, 2002.

[9] A. Joyal. *The theory of quasi-categories and its applications*. Quadern 45, Vol. II, Centre de Recerca Matemàtica Barcelona. 2008.

[10] A. Joyal and M. Tierney. Quasi-categories vs Segal spaces. In A. D. et al, editor, *Categories in Algebra, Geometry and Physics*, volume 431 of *Contemporary Mathematics*, pages 277–326. American Mathematical Society, 2007.
[11] Z. L. Low. The homotopy bicategory of $(\infty, 1)$-categories. http://arxiv.org/abs/1310.0381, 2013.
[12] J. Lurie. *Higher Topos Theory*, volume 170 of *Annals of Mathematical Studies*. Princeton University Press, Princeton, New Jersey, 2009.
[13] S. Mac Lane. *Categories for the Working Mathematician*. Springer-Verlag, New York, 1971.
[14] C. Rezk. A model for the homotopy theory of homotopy theory. *Transactions of the American Mathematical Society*, 353(3):973–1007, 2001.
[15] C. Rezk. A cartesian presentation of weak $n$-categories. *Geometry and Topology*, 2010.
[16] E. Riehl. *Categorical homotopy theory*, volume 24 of *New Mathematical Monographs*. Cambridge University Press, 2014.
[17] E. Riehl and D. Verity. The theory and practice of Reedy categories. *Theory and Applications of Categories*, 29(9):256–301, 2014. arXiv:1304.6871 [math.CT].
[18] E. Riehl and D. Verity. The 2-category theory of quasi-categories. *Adv. Math.*, 280:549–642, 2015.
[19] E. Riehl and D. Verity. Completeness results for quasi-categories of algebras, homotopy limits, and related general constructions. *Homol. Homotopy Appl.*, 17(1):1–33, 2015.
[20] E. Riehl and D. Verity. Kan extensions and the calculus of modules for $\infty$-categories. arXiv:1507.01460, 2015.
[21] E. Riehl and D. Verity. Homotopy coherent adjunctions and the formal theory of monads. *Adv. Math.*, 286:802–888, 2016.
[22] R. Street. Elementary cosmoi. I. In *Category Seminar (Proc. Sem., Sydney, 1972/1973)*, pages 134–180. Lecture Notes in Math., Vol. 420. Springer, Berlin, 1974.
[23] R. Street. Fibrations and Yoneda's lemma in a 2-category. In *Category Seminar (Proc. Sem., Sydney, 1972/1973)*, pages 104–133. Lecture Notes in Math., Vol. 420. Springer, Berlin, 1974.
[24] R. H. Street. Fibrations in bicategories. *Cahiers de Topologie et Géom. Diff.*, XXI-2:111–159, 1980.
[25] R. H. Street. Correction to "fibrations in bicategories". *Cahiers topologie et géométrie différentielle catégoriques*, 28(1):53–56, 1987.
[26] B. Toën. Vers une axiomatisation de la théorie des catégories supérieures. *K-Theory*, 34(3):233–263, 2005.
[27] D. Verity. Weak complicial sets I, basic homotopy theory. *Advances in Mathematics*, 219:1081–1149, September 2008.

Department of Mathematics, Johns Hopkins University, Baltimore, MD 21218, USA
*E-mail address*: eriehl@math.jhu.edu

Centre of Australian Category Theory, Macquarie University, NSW 2109, Australia
*E-mail address*: dominic.verity@mq.edu.au